\documentclass{article}

\usepackage[T1]{fontenc}
\usepackage{authblk}

\usepackage{tocloft}
\usepackage[numbers]{natbib}
\usepackage{enumerate, amsmath, amsfonts, amssymb, amsthm, dsfont, mathrsfs, wasysym, graphics, graphicx, xcolor, url, hyperref, hypcap, xargs, multicol, pdflscape, multirow, hvfloat, array, ae, aecompl, pifont, mathtools, a4wide, float, blkarray, overpic, nicefrac, stmaryrd, anyfontsize, yfonts, fontawesome}
\usepackage{pdflscape}
\usepackage{xargs, bbm, enumerate, paralist}
\usepackage[noabbrev,capitalise]{cleveref}
\usepackage[normalem]{ulem}
\usepackage{marginnote}

\hypersetup{colorlinks=true, citecolor=darkblue, linkcolor=darkblue}
\usepackage[all]{xy}
\usepackage{tikz}
\usepackage{tikz-cd}
\usetikzlibrary{trees, decorations, decorations.pathmorphing, decorations.markings, decorations.shapes, shapes, arrows, matrix, calc, fit, intersections, patterns, angles}
\graphicspath{{figures/}{figures/diagonals/}{figures/walks/}{figures/tubes/}{figures/blocks/}}
\makeatletter\def\input@path{{figures/}}\makeatother
\usepackage{caption}
\usepackage[export]{adjustbox}

\usepackage{algorithm2e}
\usepackage[noend]{algpseudocode}

\usepackage{paralist}

\usepackage{listings}
\DeclareFontEncoding{LY}{}{}
\DeclareFontSubstitution{LY}{yfrak}{m}{n}
\DeclareFontEncoding{LYG}{}{}
\DeclareFontSubstitution{LYG}{ygoth}{m}{n}
\DeclareFontFamily{LYG}{ygoth}{}
\DeclareFontShape{LYG}{ygoth}{m}{n}{<->ygoth}{}
\DeclareFontFamily{LY}{yfrak}{}
\DeclareFontShape{LY}{yfrak}{m}{n}{<->yfrak}{}
\DeclareFontFamily{LY}{ysmfrak}{}
\DeclareFontShape{LY}{ysmfrak}{m}{n}{<->ysmfrak}{}
\DeclareFontFamily{LY}{yswab}{}
\DeclareFontShape{LY}{yswab}{m}{n}{<->yswab}{}

\newtheorem{thmUniv}{Theorem}

\newtheorem{theorem}{Theorem}[section]
\newtheorem{corollary}[theorem]{Corollary}

\newtheorem{proposition}[theorem]{Proposition}
\newtheorem{lemma}[theorem]{Lemma}

\newtheorem*{theorem*}{Theorem}

\theoremstyle{definition}
\newtheorem{definition}[theorem]{Definition}
\newtheorem{example}[theorem]{Example}
\newtheorem{remark}[theorem]{Remark}
\newtheorem{question}{Question}

\newtheorem{problem}[theorem]{Problem}

\crefname{equation}{Equation}{Equations}
\crefname{problem}{Problem}{Problems}

\newcommand{\R}{\mathbb{R}} 
\newcommand{\N}{\mathbb{N}} 
\renewcommand{\P}{\mathbb{P}} 
\renewcommand{\c}[1]{{\mathcal{#1}}} 
\renewcommand{\b}[1]{{\boldsymbol{#1}}} 
\newcommand{\E}{\mathbb{E}} 
\newcommand{\B}{\mathbb{B}^2} 
\newcommand{\X}{\mathbb{X}} 

\renewcommand{\emptyset}{\varnothing} 
\renewcommand{\epsilon}{\varepsilon} 
\newcommand{\ssm}{\smallsetminus} 
\newcommand{\eqdef}{\mbox{\,\raisebox{0.2ex}{\scriptsize\ensuremath{\mathrm:}}\ensuremath{=}\,}} 
\newcommand{\simplex}{\polytope{\Delta}} 

\DeclareMathOperator{\conv}{conv} 
\DeclareMathOperator{\cone}{cone} 
\DeclareMathOperator{\Var}{Var} 
\DeclareMathOperator{\vis}{vis_{\varepsilon}} 

\newcommand{\ie}{\textit{i.e.},~} 
\newcommand{\eg}{\textit{e.g.},~} 
\definecolor{darkblue}{rgb}{0,0,0.7} 
\definecolor{green}{RGB}{57,181,74} 
\definecolor{violet}{RGB}{147,39,143} 
\newcommand{\darkblue}{\color{darkblue}} 
\newcommand{\defn}[1]{\textsl{\darkblue #1}} 

\usepackage{todonotes}

\makeatletter
\def\part{\@startsection{part}{1}%
\z@{.7\linespacing\@plus\linespacing}{.8\linespacing}%
{\LARGE\sffamily\centering}}
\makeatother

\newcommand{\polytope}[1]{\mathsf{#1}}
\newcommand{\polytopeP}{\mathsf{P}}
\newcommand{\polytopeB}{\mathsf{B}}
\newcommand{\polytopeD}{\mathsf{D}}
\newcommand{\pol}[1][P]{\polytope{#1}}
\newcommand{\polytopeQ}{\mathsf{Q}}
\newcommand{\polytopeF}{\mathsf{F}}
\newcommand{\FBody}{\mathsf{F}_{\varepsilon}}
\newcommand{\polytopeC}{\mathsf{C}}
\newcommand{\Lop}[1][d]{\polytope{Lop}_{#1}}

\newcommand{\CrossPol}[1][d]{\Diamond_{#1}}
\newcommand\cyc{\polytope{Cyc}}

\newcommand\Ass[1][n]{\polytope{Asso}_{#1}}
\newcommand{\permuto}[1][n]{\polytope{\Pi}_{#1}}
\newcommandx{\HypSimpl}[2][1=d,2=k]{\polytope{\Delta}(#1,#2)}
\newcommand{\HypSimplTwo}[1][d]{\HypSimpl[#1][2]}
\newcommandx{\MPP}[2][1=\polytopeP,2=\b c]{\polytope{M}_{#2}(#1)}
\newcommandx{\MPPHypSimpl}[2][1=d,2=k]{\polytope{M}(#1,#2)}

\newcommandx{\FibPol}[3][1=\polytope{P},2=\polytope{Q},3=\pi]{\Sigma_{#3}(#1,#2)}

\newcommandx{\PolProj}[3][1=\polytope{P},2=\polytope{Q},3=\pi]{#3~:~#1\to #2}

\newcommand{\inner}[1]{\left<#1\right>}

\title{Unimodality of the number of paths per length on polytopes:\\[0.2em]\smaller{}Examples, counterexamples, and a central limit theorem}

\author{Martina Juhnke, Germain Poullot}

\AtEndDocument{\bigskip{\footnotesize%
  \textsc{Universität Osnabrück, Germany} \par  
  \textit{E-mail}, M.~Juhnke: \texttt{martina.juhnke@uni-osnabrueck.de} ~;~~~
G.~Poullot: \texttt{germain.poullot@uni-osnabrueck.de} \par
}}

\date{}

\begin{document}

\maketitle

\begin{abstract}
Because of its importance in combinatorics and optimization we study the full distribution of the lengths of monotone paths of a convex polytope. De Loera had conjectured that the number of such paths counted by length is always a unimodal sequence.
We confirm this for several classes of polytopes but construct counterexamples disproving the conjecture in general.
Nevertheless, we show that for random polytopes with vertices uniformly distributed on the sphere, the length of coherent paths satisfies a central limit theorem.



\end{abstract}

\tableofcontents

\section*{Acknowledgments}
First, our deepest thanks go to Jes{\'u}s De Loera without whom this project would not have started: he asked the second author in December 2022 a question which became \Cref{prb:Unimodality}.
Furthermore, the authors want to express their gratitude to Matthias Reitzner and Christoph Th\"ale for their help with the literature on random polytopes, and to Gilles Bonnet for his clever remarks on coherent paths.
In addition, we thank Kinga Nagy for her joyful technical support on the most convoluted integral computations of this paper.
Last but not least, we are grateful to Alex Black for being an inexhaustible source of comments around the applications of the simplex method.



\section{Introduction}







For a polytope $\pol\subset \R^d$ and a direction $\b c\in \R^d$, we wonder about the \defn{monotone paths}:
paths in the graph of $\pol$ along which the scalar product with $\b c$ increases.
How many monotone paths are there?
Of which lengths?
Are there more short/long paths, or more paths of almost average length?

Not only are these questions natural to ask, but in addition, they are also of prime importance in several contexts.
First, to solve a linear program using the famous simplex method introduced by Dantzig in 1947 (see \cite{Dantzig1963-LinearProgrammingAndExtensions}), one traverses such a path.
Hence, the length of monotone paths is a key measure of the complexity of the simplex method.
Starting with Klee and Minty's \cite{KleeMinty-SimplexAlgorithm} seminal example of a family of polytopes with $m$ facets and monotone paths of length exponential in $m$, monotone paths on polytopes have received a lot of attention, \eg \cite{AthanasiadisDeLoeraZhang-EnumerativeProblems} studies the number of monotone paths on polytopes, whereas \cite{AthanasiadisEdelmanReiner-MonotonePathsOnPolytopes} investigates the connectivity of the graph of paths (where a path can be ``flipped'' into another by switching it around a 2-face).
In addition, \cite{BlanchardDeLoeraLouveaux-LengthMonotonePath} presents some results on the extremal lengths of monotone paths for specific classes of polytopes, and \cite{AthanasiadisSantos-MonotonePathsOnZonotopes} addresses the case of zonotopes.

On a combinatorial side, if the graph of $\pol$ (directed along $\b c$) is a lattice, then the monotone paths are the maximal chains in this lattice.
Nelson \cite{Nelson2017-MaximalChainsTamariLattice} unraveled monotone paths on the associahedron (maximal chains in the Tamari lattice).
This is extended to graph associahedra in \cite{DahlbergFishel2024-MaximalChainsGraphAssociahedra}.
Monotone paths on the permutahedron are known as \emph{sorting networks} \cite{AngelHolroydRomikVirag-RandomSortingNetworks,Dauvergne-ArchimedianLimitRandomSortingNetwork}.


The simplex method follows a monotone path according to a \defn{pivot rule}, where, at each vertex, this rule determines which ($\b c$-improving) adjacent vertex to choose next in the path.
As short paths are desirable, clever pivot rules were proposed to avoid ``whirling too much'' around $\pol$.
One such instance is the \defn{shadow vertex rule}, which first projects $\pol$ onto a chosen $2$-dimensional plane, and then takes one of the two boundary paths of the resulting polygon as its monotone path.
Following this idea, a monotone path is called \defn{coherent} if it is selected by the shadow vertex rule for some plane of projection. 
Billera and Sturmfels \cite{BilleraSturmfels-FiberPolytope} generalized coherent paths to coherent subdivisions via fiber polytopes.
Based on this work, coherent paths (and monotone path polytopes) were studied for simplices and cubes \cite{BilleraSturmfels-FiberPolytope}, for cyclic polytopes \cite{Athanasiadis_2000}, for $S$-hypersimplices \cite{ManeckeSanyalSo2019shypersimplices}, for cross-polytopes \cite{BlackDeLoera2021monotone} and for (usual) hypersimplices \cite{poullot2024verticesmonotonepathpolytopes}.

Borgwardt analyzed the shadow vertex rule, computing the expectation of the length of a coherent path for several classes of random polytopes:
He left open the computation of the higher probabilistic moments \cite[Ch. 0.12, Question 8]{Borgwardt1987-SimplexMethod}. 
Our starting point is the following question:

\begin{question}\label{prb:Unimodality}
Let \defn{$N_\ell$} (respectively \defn{$N^{\text{coh}}_\ell$}) be the number of monotone (resp. coherent) paths of length $\ell$ on a polytope $\pol$ for a direction $\b c$.
Are the sequences $(N_\ell)_\ell$ (resp. $(N_\ell^{\text{coh}})_\ell$) unimodal?
\end{question}

Recall that a sequence $\mathbf{a}=(a_1, \dots, a_r)\in\N^r$ is \defn{unimodal} if there is $k\in [r]$ such that $a_i \leq a_{i+1}$ for $i < k$, and $a_i\geq a_{i+1}$ for $i \geq k$. 
Jes{\'u}s de Loera (personal communication) has conjectured that \Cref{prb:Unimodality} has an affirmative answer, confirming it for special instances.
The goal of this paper is threefold.
First, we provide more examples where the answer is ``yes'', leading to (see \Cref{sec:PositiveExamples}):
\begin{thmUniv}\label{thm:Intro0}
The sequences $(N_\ell)_\ell$ and $(N_\ell^{\text{coh}})_\ell$ are unimodal for:
\begin{compactenum}
\item[(a)] the $d$-simplex, for any generic $\b c$;
\item[(b)] the standard $d$-cube, for any generic $\b c$;
\item[(c)] the $d$-cross-polytope, for any generic $\b c$;
\item[(d)] the cyclic polytopes, for $\b c = (1, 0, \dots, 0)$;
\item[(e)] the $S$-hypersimplex, for $\b c = (1, \dots, 1)$;
\item[(f)] the prism $\pol\times[0, 1]$, for $(\b c, 1)$, if the corresponding sequence for $\pol$ and $\b c$ is log-concave.
\end{compactenum}
\end{thmUniv}

Second, we show that the answer to \Cref{prb:Unimodality} is ``no'' in general (\Cref{sec:NegativeExamples}) by providing specific counterexamples for several classes of polytopes. 
In summary, these results are:

\begin{thmUniv}\label{thm:IntroA}
For the following classes, there exist polytopes and $\b c$ such that $(N_\ell)_\ell$ is not unimodal:
\begin{compactenum}
\item[(a)] $d$-dimensional polytopes, combinatorially isomorphic to the $d$-cube for $d\geq 3$ (\Cref{cor:LopsidedCube});
\item[(b)] $3$-dimensional simplicial polytopes (\Cref{thm:counterexampleSimplicial});
\item[(c)] $5$-dimensional generalized permutahedra (\Cref{thm:CounterexampleGenPer});
\item[(d)] $5$-dimensional $0/1$-polytopes, not combinatorially isomorphic to the $5$-cube (\Cref{thm:Counterexample01}).
\end{compactenum}
Moreover, for the classes $(a)$ and $(c)$, the sequence $(N_\ell^{\text{coh}})_\ell$ is not unimodal either.
\end{thmUniv}

Third, our aim is to demonstrate that, notwithstanding the negative results of \Cref{prb:Unimodality}, random polytopes admit a somewhat positive answer.
We prove that the length of a coherent path (for random polytopes with vertices on the unit sphere) admits a central limit theorem (\Cref{sec:RandomCase}):

\begin{thmUniv}[\Cref{thm:CLTmain}]\label{thm:IntroB}
Let $d\geq 4$ and let $\b c, \b \omega\in \R^d$ be linearly independent vectors (possibly randomly chosen).
Let $Z_1, \dots, Z_n$ be independent random points, chosen uniformly on the sphere $\mathbb{S}^{d-1}\subseteq\R^d$, and let $\pol_n = \conv(Z_1, \dots, Z_n)$.
Then, the length $L_n$ of the coherent $\b c$-monotone path captured by~$\b\omega$ on~$\pol_n$ obeys a central limit theorem (here, the convergence is in distribution):
\[
\frac{L_n - \E (L_n)}{\sqrt{\Var (L_n)}} \xrightarrow[n\to+\infty]{} U, ~~~~ \text{where } U\sim \c N(0, 1) \text{ is standard normally distributed}
\]
Moreover, $\E(L_n) \sim c\, n^{\frac{1}{d-1}}$, 
and: $c'\, n^{\frac{1}{d-1} - a}\leq \Var(L_n) \leq c''\, (\log n)^{3 - \frac{1}{d-1}} \, n^{\frac{1}{d-1}}$ for any $a > 0$.
\end{thmUniv}
 We refer to \Cref{sec:RandomCase} for unexplained notation and all probabilistic background. 
We want to remark that, \eg via flips over $2$-faces (see \cite{BilleraKapranovSturmfels-CellularStrings,AthanasiadisDeLoeraZhang-EnumerativeProblems}), one can usually ``slightly modify'' a path to obtain a path of similar length. One therefore expect that there are only relatively few paths of large/small length, whereas the length of most of the paths is close to the average length. Indeed, this intuition is confirmed by \Cref{thm:IntroB}.


\section{Preliminaries}

\subsection*{On monotone paths and coherent paths on polytopes}\label{subsec:monCoh}
In this subsection, we provide the necessary background on polytopes and their monotone/coherent paths. We recommend \cite{Gruenbaum-ConvexPolytopes} and \cite{Ziegler-LecturesOnPolytopes} for more details on these topics. 

A \defn{polytope} $\pol\subseteq \R^d$ is the convex hull of finitely many points, or equivalently, the bounded intersection of finitely many half-spaces.
For $\b c\in \R^d$, we let $\pol^{\b c} \eqdef \{\b x\in \pol ~;~ \inner{\b x, \b c} = \max_{\b y\in \pol}\inner{\b y, \b c}\}$.
A subset $\polytopeF\subseteq\pol$ is a \defn{face} of $\pol$ if there exists $\b c\in \R^d$ such that $\polytopeF = \pol^{\b c}$.
The \defn{dimension} of a face $\polytopeF$ is its affine dimension, \ie the dimension of the smallest affine subspace of $\R^d$ containing $\polytopeF$.
By convention, $\emptyset$ is a face of $\pol$ of dimension $-1$.
\defn{Vertices} and \defn{edges} are $0$- and $1$-dimensional faces of $\pol$, respectively. 
The vertices and edges of $\pol$ naturally form a graph $G_{\pol}$. 
A direction $\b c\in \R^d$ is \defn{generic} with respect to $\pol$ if $\inner{\b u, \b c} \ne \inner{\b v, \b c}$ for every edge $[\b u, \b v]$ of $\pol$.
Here, \defn{$[\b u,\b v]$} denotes the line segment between $\b u$ and $\b v$. 
If $\b c\in \R^d$ is generic, one can turn $G_{\pol}$ into a \defn{directed graph $G_{\pol, \b c}$} by directing an edge $\b u \b v$ of $G_{\pol}$ from $\b u$ to $\b v$ if and only if $\inner{\b u, \b c} < \inner{\b v, \b c}$.
By definition, $G_{\pol, \b c}$ has a unique source \defn{$\b v_{\min}$} and a unique sink \defn{$\b v_{\max}$}, given by the vertices of $\pol$ that minimizes and maximizes $\inner{\b x, \b c}$ for $\b x\in \pol$, respectively.

\begin{figure}[!h]
  \centering
  \includegraphics[width=0.95\linewidth]{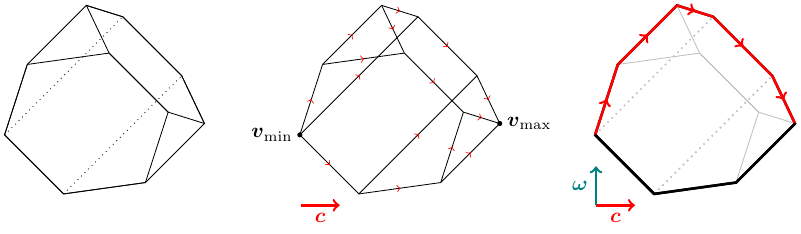}
  \caption[Monotone paths versus coherent paths on a polytope]{(Left) A $3$-dimensional polytope $\pol$.
  (Middle) The graph $G_{\pol, \b c}$ having 7 monotone paths (3 of length 3, and 2 of length 4, and 2 of length 5).
  (Right) The projection of $\pol$ onto the plane spanned by $\b c$ and $\b\omega$ with, in \textcolor{red}{red}, the path formed by the upper faces of this projection (another choice of $\b\omega$ would give rise to another coherent path).}
  \label{fig:ExamplePolytopeWithMonotonePaths}
\end{figure}

\begin{definition}
For a polytope $\pol\subseteq\R^d$ and a direction $\b c\in \R^d$, a \defn{$\b c$-monotone path} is a directed path in $G_{\pol, \b c}$ from $\b v_{\min}$ to $\b v_{\max}$, see \Cref{fig:ExamplePolytopeWithMonotonePaths} (left and middle).
The \defn{length} of a $\b c$-monotone path is its number of edges, \ie its number of vertices minus 1.
We denote by \defn{$N_\ell(\pol,\b c)$} the number of $\b c$-monotone paths of $\pol$ of length $\ell$. 
\end{definition}

If $\pol$ and $\b c$ are clear from the context, we will often simply write \defn{$N_\ell$} instead.
Similarly, we will often write ``monotone path'' instead of $\b c$-monotone path.

\begin{remark}
Balinski's theorem, see \cite[Section 3.5]{Ziegler-LecturesOnPolytopes}, ensures that the graph of a $d$-dimensional polytope is $d$-connected. 
It follows that $\sum_\ell N_\ell \geq d$ since by Menger's theorem there exist at least $d$ (internally disjoint) $\b c$-monotone paths on $\pol$ (for generic $\b c$). In \cite[Proposition 4.5]{AthanasiadisDeLoeraZhang-EnumerativeProblems}, this lower bound was improved to $\sum_\ell N_\ell \geq \left\lceil\frac{d\cdot f_0}{2}\right\rceil - f_0 + 2$, where $f_0$ is the number of vertices of $\pol$.
\end{remark}

\begin{definition}\label{def:CoherentPaths}
Let $\pol\subseteq\R^d$ be a polytope and $\b c, \b \omega \in \R^d$ be linearly independent directions. Let \defn{$\pol_{\b c, \b \omega}$} be the polygon that is the orthogonal projection of $\pol$ onto the plane spanned by $\b c$ and $\b \omega$:
\[
 \pol_{\b c, \b \omega} \eqdef \pi_{\b c, \b\omega}(\pol)\eqdef \Bigl\{\bigl(\inner{\b x, \b c}, \inner{\b x, \b \omega}\bigr) ~;~ \b x\in \pol\Bigr\}
\]
(see \Cref{fig:ExamplePolytopeWithMonotonePaths}, right). 
An edge $\polytope{F}$ of $\pol_{\b c, \b \omega}$ is an \defn{upper edge} if it has an outer normal vector with positive second coordinate, \ie if $(x_1, x_2) + (0, \varepsilon) \notin \pol_{\b c, \b \omega}$ for all $(x_1, x_2) \in \polytope{F}$ and all $\varepsilon > 0$.
The upper edges of $\pol_{\b c,\b \omega}$ form a path on the boundary of $\pol_{\b c,\b \omega}$, referred to as the \defn{upper path}. The remaining edges form the so-called \defn{lower path} of $\pol_{\b c,\b \omega}$. 
A $\b c$-monotone path $\c L$ on $\pol$ is \defn{coherent} if there exists $\b \omega \in \R^d$ such that the projected path $\pi_{\b c, \b\omega}(\c L)$ is the  upper path of $\pol_{\b c, \b\omega}$. In this case, such a $\b\omega$ is said to \defn{capture} the coherent path $\c L$.
The \defn{length} of a coherent path is its number of edges.
We use \defn{$N^{\text{coh}}_\ell(\pol, \b c)$} to denote the number of coherent paths of length $\ell$ on $\pol$ (with respect to $\b c$). 
We write \defn{$N^{\text{coh}}_\ell$} if $\pol$ and $\b c$ are clear from the context.
\end{definition}

\begin{remark}
On the one hand, we have $N^{\text{coh}}_\ell \leq N_\ell$. 
On the other hand, since the lower path of  $\pol_{\b c, \b \omega}$ is the upper path of $\pol_{\b c,-\b\omega}$ and since these two paths are internally disjoint, we have $\sum_\ell N^{\text{coh}}_\ell \geq 2$.
Using the theory of monotone path polytopes, this lower bound can be improved (\eg to $\sum_\ell N_\ell^{\text{coh}} \geq d$). However, this is outside the scope of this article.
\end{remark}

\subsection*{On tools around unimodality}
In this subsection, we gather several tools to prove  unimodality.
We do not intend to provide a handbook on unimodality nor log-concavity and refer the interested reader to \eg \cite{Braenden}.


A sequence  $\mathbf{a}=(a_1, \dots, a_r)\in \mathbb{N}^r$ is \defn{log-concave} if $a_{i-1}a_{i+1} \leq a_i^2$ for all $2\leq i\leq r-1$. The sequence $\mathbf{a}$ is \defn{symmetric} if $a_{r-i} = a_i$ for all $i\in [r]$.
Any (not necessarily unique) index $k$ is a \defn{mode} of $\mathbf{a}$ if $a_k = \max(a_\ell~;~1\leq \ell\leq r)$. 
It is well-known that any log-concave sequence is also unimodal. 
An example of a sequence that is symmetric and log-concave (and, hence,  unimodal) is provided by the binomial coefficients $(\binom{n}{i} ~;~ i\in [n])$. 
We collect some helpful results that we will use to prove log-concavity or unimodality, omitting their proofs since they are well-known and easy.

\begin{lemma}
\label{lem:OnLogConcaveSequences}
Let $\mathbf{a}^{(1)}=(a^{(1)}_i)_{i\in[r]}, \ldots, \mathbf{a}^{(m)}=(a^{(m)}_i)_{i\in[r]}$ be sequences of non-negative integers.
\begin{compactenum}
\item[(1)] If $\mathbf{a}^{(1)}, \ldots, \mathbf{a}^{(m)}$ are log-concave, then so is their product $(a^{(1)}_i\cdots a^{(m)}_i)_{i\in[r]}$. 
\item[(2)] If $\mathbf{a}$ is unimodal and symmetric, then $\left\lfloor\frac{r+1}{2}\right\rfloor$ and $ \left\lceil\frac{r+1}{2}\right\rceil$ are among its modes.
\item[(3)] If $\mathbf{a}^{(1)}, \ldots, \mathbf{a}^{(m)}$ are unimodal sequences with mode $k$, then  $\left(\sum_{j=1}^m a^{(j)}_i\right)_{i\in[r]}$ is unimodal and has $k$ as one of its modes.
\item[(4)] If $\mathbf{a}^{(1)}$ and $\mathbf{a}^{(2)}$ are unimodal with respective modes $k$ and $k+1$, then their sum $(a_i^{(1)}+a_i^{(2)})_{i\in [r]}$ is unimodal and at least one of $k$ or $k+1$ is a mode.
\end{compactenum}
\end{lemma}

\subsection*{On Landau notation}

Throughout this paper, we use the usual Landau notation (see \eg \cite[Chapter IV.3]{AE-Landau} for more details).
For $a \in \R\cup \{\pm \infty\}$, and for functions $f, g : \R \to \R\cup\{\pm\infty\}$, we write
\begin{compactitem}
\item \defn{${f(x) = O\bigl(g(x)\bigr)}$} if and only if there exists $k\in \R$ such that $|f(x)| \leq k \cdot |g(x)|$ for all $x$ in a neighborhood of $a$;
\item \defn{${f(x) = o\bigl(g(x)\bigr)}$} if and only if $\frac{f(x)}{g(x)} \to 0$ if $x\to a$;
\item \defn{${f(x) \sim g(x)}$} if and only if $\frac{f(x)}{g(x)} \to 1$ if $x\to a$.
\end{compactitem}

\section{Positive examples}\label{sec:PositiveExamples}
The aim of this section is to provide examples for which \Cref{prb:Unimodality} has a positive answer, \ie for which $(N_\ell)_\ell$ and $(N^{\text{coh}}_\ell)_\ell$ are unimodal sequences.
In addition to previously known examples, this also includes some new ones. 
First, we recall the definitions of several  polytopes. 

\begin{center}
\begin{tabular}{l|c|l}
polytope & notation & definition \\ \hline
$d$-simplex & $\simplex_d$ & convex hull of $d+1$ affinely independent points \\
standard $d$-cube & $[0, 1]^d$ & $[0, 1]^d$ \\
$d$-cross-polytope & $\CrossPol$ & $\conv(\pm\b e_i ~;~ 1\leq i\leq d)$ with $\b e_i$ the $i$-th unit vector in $\R^d$ \\
cyclic polytope & $\cyc_d(\b t)$ & $\conv\bigl((t_i, t_i^2, \dots, t_i^d) ~;~ 1\leq i\leq n\bigr)$ with $\b t = (t_1, \dots, t_n)\in \R^n$ \\
$S$-hypersimplex & $\simplex_d(S)$ & $\conv(\b x\in \{0, 1\}^d ~;~ \sum_i x_i \in S)$ with $S \subseteq [d]$ \\
\end{tabular}
\end{center}

For the above examples, the total number of monotone paths $\sum_\ell N_\ell$ can be found in the literature or is not hard to deduce. 
In each case, monotone paths are associated with a combinatorial object: we refine this count to deduce the number of monotone paths of length $\ell$. 
The next table lists these results (with the direction $\b c$).
The first column provides a reference to an article and/or a remark/example below, where unexplained notions are defined and unimodality is studied. 

\begin{center}
\begin{tabular}{l|c|c|c|c}
reference& $\pol$ & $\b c$ & $\sum_\ell N_\ell$ & $N_\ell$ \\ \hline
\Cref{exm:CompleteGraph} & any $d$-simplex &  generic & $2^{d-1}$ & $\binom{d-1}{\ell-1}$ \\ 
\Cref{exm:GradedGraph} & $[0, 1]^d$ &  generic & $d!$ & $d!$ iff $\ell = d$ \\ 
\cite{BlackDeLoera2021monotone}, \Cref{rmk:CrossPolytopeMonotone} & $\CrossPol$ & generic & $\frac{1}{3}(2^{2d-1} - 2)$ & $2\sum_{k=0}^{d-2} \binom{2k}{\ell-2}$ \\
\Cref{exm:CompleteGraph} & $\cyc_d(\b t)$,\, $d\geq 4$ & $(1, 0, \dots, 0)$ & $2^{n-2}$ & $\binom{n-2}{\ell-1}$ \\
\cite{ManeckeSanyalSo2019shypersimplices}, \Cref{rmk:SHypersimplexPaths} & $\simplex_d(S)$,\, $|S| = r$ & $(1, 1, \dots, 1)$ & $\binom{d}{\Tilde{s}_1, \Tilde{s}_2, \dots, \Tilde{s}_r}$ & $\binom{d}{\Tilde{s}_1, \Tilde{s}_2, \dots, \Tilde{s}_r}$ iff $\ell = r$ \\
\end{tabular}
\end{center}

The next table provides  analogous information for coherent paths.

\begin{center}
\begin{tabular}{l|c|c|c|c}
reference & $\pol$ & $\b c$ & $\sum_\ell N^{\text{coh}}_\ell$ & $N^{\text{coh}}_\ell$ \\ \hline
\cite{BilleraSturmfels-FiberPolytope}, \Cref{exm:CompleteGraph} & any $d$-simplex &  generic & $2^{d-1}$ & $\binom{d-1}{\ell-1}$ \\ 
\cite{BilleraSturmfels-FiberPolytope}, \Cref{exm:GradedGraph} & $[0, 1]^d$ &  generic & $d!$ & $d!$ iff $\ell = d$ \\ 
\cite{BlackDeLoera2021monotone}, \Cref{rmk:CrossPolytopeCoherent} & $\CrossPol$ &  generic & $3^{d-1} - 1$ & $\binom{d-1}{\ell-1} 2^{\ell-1}$ \\
\cite{Athanasiadis_2000}, \Cref{rmk:CycCoherent} & $\cyc_d(\b t)$,\, $d\geq 4$ & $(1, 0, \dots, 0)$ & $2\sum_{j=0}^{d-2}\binom{n-3}{j}$ & \Cref{rmk:CycCoherent} \\
\cite{ManeckeSanyalSo2019shypersimplices}, \Cref{rmk:SHypersimplexPaths} & $\simplex_d(S)$,\, $|S| = r$ & $(1, 1, \dots, 1)$ & $\binom{d}{\Tilde{s}_1, \Tilde{s}_2, \dots, \Tilde{s}_r}$ & $\binom{d}{\Tilde{s}_1, \Tilde{s}_2, \dots, \Tilde{s}_r}$ iff $\ell = r$ \\
\end{tabular}
\end{center}

\begin{example}\label{exm:CompleteGraph}
If the graph of a polytope $\pol$ is the complete graph on $n$ vertices (as  \eg for  simplices, cyclic polytopes for $d\geq 4$, and more generally, neighborly polytopes), then its directed graph $G_{\pol, \b c}$ yields an acyclic orientation of the underlying complete graph if $\b c $ is generic.
Any monotone path of length $\ell$ hence corresponds to an $(\ell+1)$-element subset of the vertices of $\pol$ containing $\b v_{\min}$ and $\b v_{\max}$.
This implies that $N_\ell = \binom{n-2}{\ell-1}$, which clearly is a unimodal sequence, and thus provides a partial positive answer to \Cref{prb:Unimodality}.
For the simplex, it follows from the paragraph preceding Example 5.4 in \cite{BilleraSturmfels-FiberPolytope} that every monotone path is coherent.
\end{example}

\begin{remark}\label{rmk:CrossPolytopeMonotone}
Let $\b c$ be a generic direction for the $d$-cross-polytope (which means $c_i \ne c_j$ for $i\ne j$). 
By \cite[page 11]{BlackDeLoera2021monotone} $\b c$-monotone paths of length $\ell$ on $\CrossPol$ are in bijection with subsets $X\subseteq\{-(d-1), \dots, -1, +1, \dots, d-1\}$ with $|X|=\ell-1$ such that if $\{-k,k\}\subseteq X$, then there exists $j\in X$ with $-k < j < +k$.
We now count such subsets for fixed $\ell$. 
First, note that there are $\binom{2(d-1-k)}{\ell-3}$ subsets of  size $\ell-1$  that contain $\{-k,+k\}$ but no value in between. Since the total number of subsets of size $\ell-1$ of $\{-(d-1), \dots, -1, +1, \dots, d-1\}$ equals $\binom{2(d-1)}{\ell-1}$, it follows that 
$N_{d, \ell} = \binom{2(d-1)}{\ell-1} - \sum_{k=1}^{d-1} \binom{2(d-1-k)}{\ell-3}$, where we write $N_{d,\ell}$ instead of $N_\ell$ to account for the dimension. 
Using $\binom{2d}{\ell-1} = \binom{2(d-1)}{\ell-1} + 2 \binom{2(d-1)}{\ell-2} + \binom{2(d-1)}{\ell-3}$, we get $N_{d+1, \ell} = N_{d, \ell} + 2\binom{2(d-1)}{\ell-2}$, yielding
$$N_{d, \ell} = 2\sum_{k=0}^{d-2} \binom{2k}{\ell-2}.$$

Using \Cref{rmk:CrossPolytopeMonotone} we can answer \Cref{prb:Unimodality} in the affirmative for $d$-cross-polytopes.

\begin{lemma}
If $d\geq 3$, the sequence $(N_{d,\ell})_\ell$ is unimodal with  mode $d$.
If $d \geq 4$, the mode is unique.
\end{lemma}
\begin{proof}
We prove the claim by induction. For $d=3$, $(N_{3, \ell})_\ell = (0, 4, 4, 2)$ is unimodal with modes $2$ and $3$.
Suppose, by induction, that $(N_{d, \ell})_\ell$ is unimodal with mode $d$. As $\left(2\binom{2(d-1)}{\ell-2}\right)_\ell$ is unimodal with mode $d+1$, \Cref{lem:OnLogConcaveSequences} (3) and (4) imply that $N_{d+1, \ell}$ is unimodal with mode $d$ or $d+1$.

We show $N_{d+1, d} < N_{d+1, d+1}$.
By the hockey-stick identity $\sum_{k\leq 2(d-1)} \binom{k}{\ell-2} = \binom{2d-1}{\ell-1}$, we get
\begin{equation}\label{eq:unimodal}
\sum_{\substack{k \leq 2(d-1)\\ k\text{ odd}}} \binom{k}{d-2} + \frac{1}{2}N_{d+1, d} = \binom{2d-1}{d-1} = \binom{2d-1}{d} = \frac{1}{2}N_{d+1, d+1} + \sum_{\substack{k \leq 2(d-1)\\ k\text{ odd}}} \binom{k}{d-1}.
\end{equation}
Since  $\binom{k}{d-2} \geq \binom{k}{d-1}$ if $k$ is odd and $k\leq 2(d-1)$ and since this inequality is strict if $k \ne 2d-3$, we conclude
 $\sum_{\substack{k \leq 2(d-1)\\ k\text{ odd}}} \binom{k}{d-2} > \sum_{\substack{k \leq 2(d-1)\\ k\text{ odd}}} \binom{k}{d-1}$ and $N_{d+1, d} < N_{d+1, d+1}$ by \eqref{eq:unimodal}.
\end{proof}
\end{remark}

\begin{remark}\label{rmk:SHypersimplexPaths}
Let $S = \{s_1 < \dots < s_r\}$. By \cite[Corollary 4.1]{ManeckeSanyalSo2019shypersimplices} the $(1, \dots, 1)$-monotone paths on the $S$-hypersimplex are in bijection with chains $A_1\subset A_2\subset\dots\subset A_r\subseteq [d]$ with $|A_i| = s_i$ for $i\in [r]$ and there are $\binom{d}{\Tilde{s}_1, \Tilde{s}_2, \dots, \Tilde{s}_r}$ such chains, where $\Tilde{s}_1 = s_1$ and $\Tilde{s}_i = s_i - s_{i-1}$ for $i > 1$. As this also implies that  all monotone paths have length $r = |S|$, the sequence $(N_\ell)_\ell$ is unimodal. 
Since in \cite{ManeckeSanyalSo2019shypersimplices}, it is also shown that all monotone paths are coherent, $(N^{\text{coh}}_\ell)_\ell$ is also unimodal.
This provides another class of polytopes for which \Cref{prb:Unimodality} has a positive answer.
The example of $S$-hypersimplices generalizes to (rank selected) matroid polytopes \cite[Section 8]{BlackSanyal-FlagPolymatroids} and non-generic directions.
However, the number of (coherent) paths per length does not seem easy to derive.
\end{remark}

\begin{example}\label{exm:GradedGraph}
If we let $\pol$ be  the $d$-cube  and $\b c$ be  a generic direction, then all directed paths in $G_{\pol, \b c}$ from $\b v_{\min}$ to $\b v_{\max}$ have the same length. In particular, both sequences $(N_\ell)_\ell$ and $(N^{\text{coh}}_\ell)_\ell$ are unimodal in this case as they contain only one term. This also happens for the permutahedron $\permuto = \conv\bigl((\sigma(1), \sigma(2), \dots, \sigma(n) ~;~ \sigma\in S_n\bigr)$ with $\b c = (1, 2, \dots, n)$, for all Coxeter permutahedra, and for zonotopes for generic directions. 
For the cube, it turns out (see \cite[Example 5.4]{BilleraSturmfels-FiberPolytope}) that every monotone path is coherent, hence $N_\ell=N_\ell^{\text{coh}}$ in this case.
\end{example}

\begin{remark}\label{rmk:CrossPolytopeCoherent}
Let $\b c$ be a generic direction for the $d$-cross-polytope. 
By \cite[Corollary 3.5]{BlackDeLoera2021monotone} the  coherent paths of length $\ell$ on $\CrossPol$  are in bijection with sequences in $\{-, +, 0\}^{d-1}\ssm\{\b 0\}$ having $\ell-1$ non-zero entries.
Since there are $\binom{d-1}{\ell-1}$ ways to choose these non-zero positions, and since for each such position there are two possible entries, we obtain~$N_\ell^{\text{coh}} = \binom{d-1}{\ell-1}2^{\ell-1}$.
This sequence is easily seen to be log-concave, adding another class of polytopes answering \Cref{prb:Unimodality} affirmatively.
\end{remark}

\begin{remark}\label{rmk:CycCoherent}
Let for $d\geq 4$ and  $\b c = \b e_1$. 
By \cite[Corollary 3.5]{Athanasiadis_2000}, the coherent paths on $\cyc_d(\b t)$ of length $\ell$ are in bijection with sign sequences $\{+, -\}^{n-2}$ with $\ell-1$ occurrences of $+$ and at most $d-1$ plateaus (\ie maximal subsequences of constant sign).
The number of sequences  in $\{+, -\}^{n-2}$ with $\ell-1$ occurrences of $+$ and \textbf{exactly} $d-1$ plateaus is equal to the following:
\begin{compactenum}
\item[$\bullet$] $2\binom{\ell-2}{\delta-1}\binom{n-\ell-2}{\delta-1}$, if $d-1 = 2\delta$ is even. 
\item[$\bullet$]  $\binom{\ell-2}{\delta}\binom{n-\ell-2}{\delta-1} + \binom{\ell-2}{\delta-1}\binom{n-\ell-2}{\delta}$, if $d-1 = 2\delta+1$ is odd.
\end{compactenum}

Hence, the number of coherent paths of length $\ell$ is given by:
\[
\sum_{2\delta \leq d-1} 2\binom{\ell-2}{\delta-1}\binom{n-\ell-2}{\delta-1} ~~+~~ \sum_{2\delta+1 \leq d-1} \left(\binom{\ell-2}{\delta}\binom{n-\ell-2}{\delta-1} + \binom{\ell-2}{\delta-2}\binom{n-\ell-2}{\delta}\right).
\]
Since $(\binom{\ell}{\delta})_\ell$ (for fixed $\delta$) is log-concave, \Cref{lem:OnLogConcaveSequences} (1) implies that so is  $\binom{\ell-2}{\delta-1}\binom{n-\ell-2}{\delta-1}$. 
Moreover, these products are symmetric with the same center of symmetry. This implies that $N_\ell$ is a sum of symmetric and unimodal sequences and hence  symmetric and unimodal by \Cref{lem:OnLogConcaveSequences} (3).
\end{remark}

\begin{example}\label{exm:SecondHypersimplex}
Let $\HypSimplTwo \eqdef \conv(\b x\in \{0, 1\}^d ~;~ \sum_{i=1}^d x_i = 2)$ be the \defn{second hypersimplex} and let~$\b c $ be a generic direction for $\HypSimplTwo$.
By \cite[Prop. 5.4]{poullot2024verticesmonotonepathpolytopes}, the number $N_\ell^{\text{coh}}$ is the coefficient of~$z^\ell$ in the polynomial $T_d + Q_d + C_d$ defined by:
$$\text{for } d\geq 4, \begin{pmatrix}T_{d+1} \\ Q_{d+1} \\ C_{d+1}\end{pmatrix} = \c M\begin{pmatrix}T_d \\ Q_d \\ C_d\end{pmatrix} ~~\text{with}~~ \c M = \begin{pmatrix}
  z&1+z&1+z \\ 0&1+z&z \\ z+z^2&0&1+z
\end{pmatrix} ~,~~ \begin{pmatrix}T_4 \\ Q_4
\\ C_4\end{pmatrix} = \begin{pmatrix} z^4+2z^3 \\ z^4 \\ 2z^4+2z^3 \end{pmatrix}.$$

It is conjectured (see \cite[Conj. 6.2]{poullot2024verticesmonotonepathpolytopes}) that these sequences are unimodal for $d\geq 4$. This has been confirmed, via computer experiments, for all $d\leq 150$, but the conjecture is open in general.

\begin{problem}
Is the sequence $(N^{\text{coh}}_\ell)_\ell$ defined above unimodal (and log-concave)?
\end{problem}
\end{example}

\begin{example}\label{exm:BlackOperations}
In his PhD Thesis \cite{Black-PhD}, Black derived formulas for the (total) number of monotone and coherent paths for the product $\pol\times \simplex_n$ of a polytope $\pol$ with a simplex and for the pyramid $\pol[Pyr](\pol)\eqdef \conv(\{\b 0\}\,\cup\,\pol\times\{1\})$ in terms of the corresponding numbers for $\pol$.
It is easy to refine these numbers accounting for the lengths of the paths.

In particular, if there are $N_\ell$ and $N^{\text{coh}}_\ell$ many $\b c$-monotone and coherent paths of length~$\ell$ on~$\pol$, respectively, then there are $\ell\cdot N_{\ell-1}$ and $\ell\cdot N^{\text{coh}}_{\ell-1}$ many $(\b c, 1)$-monotone and coherent paths of length~$\ell$ on $\pol\times [0, 1]$, respectively.
By combining \cite[Cor~3.2.2]{Black-PhD} and \Cref{lem:OnLogConcaveSequences}~(1), we obtain:

\begin{theorem}
If the sequence $(N_\ell)_\ell$ of the numbers of $\b c$-monotone paths (respectively coherent paths) on~$\pol$ of length $\ell$ is log-concave, then so is the sequence $(N'_\ell)_\ell$ of the number of $(\b c, 1)$-monotone paths (respectively coherent paths) on $\pol\times[0, 1]$.

Furthermore, if the polynomial $\sum_\ell N_\ell\, x^\ell$ is real-rooted, then so is the polynomial $\sum_\ell N'_\ell\, x^\ell$.
\end{theorem}

This provides another positive answer to \Cref{prb:Unimodality}.
However, the cases of $\pol\times\simplex_n$ (for $n\geq 2$) and $\pol[Pyr](\pol)$ are more convoluted.
By \cite[Thm~3.3.1]{Black-PhD}, the number of monotone (coherent) paths on $\pol[Pyr](\pol) $ is the sum over the vertices $\b v$ of $\pol$ of the number of monotone (coherent) paths from $\b v_{\min}$ to $\b v$.
This kind of sum may produce a unimodal sequence, even if the sequence $(N_\ell)_\ell$ (or $(N_\ell^{\text{coh}})_\ell$) for $\pol$ was not.
To motivate future research, we propose:

\begin{problem}
For a polytope $\pol$ and a direction $\b c$, let $\pol[Pyr]^0(\pol) \eqdef \pol$ and $\pol[Pyr]^{k+1}(\pol) \eqdef \pol[Pyr]\bigl(\pol[Pyr]^k(\pol)\bigr)$. Moreover, let $\b c^0 \eqdef \b c$ and $\b c^{k+1} \eqdef \bigl(\b c^k, 1 - \min_{\b x\in \pol[Pyr]^k(\pol)}\inner{\b x, \b c^k}\bigr)$.
For which $\pol$ and $\b c$, does there exists $k$ such that the number of $\b c^k$-monotone (and coherent) paths on $\pol[Pyr]^k(\pol)$ is unimodal?
\end{problem}

We conjecture that the answer to this question is ``for all $\pol$ and generic $\b c$''.

Instead of considering products in full generality, one could restrict to products $\simplex_n\times\simplex_m$  of simplices. 
For generic $\b c$, all $\b c$-monotone paths are coherent in this case and, by \cite[Corollary 3.2.3]{Black-PhD}, the number of such paths of length $\ell$ equals $N_\ell(n, m) \eqdef \sum_{k\geq 1} \binom{n-2}{k-1}\binom{m-2}{\ell-k-1}\binom{\ell}{k}$. 
We verified with a computer that the sequences $\bigl(N_\ell(n, m)\bigr)_\ell$ are log-concave for $n, m \leq 100$.
From \cite[Proposition 3.2.1]{Black-PhD}, one can also deduce a formula for any product of simplices.
This motivates:

\begin{problem}
 Is the number of $\b c$-monotone paths on $\simplex_{n_1}\times\dots\times\simplex_{n_r}$ log-concave (for $\b c$ generic)?
\end{problem}

One might try to tackle this problem either by brute force (\ie cleverly manipulating inequalities), or by finding an injection from the pairs of paths of length $\ell$ to the pairs formed by a path of length $\ell-1$ and a path of length $\ell+1$.
We also want to suggest another method: 
Namely, using log-concavity of the generalized hypergeometric functions.
Indeed, using the generalized hypergeometric function $_3F_2$, we have $N_\ell(n, m) = \ell \,_3F_2(1 - \ell, 2 - m, 2 - n;1, 2;-1)$.
The works of Kalmykov, Karp, Sitnik, and others shed light on the domains of log-concavity of the functions~$_pF_q$, see~\cite{KarpSitnik2010-LogConcavityHypergeometricFunctions,KalmykovKarp2017-LogConcavityAndTuranInequalityHypergeometricFunctions} and the references therein:
one should try to deduce log-concavity for sequences of sums of products of binomial coefficients, from the log-concavity of such functions.
\end{example}

\section{Negative examples}\label{sec:NegativeExamples}

We now prove \Cref{thm:IntroA}, providing polytopes for which \Cref{prb:Unimodality} has a negative answer.

As a warm-up, consider the $2$-dimensional case: 
For a polygon $\pol$ and a generic direction $\b c$, there exist exactly two $\b c$-monotone paths, both of which are coherent.
It is easy to construct polygons (and directions) for which the lengths of these two paths differ by at least 2 (see \eg the projected $3$-polytope in \Cref{fig:ExamplePolytopeWithMonotonePaths}).
Consequently, none of the sequences $(N_\ell)_\ell$ and $(N^{\text{coh}}_\ell)_\ell$ is unimodal, as both contain two non-consecutive $1$s.
This answers \Cref{prb:Unimodality} in the negative for $d = 2$.

Let $d \geq 3$. 
In personal communication (posterior to the writing of this section), Alexander Black informed us that Christopher Eur had found  an example of a $3$-dimensional polytope with $7$ vertices, $13$ edges and $8$ facets ($6$ triangles, $2$ squares), and a  direction $\b c\in\R^3$, for which the sequence $N_\ell$ is unimodal but not log-concave. 
This shows that the strengthening of \Cref{prb:Unimodality} already fails in dimension $3$. 
Throughout the rest of this section, we make \Cref{thm:IntroA} explicit.

\subsection{Lopsided cubes}\label{ssec:LopsidedCube}
The goal is to  prove \Cref{thm:IntroA} (a) by providing an explicit construction.
The main idea is based on the following observation: 
The monotone paths of any polytope with a $2$-colorable graph (\eg the standard $d$-cube, the permutahedron) are either all of even length or all of odd length. 
In particular, if there are monotone paths of \emph{different} lengths, then, as in dimension $2$, the sequence $(N_\ell)_\ell$ is not unimodal since it has internal $0$s. 
The same reasoning applies to coherent paths.

Since the monotone paths of  the standard cube $[0,1]^d$ and the permutahedron all have the same length (see \Cref{exm:GradedGraph}), we cannot use these polytopes directly to make the described idea work, but we need to modify them slightly.

\begin{figure}
\centering
\includegraphics[width=0.8\linewidth]{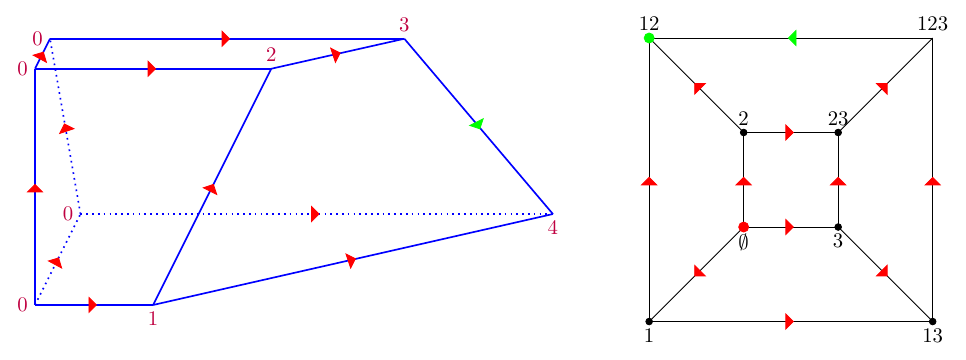}
\caption[Lopsided 3-cube]{(Left) The (oriented) lopsided $3$-cube $\Lop[3]$, its vertices labeled by their first coordinate.
(Right) The graph of $\Lop[3]$ directed by $(1, 1, 1)$, its vertices $\b u_X$ labeled by $X$, its source and sink encircled in \textcolor{red}{red} and \textcolor{green}{green}, respectively.}
\label{fig:Lopsided3Cube}
\end{figure}

\begin{example}\label{exm:3Lopsided}
The polytope $\Lop[3]$ defined as the convex hull of the points
$$\begin{array}{rclcrclcrclcrcl}
\b u_{\emptyset} & = & (0, 0, 0), & ~ & \b u_{\{2\}} & = & (0, 1, 0), & ~ & \b u_{\{3\}} & = & (0, 0, 1), & ~ & \b u_{\{2, 3\}} & = & (0, \frac{1}{3}, 1), \\
\b u_{\{1\}} & = & (1, 0, 0), & ~ & \b u_{\{1, 2\}} & = & (\textcolor{red}{4}, 1, 0), & ~ & \b u_{\{1, 3\}} & = & (2, 0, 1), & ~ & \b u_{\{1, 2, 3\}} & = & (\textcolor{red}{3}, \frac{1}{3}, 1)
\end{array}$$
is called the \defn{lopsided $3$-cube} (see \Cref{fig:Lopsided3Cube} (left) for an illustration). 
The polytope $\Lop[3]$ is combinatorially isomorphic to a $3$-cube and, for $\b c = (1, 1, 1)$, its directed graph is obtained from the directed graph of $[0,1]^3$ by reversing the direction of the edge $\b u_{\{1, 2\}} \to \b u_{\{1, 2, 3\}}$, see \Cref{fig:Lopsided3Cube} (right).
It has 2 and 4 monotone paths of length 2 and 4, respectively, and all of these are coherent.
Hence, the sequences $(N_\ell)_\ell$ and $(N_\ell^{\text{coh}})_\ell$ are \textbf{not} unimodal.
\end{example}

Using prisms, we now extend the previous construction to an arbitrary dimension.

\begin{definition}
Let $\pol\subseteq\R^d$ be a polytope. The \defn{$k$-fold (standard) prism $\pol[Prism]_k \pol$} of $\pol$ is the polytope $\pol[Prism]_k \pol \eqdef \conv\Bigl((\b v, \b e_X) ~;~ \b v\in \pol,~ X\subseteq[k]\Bigr) \subseteq\R^{d+k}$, where $\b e_X = \sum_{i\in X}\b e_i\in \R^k$ for $X\subseteq[k]$.

Let $d\geq 4$. The \defn{lopsided $d$-cube $\Lop$} is the $(d-3)$-fold standard prism over $\Lop[3]$, that is\linebreak $\Lop \eqdef \conv(\b u_X ~;~ X\subseteq[d])$, where $\b u_X \eqdef(\b u_{X_{\geq d-2}},\b e_{X_{\leq d-3}})$ for $X\subseteq[d]$,   
 $X_{\leq d-3} = X\cap [d-3]$, and $X_{\geq d-2} = \bigl\{i-d+3 ~;~ i\in X\cap \{d-2, d-1, d\}\bigr\}$.
\end{definition}

In the following, we count (coherent) monotone paths on $\Lop$, 
iterating \cite[Prop. 3.2.1]{Black-PhD}. 

\begin{lemma}\label{lem:Prims}
Let $\pol\subseteq \R^d$ be a polytope. 
Let $\b c\in \R^d$ and  let $N_\ell$ be the number of $\b c$-monotone paths of length $\ell$ on $\pol$.
Then the number of $(\b c, 1, \dots, 1)$-monotone paths of length $k+\ell$ on $\pol[Prism]_k\pol$ is equal to $\frac{(k+\ell)!}{\ell!}N_\ell$. 
\end{lemma}

\begin{proof}
%
To keep this article self-contained, we provide a simple proof of this lemma.
Let $\c L$ be a $(\b c, 1,\dots, 1)$-monotone path  of length $k+\ell$ on $\pol[Prism]_k\pol$.
Then, by definition, $\c L$ has exactly one edge  parallel to $\b e_i$ for every $i\in [d+1, d+k]$, and the remaining $\ell$ edges are  parallel to an edge of $\pol$ (directed according to $\b c$). 
Consequently, $\c L$ is obtained by shuffling a $\b c$-monotone path of length $\ell$ on $\pol$ and a $(1,\ldots,1)$-monotone path on $[0,1]^k$ together (while maintaining the relative order of their edges). 
Vice versa, any shuffle of a $\b c$-monotone path of length $\ell$ on $\pol$ and a $(1,\ldots,1)$-monotone  on $[0,1]^k$ gives rise to a $(\b c, 1, \dots, 1)$-monotone path on $\pol[Prism]_k\pol$.
The number of such shuffles is $\binom{k+\ell}{\ell}\cdot N_\ell\cdot k! = \frac{(k+\ell)!}{\ell!}\cdot N_\ell$, which shows the claim.
\end{proof}

\begin{example}\label{exm:SimpleNonExample}
The lopsided cube $\Lop$ is combinatorially isomorphic to $[0,1]^d$, and for $\b c = (1, \dots, 1)$ its directed graph $G_{\Lop,\b c}$  differs from the one of the standard cube $[0,1]^d$ by reversing the edges $\b u_X \to \b u_Y$ for which $X_{\geq d-2} = \{1, 2\}$ and $Y_{\geq d-2} = \{1, 2, 3\}$.
W have $\b v_{\min} = \b u_{\emptyset} = (0, \dots, 0)$ and $\b v_{\max} = \b u_{[d-1]} = (4, 1, 0,1, \dots, 1)$.
Applying \Cref{lem:Prims} to \Cref{exm:3Lopsided} yields that the number of (coherent) monotone paths on $\Lop$ of length $\ell$ is given by the following non-unimodal sequence:
\begin{center}
\begin{tabular}{c|ccc|c}
   $\ell$ & $d-1$ & $d$ & $d+1$ & total \\ \hline
   $N_\ell = N^{\text{coh}}_\ell$ & $(d-1)!$ & $0$ & $\frac{1}{6}(d+1)!$ & $(d-1)! \cdot \left(1+\frac{d(d+1)}{6}\right)$
\end{tabular}.
\end{center}

\end{example}

The previous lemma also shows that every $\b c$-monotone path on $\Lop$ is coherent.
This counterexample proves the following statement, which makes \Cref{thm:IntroA} (a) more explicit.

\begin{theorem}\label{cor:LopsidedCube}
For all $d\geq 3$, there exists a $d$-dimensional polytope $\pol\subseteq\R^d$, combinatorially isomorphic to a $d$-cube, and a direction $\b c\in \R^d$, such that the sequences $(N_\ell(\pol,\b c))_\ell$ and $(N_\ell(\pol,\b c)^{\text{coh}})_\ell$ are \textbf{not} unimodal.
\end{theorem}

\begin{figure}
\centering
\includegraphics[width=0.38\linewidth]{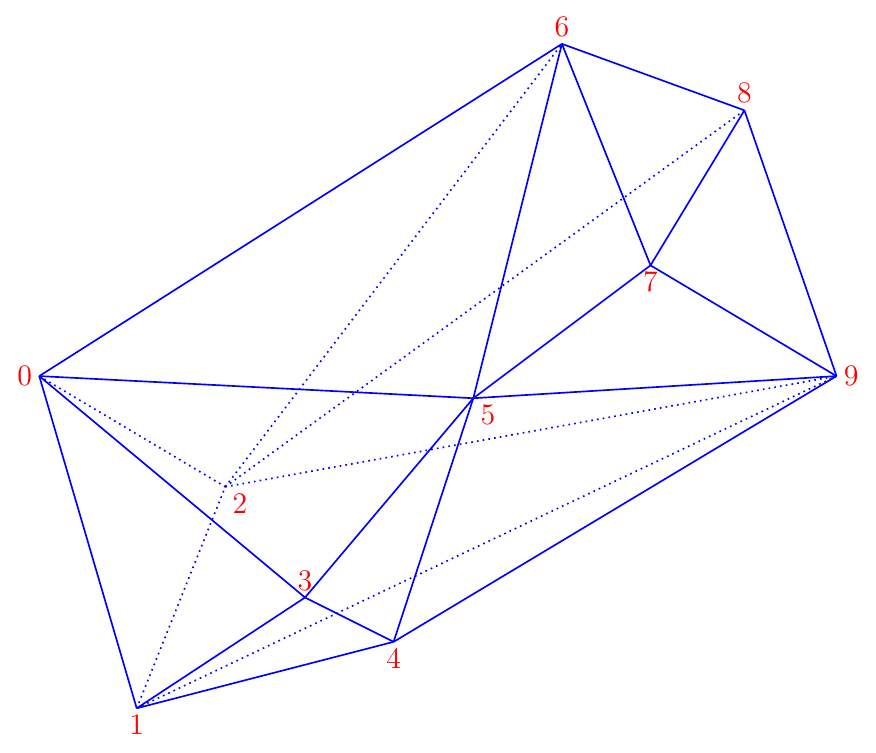}\hspace{0.5cm}
\includegraphics[width=0.48\linewidth]{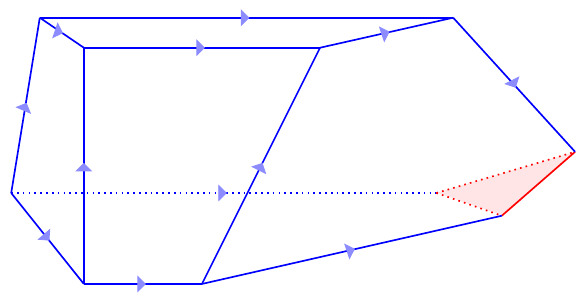}
\caption[Simplicial and simple polytope with non-unimodal $(N_\ell)_\ell$]{(Left) The $3$-dimensional simplicial polytope $\pol_{10}$ from \Cref{exm:SimplicialNonExample}. (Right) The $3$-dimensional simple polytope from \Cref{exm:SimpleNonExample}.}
\label{fig:Simplicial3D}
\end{figure}

\begin{remark}
The constructed counterexamples have internal $0$s. 
We now sketch a method that seems suitable for constructing counterexamples without internal $0$s.
We start with $\Lop$ and intersect it with a half-space $H_{\b a, b}^- = \{\b x\in \R^d~;~ \inner{\b x, \b a} \leq b\}$ that contains all vertices of $\Lop$ except $\b v_{\max}=\b u_{[d-1]}$. 
We set $\pol_{\b a, b} = \Lop\cap H_{\b a, b}^-$.
If $\b a$  and $(1, \dots, 1)$ are linearly independent, then $(1, \dots, 1)$ is generic for $\pol_{\b a, b}$.
As $\Lop$ is a simple polytope, $G_{\pol_{\b a, b}, (1, \dots, 1)}$ is obtained from $G_{\Lop, (1, \dots, 1)}$ by replacing the vertex $\b u_{[d-1]}$ by an oriented clique and connecting its vertices to the neighbors of $u_{[d-1]}$.  
Performing this procedure for $d=3$ and moving another vertex, yields a $3$-dimensional simple polytope with $(N_\ell)_{2\leq\ell\leq 5} = (1, \textcolor{red}{2}, \textcolor{red}{1}, \textcolor{red}{3})$ (see \Cref{fig:Simplicial3D} (right)).
For $d = 4$, taking \eg $\b a = (2, 4, 3, 3)$, $b = 20.5$ produces $\pol_{\b a, b}$ with $(N_\ell)_{4\leq\ell\leq 8} = (6, \textcolor{red}{22}, \textcolor{red}{6}, \textcolor{red}{8}, 4)$. 
Computer experiments suggest that this method tends to create non-unimodal sequences without internal $0$s.
\end{remark}

\subsection{Simplicial polytopes}
Until now, the examples that provided a negative answer to \Cref{prb:Unimodality} were exclusively simple polytopes.
One might wonder what happens for simplicial polytopes. 
The purpose of this subsection is to show that the answer to  \Cref{prb:Unimodality} can be negative for simplicial polytopes.
More precisely, \Cref{exm:SimplicialNonExample} proves the following refinement of \Cref{thm:IntroA} (b):

\begin{theorem}\label{thm:counterexampleSimplicial}
Let  $\b c = (1, 0, 0)$. There exists a $3$-dimensional polytope $\pol\subseteq\R^3$ with $10$ vertices such that $(N_\ell(\pol,\b c))_\ell$ is \textbf{not} unimodal.
\end{theorem}


\begin{example}\label{exm:SimplicialNonExample}
Let $\pol_{10}\subseteq \R^3$ be the polytope with vertices:
$(0, 0, 0)$, $(1, -5, -5)$, $(2, 0, -5)$, $(3, -5, 0)$, $(4, -6, 0)$, $(5, -3, 5)$, $(6, 5, 5)$, $(7, 0, 5)$, $(8, 5, 2)$ and $(9, 0, 0)$ (see \Cref{fig:Simplicial3D} (Left), where the vertex labels are equal to the first coordinates).  
The next table displays the number of $(1,0,0)$-monotone paths on $\pol_{10}$ 
(from $(0, 0, 0)$ to $(9, 0, 0)$).
\begin{center}
\begin{tabular}{c|ccccccc|c}
  $\ell$ & 2 & 3 & 4 & 5 & 6 & 7 & 8 & total \\ \hline
  $N_\ell$ & 3 & 8 & \textcolor{red}{{12}} & \textcolor{red}{{11}} & \textcolor{red}{{12}} & 6 & 1 & 53
\end{tabular}.
\end{center}
\end{example}

\begin{remark}\label{rmk:SimplicialOnSphereNotUnimodal}
Let $\b c=(1,0,0)$ and let $\b b$ be the barycenter of $\pol_{10}$.
The vertices of $\pol_{10}$ can be embedded in the 2-dimensional sphere $\mathbb{S}^2 = \{\b x\in \R^3 ~;~ \|\b x\| = 1\}$ using the map $\b x \mapsto \frac{\b x - \b b}{\|\b x - \b b\|}$.
Let $\overline{\pol}_{10}$ be the convex hull of the images of the vertices of $\pol_{10}$. 
This polytope is simplicial but its directed graph $G_{\overline{\pol}_10,\b c}$  is different from  $G_{\pol_{10},\b c}$.  
One can verify that  $(N_\ell(\overline{\pol}_{10},\b c))_{1\leq \ell \leq 8}=(4, 8, \textcolor{red}{10}, \textcolor{red}{8}, \textcolor{red}{11}, 6, 1)$, which is not unimodal.

As slightly modifying the coordinates of the vertices of $\overline{\pol}_{10}$ does not change its directed graph, there exists a subset $\c A\subseteq (\mathbb{S}^2)^{10}$, which is not of measure $0$, such that if $(\b x_1, \dots, \b x_{10})\in \c A$, the sequence of the numbers of $\b c$-monotone paths on $\conv(\b x_1, \dots, \b x_{10})$ of length $\ell$ is not unimodal.
Equivalently, when constructing a polytope as the convex hull of $10$ points chosen uniformly at random on $\mathbb{S}^2$, there is a strictly positive probability that the sequence of the number of $\b c$-monotone paths of length $\ell$ is not unimodal (\ie that the answer to \Cref{prb:Unimodality} is ``no'').
The reader should keep this in mind while reading \Cref{sec:RandomCase}.
\end{remark}

\begin{problem}
Find a simplicial polytope $\pol$ and a direction $\b c$ such that  $(N_\ell^{\text{coh}}(\pol, \b c))_\ell$ is not unimodal (we found a non-log-concave example, but do not present it here).
\end{problem}

\subsection{Loday's associahedron of dimension 5}\label{ssec:Asso5}
For $n\geq 3$, \defn{Loday's $n$-associahedron $\Ass\subseteq \R^n$} is the $(n-1)$-dimensional generalized permutahedron (\ie its edge directions are $\b e_i - \b e_j$ for some $1\leq i<j\leq n$) with the following facet-description:
$$\Ass = \left\{\b x\in \R^n ~;~ \begin{array}{l} \sum_{i=1}^n x_i = 0 \\ \sum_{i\in I} x_i \geq \binom{|I|+1}{2} ~~\text{ for }~~ \emptyset\ne I=[a,b]\subsetneq[n] \end{array} \right\}.$$
It is known that the graph of $\Ass$, directed by $\b c = (1, 2, \dots, n)$, is the Hasse diagram of the Tamari lattice (see \cite{Loday2004-Associahedron,PilaudSantosZiegler-CelebratingLodayAssociahedron}).
Hence, $N_\ell=N_\ell(\Ass,\b c)$ is equal to the number of maximal chains in the Tamari lattice.
The latter was computed by Nelson \cite[Thm. 5.9]{Nelson2017-MaximalChainsTamariLattice} and discussed in the general context of graph associahedra in \cite{DahlbergFishel2024-MaximalChainsGraphAssociahedra}.
Nelson gives the following sequence $N_\ell$ for $n = 6$, see also \href{https://oeis.org/A282698}{OEIS A282698}, which we completed by computing~$N^{\text{coh}}_\ell=N^{\text{coh}}_\ell(\Ass,\b c)$:

\begin{small}
\begin{center}
\begin{tabular}{c|ccccccccccc|c}
 $\ell$ & 6 & 7 & 8 & 9 & 10 & 11 & 12 & 13 & 14 & 15 & 16 & total \\ \hline
$N_\ell$ & 1 & 20 & 112 & 232 & \textcolor{red}{382} & \textcolor{red}{348} & \textcolor{red}{456} & 390 & 420 & 334 & 286 & 2981 \\
$N^{\text{coh}}_\ell$ & 1 & 20 & 105 & 206 & \textcolor{red}{332} & \textcolor{red}{274} & \textcolor{red}{332} & 270 & 206 & 122 & 142 & 2010
\end{tabular}.
\end{center}
\end{small}

Both sequences are \textbf{not} unimodal. 
This example proves \Cref{thm:IntroA} (c):

\begin{theorem}\label{thm:CounterexampleGenPer}
There exists a $5$-dimensional generalized permutahedron $\pol$ and a direction $\b c\in \R^6$ such that $(N_\ell(\pol,\b c))_\ell$ and $(N_\ell^{\text{coh}}(\pol,\b c))_\ell$ are \textbf{not} unimodal.
\end{theorem}

\begin{remark}
Similar computations yield $N_\ell$ and $N_\ell^{\text{coh}}$ for the $4$-dimensional polytope $\Ass[5]$ with respect to $\b c = (1,2,3,4,5)$.
The sequence $(N_\ell)_\ell$ is unimodal while the sequence $(N_\ell^{\text{coh}})_\ell$ is not.
\begin{center}
\begin{tabular}{c|c|ccccccc|c}
  Source & $\ell$ & 4 & 5 & 6 & 7 & 8 & 9 & 10 & total \\ \hline
  \cite{Nelson2017-MaximalChainsTamariLattice}, \href{https://oeis.org/A282698}{A282698} & $N_\ell$ & 1 & 10 & 22 & 22 & 18 & 13 & 12 & 98 \\
  Our computation & $N^{\text{coh}}_\ell$ & 1 & 10 & 21 & 21 & \textcolor{red}{{18}} & \textcolor{red}{{9}} & \textcolor{red}{{10}} & 90
\end{tabular}.
\end{center}
\end{remark}

Moreover, for $n = 7$, the sequence $(N^{\text{coh}}_\ell)_\ell$ has several non-unimodal sub-triples; but for $n = 8$, the sequence $(N^{\text{coh}}_\ell)_\ell$ is unimodal.
To our knowledge, counting coherent paths of Loday's associahedron is an open problem.
Nelson interprets maximal chains in the Tamari lattice as tableaux. Unfortunately,  the notion of the coherence of a path is not equivalent to the realizability of Young tableaux developed in \cite{MallowsVanderbei-YoungTableauxAsSums,AraujoBlackBurcroffGaoKruegerMcDonough-RealizableStandardYoungTableaux} (see also \cite[Section 8]{BlackSanyal-FlagPolymatroids} for a more general perspective around coherent paths on the permutahedron).

\begin{problem}
Describe and count the coherent paths on Loday's associahedron for $\b c = (1, 2, \dots, n)$.
\end{problem}

\subsection{Polytopes with 0/1-coordinates}
The counterexamples from the previous subsections shared the feature that their vertices had many different coordinate values.
This was necessary to construct  the quite convoluted behaviors of their paths.
In this section, we  present a counterexample only using $0$ and $1$ as coordinates. 
In particular, we make \Cref{thm:IntroA} (d) more explicit.


A polytope $\pol\subseteq\R^n$ is a \defn{$0/1$-polytope} if $\b v\in\{0,1\}^n$ for each vertex $\b v$ of $\pol$.
As every $0/1$-vector is determined by its non-zero coordinates, every  $0/1$-polytope is of the form 
\[
\pol_{\c X} \eqdef \conv\left(\b e_X ~;~ X\in \c X\right),
\]
where $\b e_X \eqdef\sum_{i\in X} \b e_i$ and $\c X\subseteq 2^{[n]}$. 

Let $\b c_{\text{lex}} = (2^1, 2^2, \dots, 2^n)\in \R^n$. 
The orientation of the graph $G_{\pol_{\c X}, \b c_{\text{lex}}}$ is induced by the  lexicographic order:
if $[\b e_X, \b e_Y]$ is an edge of $\pol_{\c X}$, then $\inner{\b e_X, \b c_{\text{lex}}} < \inner{\b e_Y, \b c_{\text{lex}}}$ if and only if\linebreak $X = \{x_r > \dots > x_1\}$ is lexicographically smaller than $Y = \{y_s > \dots > y_1\}$.
The following example provides a $0/1$-polytope for which the sequence $(N_\ell)_\ell$ of the number of $\b c_{\text{lex}}$-monotone paths of length $\ell$ is not unimodal.

\begin{example}\label{exm:SimplicialComplexNotUnimodal}
Let $\c X$ be the collection of all sets contained in $\{1,4\}$, $\{1,2,3,5\}$ or $\{2,3,4,5\}$. 
The polytope $\pol_{\c X}\subseteq\R^5$ is neither simple nor simplicial and its $f$-vector is equal to $(1, 25, 75, 90, 51, 13, 1)$.
For the direction $\b c_{\text{lex}} = (2^1, 2^2, 2^3, 2^4, 2^5)$, the sequence $(N_\ell(\pol_{\c X},\b c_{\text{lex}}))_\ell$ is the following:
\begin{center}
\begin{tabular}{c|cccccc|c}
  $\ell$ & 3 & 4 & 5 & 6 & 7 & 8 & total \\ \hline
  $N_\ell$ & 2 & 36 & \textcolor{red}{{96}} & \textcolor{red}{{76}} & \textcolor{red}{{84}} & 36 & 330
\end{tabular}.
\end{center}
\end{example}

This example proves \Cref{thm:IntroA} (d), but leaves the next problem open.

\begin{theorem}\label{thm:Counterexample01}
There exists a $5$-dimensional $0/1$-polytope $\pol\subseteq\R^5$ with $25$ vertices such that, for $\b c_{\text{lex}} = (2^1, 2^2, 2^3, 2^4, 2^5)$, the sequence $(N_\ell(\pol, \b c_{\text{lex}} ))_\ell$ is \textbf{not} unimodal.
\end{theorem}

\begin{problem}
Find a $0/1$-polytope $\pol$ and a direction $\b c$ such that $(N^{\text{coh}}_\ell(\pol,\b c))_\ell$ is not  unimodal.
\end{problem}

Note that linear optimization on $0/1$-polytopes has been largely studied; see \cite{BlackDeLoeraKaferSanita-SimplexMethodOn01polytopes} and its section ``Prior work and context''.
Since there are polynomial algorithms for finding short paths, the above problem is more of theoretical importance than practical relevance.

\section{Random case}\label{sec:RandomCase}

We have seen that even for sufficiently nice classes of polytopes, including simple and simplicial polytopes, generalized permutahedra and $0/1$-polytopes, \Cref{prb:Unimodality} has a negative answer in general.
However, all counterexamples that we gave were rather special in the way we constructed them. 
Moreover, after experimenting with random polytopes with vertices on the sphere, it seems rather hard to find an example of a polytope that contradicts \Cref{prb:Unimodality} (see \Cref{rmk:SimplicialOnSphereNotUnimodal}). 
In the following, we make this intuition precise.
Although the question of \emph{monotone} paths is a problem in dimension $d$, understanding \emph{coherent} paths amounts to studying 2-dimensional projections of a $d$-dimensional polytope.
Since the latter seems to be more tractable, we focus on coherent paths.
We start by formulating our main result (\Cref{thm:IntroB}).

\begin{theorem}\label{thm:CLTmain}
Let $d\geq 4$. 
Fix (deterministically or at random) linearly independent vectors $\b c, \b \omega\in \R^d$.
Let $Z_1, \dots, Z_n$ be points taken uniformly independently and uniformly at random on the sphere $\mathbb{S}^{d-1}\subseteq\R^d$ and let $\pol_n = \conv(Z_1, \dots, Z_n)$.
Then, the length $L_n$ of the coherent $\b c$-monotone path captured by $\b\omega$ on $\pol_n$ admits a central limit theorem, \ie
\[
\frac{L_n - \E (L_n)}{\sqrt{\Var (L_n)}} \xrightarrow[n\to+\infty]{} U ~~~~\text{ (convergence is in distribution)},
\]
with $U\sim \c N(0, 1)$ a standard normally distributed random variable with expectation $0$ and variance $1$. 
Moreover, $\E(L_n) \sim c\, n^{\frac{1}{d-1}}$ for some $c > 0$ and $c'\, n^{\frac{1}{d-1} - a}\leq \Var(L_n) \leq c''\, (\log n)^{3 - \frac{1}{d-1}} \, n^{\frac{1}{d-1}}$ for any $a > 0$ and some $c', c'' > 0$.
\end{theorem}

We proof this theorem at the end of \Cref{sssec:CLT}, using \Cref{cor:LnIsF0Qn,cor:expectation,cor:VarianceLower,cor:VarianceUpper,cor:CLT}.

The methods we use to prove our results  come mainly from the theory of random convex bodies and polytopes and can be considered standard in this area.
In particular, they have found multiple applications; see \cite[Section 6]{ReyPecati}, and also \cite{Thale-CentralLimitTheorem,ThaleTurchiWespi-IntrinsicVolumes,BesauRosenThale-Projective}. 
Moreover, they frequently rely on combinatorial and geometric arguments and ideas, and we hope to be able to convince the reader that they might be helpful for similar problems.

As our background also lies in combinatorics and polytope theory, we tried to keep this section as accessible as possible to a combinatorial audience.
In particular, we will provide a detailed description of the methods, including illustrations to make the  ideas more clear. 
We also included a cheat sheet of important formulas as a service to the reader in \Cref{ssec:CheatSheet}. 
\Cref{ssec:ProbabilisticModel} explains the probabilistic model at stake by detailing the interaction between coherent paths (\ie projections to $\R^2$) and the uniform distribution on the sphere $\mathbb{S}^{d-1}$.
\Cref{ssec:CLT} proves \Cref{thm:CLTmain} by analyzing the behavior of \defn{$\beta$-polygons} in the plane, for $\beta > 0$.

\subsection{The probabilistic model}\label{ssec:ProbabilisticModel}

Let $Z_1,\dots,Z_n$ be independently uniformly distributed points on the $(d-1)$-dimensional sphere $\mathbb{S}^{d-1} = \{x\in \R^d ~;~ \|x\|=1\}$, and let \defn{$\pol_n=\conv(Z_1,\dots,Z_n)$} be the induced random polytope.
Let~$\b c\in \R^d$ be a fixed direction. By rotational symmetry, we assume that $\b c = \b e_1$.
We are interested in the number $N_\ell^{\text{coh}}=N_\ell^{\text{coh}}(\pol_n,\b c)$ of coherent $\b c$-monotone paths of length~$\ell$ on $\pol_n$. 
To show that the sequence $(N^{\text{coh}}_\ell)_\ell$ is ``statistically unimodal'', we study the probability distribution of the length of a random coherent path.

Let \defn{$L(\b \omega, \pol_n)$} denote the random variable that is the length of the coherent monotone path on $\pol_n$ captured by $\b\omega$, where $\b\omega\in \mathbb{S}^{d-1}$.
By \Cref{def:CoherentPaths}, $L(\b\omega, \pol_n)$ is the length, \ie the number of edges, of the upper path of the polygon obtained by projecting $\pol_n$ onto the plane spanned by $\b c$ and $\b\omega$.
We aim at understanding the distribution of $L(\b\omega, \pol_n)$ for large $n$.
In the following, we denote by $\pi_{\b c, \b\omega}$ the orthogonal projection from $\R^d$ to the plane spanned by $\b c$ and $\b\omega$, see \Cref{fig:Samples} (Left).
The next lemma, which is a special case of \cite[Lemma 4.3(a)]{KabluchkoTemesvariThale-IntrinsicVolumeBetaPolytopes}, shows that the random variables $\pi_{\b c, \b\omega}(Z_1),\dots, \pi_{\b c, \b\omega}(Z_n)$ follow a \defn{($2$-dimensional) $\beta$-distribution} (for a specific value $\beta$).

\begin{lemma}[{Adapted from \cite[Lemma 4.3(a)]{KabluchkoTemesvariThale-IntrinsicVolumeBetaPolytopes}}]\label{lem:Kabluchko}
Let $\b E$ be a $2$-dimensional plane in $\R^d$ and let $\pi_{\b E}:\R^d \to \b E$ be the orthogonal projection onto $\b E$. If a random variable $Z$ is distributed according to the uniform distribution on the sphere $\mathbb{S}^{d-1}$, then the projected random variable $\pi_{\b E}(Z)$ is distributed according to the probability density
\[
f_{2,\beta_d}(\b x) = C_{2,\beta_d}\, \left(1-\|\b x\|^2\right)^{\beta_d} ~~~\text{ for }\b x\in \B,
\]
where $\beta_d=\frac{d}{2} - 2$ and $C_{2,\beta_d}=\frac{1}{\pi}\frac{\Gamma(\beta_d+2)}{\Gamma(\beta_d+1)}$.
\end{lemma}


\begin{example}
For $d = 3$, we have $\beta_3=-\frac{1}{2}$ and the density $f_{2,\beta_3}(\b x)$ becomes higher the closer a point~$\b x$ is to the boundary of the disk $\B=\{x\in \R^2~;~\|x\|=1\}$, see \Cref{fig:3Dplots} (Left).

For $d = 4$, we have $\beta_4 = 0$ and the $2$-dimensional $\beta_4$-distribution is the uniform distribution on~$\B$.
Reitzner \cite{Reitzner-CLT} proved that for the convex hull of $n$ uniformly distributed independent random points in a convex set (in any dimension), the number of $k$-faces satisfies a central limit theorem (for any $k$).

For $d \geq 5$, we have $\beta_d> 0$ and the density $f_{2,\beta_d}(\b x)$ gets lower the closer a point $\b x$ is to the boundary of $\B$. 
In particular, the higher the dimension, the more concentrated the distribution is around the center of the disk, and the sparser it gets towards the boundary, see \Cref{fig:3Dplots,fig:Samples}.

Due to this different behavior of the density for $d\leq 3$ and $d\geq 4$ we will only consider the case $d\geq 4$, in the following.
\begin{figure}
\centering
\includegraphics[width=0.95\linewidth]{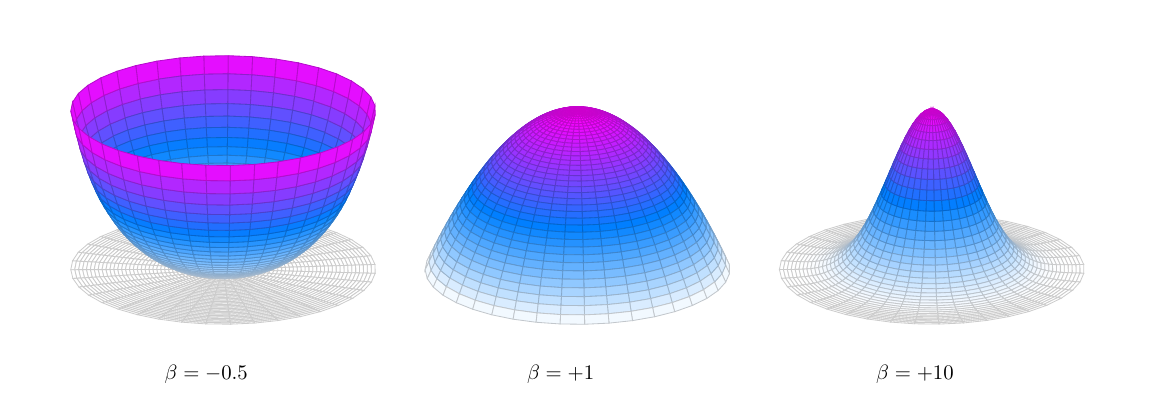}
\caption[$\beta$-distributions for negative, zero and positive $\beta$]{Examples of $\beta$-distributions for $\beta\in \{-0.5, +1, +10\}$, or equivalently $d\in\{3, 6, 22\}$.}
\label{fig:3Dplots}
\end{figure}
\end{example}

By \Cref{lem:Kabluchko}, the projected points $X_i \eqdef\pi_{\b c, \b\omega}(Z_i)$ are independently and identically distributed (\defn{i.i.d.} for short) with probability density function $f_{2,\frac{d}{2} - 2}$.
The next definition is essential:

\begin{definition}
For $X_1, \dots, X_n$ independently and identically $\beta$-distributed random points on the disk~$\B$ with $\beta = \frac{d}{2} - 2$, we set \defn{$\polytopeQ_n = \conv(X_1, \dots, X_n)$}.
\end{definition}

We use \defn{$f_1^{\text{up}}(\polytopeQ_n)$} and \defn{$f_1^{\text{low}}(\polytopeQ_n)$}  to denote the number of edges of the upper and the lower path of $\polytopeQ_n$, respectively. 
By \Cref{def:CoherentPaths}, we have $L(\b\omega,\pol_n) = f_1^{\text{up}}(\polytopeQ_n)$ for $\polytopeQ_n = \pi_{\b c, \b\omega}(\pol_n)$.
By symmetry, $f_1^{\text{low}}$ and $f_1^{\text{up}}$ are identically distributed, but \textbf{not} independent.
One can show that, after normalization, they converge to a normal distribution if and only if \defn{$f_1(\polytopeQ_n)$} $=f_1^{\text{up}}(\polytopeQ_n)+f_1^{\text{low}}(\polytopeQ_n)$ satisfies a central limit theorem (see the end of \Cref{sssec:CLT} for details).
As $\polytopeQ_n$ is a polygon, its number of vertices \defn{$f_0(\polytopeQ_n)$} satisfies $f_0(\polytopeQ_n) = f_1(\polytopeQ_n)$.
We therefore need to show that $f_0(\polytopeQ_n)$ obeys a central limit theorem.
We summarize this in the next corollary:

\begin{corollary}\label{cor:LnIsF0Qn}
Let $\b c, \b \omega\in \mathbb{S}^{d-1}$, and let $\pol_n = \conv(Z_1, \dots, Z_n)$, where $Z_1, \dots, Z_n$ are i.i.d. points on $\mathbb{S}^{d-1}$. Let $\polytopeQ_n = \conv(X_1, \dots, X_n)$, where $X_1, \dots, X_n$ are independently $\beta$-distributed with $\beta = \frac{d}{2} - 2$. 
Then, the random variables $L(\b\omega, \pol_n)$, $f_1^{\text{up}}(\polytopeQ_n)$ and $f_1^{\text{low}}(\polytopeQ_n)$ have the same distribution.
\end{corollary}

\begin{proof}
By \Cref{lem:Kabluchko}, the random variables $L(\b\omega, \pol_n)$ and $f_1^{\text{up}}(\polytopeQ_n)$ are identically distributed. 
Moreover, due to rotational symmetry, $f_1^{\text{up}}(\polytopeQ_n)$ and $f_1^{\text{low}}(\polytopeQ_n)$ also have the same distribution.
%
\end{proof}


\begin{example}
With a computer, we can sample $n$ points on a $(d-1)$-dimensional sphere, project them into dimension $2$ (by forgetting all but their first two coordinates) and construct the convex hull of the projected points.
Varying $d$, and hence $\beta_d$, gives rise to different pictures, see \Cref{fig:Samples}.
Again, one observes that the points get more concentrated around the center if $\beta$ grows.

\begin{figure}
  \centering
\includegraphics[width=0.28\linewidth]{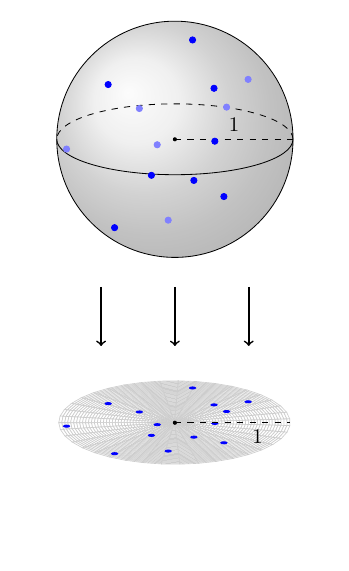}\includegraphics[width=0.71\linewidth]{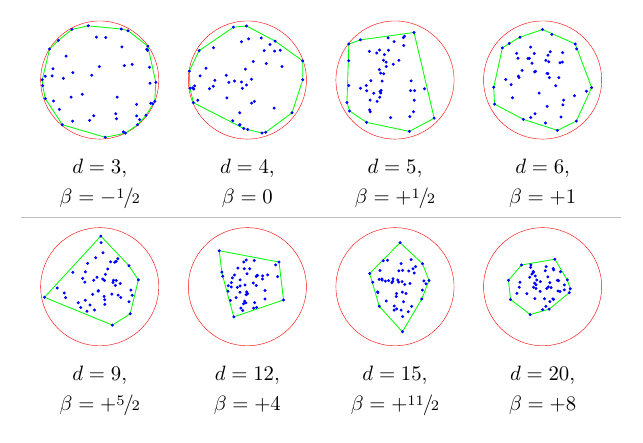}
\caption[Convex hull of $\beta$-distributed points in the disk for various $\beta$]{For each $d$, we sampled $n = 50$ points uniformly at random on the $(d-1)$-sphere and projected them onto the disk $\B$.
The convex hull $\polytopeQ_n$ of the projected points is shown in \textcolor{green}{green}.}
\label{fig:Samples}
\end{figure}


\end{example}

\subsection{The number of vertices of $\beta$-polygons in the plane}\label{ssec:CLT}

Let $d\geq 4$ and let $\polytopeQ_n$ be the  convex hull of the orthogonal projections to the disk $\B$ of  $n$ i.i.d. random points on the sphere $\mathbb{S}^{d-1}$.
We have already seen in the previous section that this gives rise to $n$ points in the plane, distributed according to the density function $f_{2, \frac{d}{2} - 2}$ .

Note that $\polytopeQ_n$ is a random variable.
Our first aim is to determine the asymptotics of the expectation and the variance of $f_0(\polytopeQ_n)$ if $n$ tends to $+\infty$. 
 

In what follows, a \defn{constant} will refer to a positive non-zero number that only depends on  $\beta$ or  $d$, rather than anything else, such as the number of vertices $n$, a small value $\varepsilon$ or a radius $R$. 
Apart from the constant $c_0$ defined in $\varepsilon = c_0\frac{\log n}{n}$ (see \Cref{sssec:Variance}), all other constants are of no importance and, slightly abusing notation, will simply write $c$ for a constant. 
In particular, this means that the exact value of $c$ might change from one line to the next.

If $\beta$ is clear from the context, we use \defn{$\mu(A)$} to denote the measure of a set $A\subseteq \B$ with respect to the density $f_{2, \beta}(\b x)$; see \Cref{ssec:ProbabilisticModel}.
Throughout this section, we assume: $\beta = \frac{d}{2}-2$.

All limits, equivalences, approximations, are taken assuming $n \to +\infty$, $\varepsilon\to 0$ or $R\to 1$.

\subsubsection{Expectation}\label{sssec:Expectation}

The expectation of $f_0(\polytopeQ_n)$ and its asymptotics directly follows from several papers on $\beta$-polytopes, \eg \cite[Thm. 1.2, Rmk. 1.4]{KabluchkoThaeleZaporozhets:Beta}.
More precisely, \cite[Thm. 1.8]{KabluchkoThaeleZaporozhets:Beta} yields:

\begin{corollary}[{From \cite[Thm. 1.8]{KabluchkoThaeleZaporozhets:Beta}}] \label{cor:expectation}
Let $d\geq 5$.
The expected number of vertices of $\polytopeQ_n$ is:
$$\E\bigl(f_0(\polytopeQ_n)\bigr) \sim c\, n^{\frac{1}{d-1}}.$$
\end{corollary}

\begin{remark}
In his famous book \cite{Borgwardt1987-SimplexMethod}, Borgwardt showed that \Cref{cor:expectation} also holds in the dual case.
More precisely, for $\b a_1, \dots, \b a_m \in \R^d$ i.i.d. random vectors taken uniformly on the sphere, he considers the polytope $\overline{\polytopeQ}_m\eqdef \bigl\{\b x\in \R^d ~;~ \inner{\b x, \b a_i} \leq 1 \text{ for all } i \bigr\}$.
He shows that the expected number of steps needed by the simplex method on $\overline{\polytopeQ}_m$ with the shadow vertex rule (\ie the expected length of a coherent path on $\overline{\polytopeQ}_m$) can be lower- and upper-bounded by $c\, m^{\frac{1}{d-1}}$ and $c'\, d^2 \,m^{\frac{1}{d-1}}$, respectively, for two constants $c, c' > 0$ (independent of both $m$ and $d$). (See \cite[0.5.7 \& 0.5.8]{Borgwardt1987-SimplexMethod}. For disambiguation: Borgwardt's $n$ is our $d$; his $E_{m, n}(X)$ is the expectation of $X$ over all instances of his probabilistic model (\eg random $\b a_1, \dots, \b a_m$ on $\mathbb{S}^{n-1}$) with dimension $n$ and $m$ facets; his~$s$ is the number of steps of the simplex method which is our length of a path; his $S$ is a good proxy for his $s$.) 

On the other hand, Kelly \& Tolle proved in \cite{KellyTolle1981-ExpectedNbVerts} that, for fixed dimension $d$, the expected number of vertices of $\overline{\polytopeQ}_m$ is linear in $m$.
Moreover, if $m$ is large, these vertices are ``not far'' from lying on the unit sphere (on purpose, we do not make this ``not far'' precise).

Hence, one may think that the asymptotic behavior of $\E(f_0(\polytopeQ_n))$ could be retrieved from Borgwardt's result $m^{\frac{1}{d-1}}$ by replacing the number of facets $m$ by the number of vertices $n$, since they are proportional to each other according to Kelly \& Tolle.
However, since we did not find a way to make this belief rigorous (especially, because the probabilistic models are not the same, as the uniform distribution on the dual does not have an immediate translation to the primal), we decided to dive into the technical details of $\beta$-distributions instead.
\end{remark}

\subsubsection{Variance}\label{sssec:Variance}

To derive the asymptotics of the variance of $f_0(\polytopeQ_n)$, we will provide lower and upper bounds and show that they match asymptotically. 
The proofs of both bounds will heavily rely on $\varepsilon$-caps.

\paragraph{$\varepsilon$-caps}

Intuitively, when $n$ is large, more  random points $Z_1,\ldots,Z_n$ tend to lie close to the circle $\mathbb{S}^1$ and, consequently, so will the vertices of $\polytopeQ_n$. 
$\varepsilon$-caps can be used to make this intuition mathematically rigorous: the reader should think of an $\varepsilon$-cap as a region of the disk $\B$, that is  close to $\mathbb{S}^1$ and  small enough to ensure that it only contains only a small number of vertices of $\polytopeQ_n$.

\begin{figure}
\centering
\includegraphics[width=0.28\linewidth]{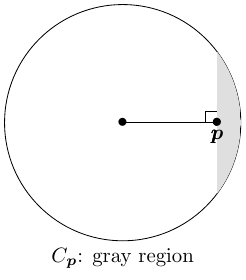}
\includegraphics[width=0.32\linewidth]{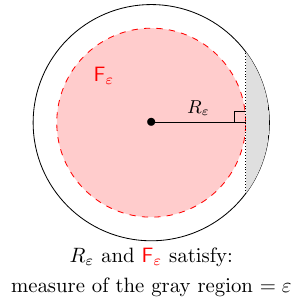}
\includegraphics[width=0.28\linewidth]{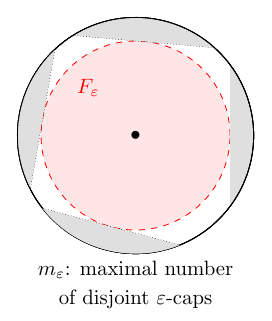}
\caption[$\varepsilon$-caps and independent $\varepsilon$-caps]{(Left) Cap induced by $\b p \in \B$.
(Middle) The $\varepsilon$-floating body is the complement of the union of all $\varepsilon$-caps.
(Right) The maximal number of independent $\varepsilon$-caps is denoted by $m_\varepsilon$.}
\label{fig:esp-capAndM-eps}
\end{figure}

\begin{definition}
Let $\b p\in \B$. The set \defn{$\polytopeC_{\b p}$} $\eqdef \{\b x\in \B ~;~ \inner{\b x - \b p,\, \b p} \geq 0\}$ is called the \defn{cap} induced by $\b p$ (see \Cref{fig:esp-capAndM-eps} Left). 
We call $\|\b p\|$  the  \defn{radius} of $\polytopeC_{\b p}$ and a cap $\polytopeC_{\b p}$ is called  an \defn{$\varepsilon$-cap} if $\mu(\polytopeC_{\b p}) = \varepsilon$. (As before, $\mu$ denotes the measure for the $\beta$-distribution with $\beta=\frac{d}{2}-2$.)
\end{definition}

\begin{lemma}\label{lem:FloatingRadius}
Let $d\geq 4$ and let $\polytopeC$ be an $\varepsilon$-cap. If $\varepsilon\to 0$, then the radius $R_\varepsilon$ of $\polytopeC$ satisfies the following:
$$1 - R_\varepsilon \sim c\, \varepsilon^{\frac{2}{d-1}}.$$
\end{lemma}

\begin{proof}
Let $\polytopeC$ be a cap of radius $R$. 
We first compute $\mu(\polytopeC)$ and then invert the formula to get an estimate for $R_\varepsilon$. (Note that this is possible since, due to rotational symmetry, $\mu(\polytopeC)$ depends only on $R$.) We have

\begin{minipage}[c]{0.55\linewidth}
\begin{align*}
\mu(R) &= C\, \int_{x=R}^1\int_{y = -\sqrt{1-x^2}}^{+\sqrt{1-x^2}} \left(1-x^2-y^2\right)^\beta \, \mathrm{d}y\,\mathrm{d}x \\
&= C\, \int_R^1 \int_{\theta = -\arccos{r}}^{+\arccos{r}} (1-r^2)^\beta \, r\mathrm{d}r\,\mathrm{d}\theta \\
&= 2\, C\, \int_R^1 r\,(1-r^2)^\beta\, \arccos{r} \, \mathrm{d}r.
\end{align*}
\end{minipage} 
\hspace{1cm}\begin{minipage}[r]{0.4\linewidth}
\begin{flushleft}
\includegraphics[width=0.45\linewidth]{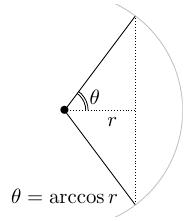}
\end{flushleft}
\end{minipage}
Note that  $\arccos r = \sqrt{2}(1-r)^\frac{1}{2} + O\left((1-r)^\frac{3}{2}\right)$, if $r\to 1$.
Moreover,  we have $R(1+R)^\beta \cdot (1-r)^\beta \leq r(1-r^2)^\beta \leq 2^\beta\cdot (1-r)^\beta$ for all $r\in [R, 1]$.
We conclude

\begin{scriptsize}
$$\begin{array}{rcl}
\displaystyle 2^{\frac{3}{2}} C\, R(1+R)^\beta \int_R^1 \left((1-r)^{\beta+\frac{1}{2}}\, + O\left((1-r)^{\beta + \frac{3}{2}}\right) \right) \, \mathrm{d}r \!\!\!\!&\leq \mu(R) \leq&\!\!\!\! \displaystyle 2^{\beta + \frac{3}{2}} C\, \int_R^1 \left((1-r)^{\beta+\frac{1}{2}}\, + O\left((1-r)^{\beta + \frac{3}{2}}\right) \right) \, \mathrm{d}r 
\end{array}
$$
\end{scriptsize}
and thus
\begin{scriptsize}
$$
\begin{array}{rcl}
\displaystyle \frac{\sqrt{2}R(1+R)^\beta}{2\beta + 3}\, C\, (1-R)^{\beta+\frac{3}{2}}\, + O\left((1-R)^{\beta + \frac{5}{2}}\right) \!\!\!\!&\leq \mu(R) \leq&\!\!\!\! \displaystyle \frac{2^{\beta+\frac{1}{2}}}{2\beta + 3}\, C\, (1-R)^{\beta+\frac{3}{2}}\, + O\left((1-R)^{\beta + \frac{5}{2}}\right) .
\end{array}$$
\end{scriptsize}
Since the lower and upper bounds coincide if $R \to 1$, we obtain $\mu(R) \sim C'\, (1-R)^{\beta + \frac{3}{2}}$ for some constant $C'$ (independent of $R$).
Finally, as $\varepsilon\to 0$ if and only if $R_\varepsilon\to 1$, it follows from $\mu(R_\varepsilon) = \varepsilon$ that $1-R_\varepsilon = c\, \varepsilon^{\frac{1}{\beta + \frac{3}{2}}}$.
Substituting $\beta = \frac{d}{2}-2$ completes the proof.
\end{proof}

In the following, we will need to work with several independent $\varepsilon$-caps.

\begin{definition}
Two caps $\polytopeC, \polytopeC'$ are said to be \defn{independent} if they are disjoint \ie $\polytopeC\cap\polytopeC' = \emptyset$.
We denote by \defn{$m_\varepsilon$} the maximal number of pairwise independent $\varepsilon$-caps.
\end{definition}

\begin{proposition}\label{prop:m_espsilon}
Let $d\geq 4$. 
If $\varepsilon \to 0$, then we have $m_\varepsilon \sim c\, \left(\frac{1}{\varepsilon}\right)^{\frac{1}{d-1}}$.
\end{proposition}

\begin{proof}
Let $\alpha_\varepsilon$ be the half-angle spanned by an $\epsilon$-cap, see picture below-right.

\begin{minipage}[c]{0.61\linewidth}
Then $m_\varepsilon \sim \frac{2\pi}{2\alpha_\varepsilon}$, and $\alpha_\varepsilon = \arccos{R_\varepsilon}$.
Using \Cref{lem:FloatingRadius}, and $\arccos{x} = \sqrt{2(1-x)} + o\left(1-x\right)$ if $x\to 1$, we get:
$$m_\varepsilon ~\sim~ \pi\, \left(\arccos{R_\varepsilon}\right)^{-1}
~\sim~ \frac{\pi}{\sqrt{2}}\, \left(\sqrt{1 - R_\varepsilon}\right)^{-1}
~\sim~ \frac{\pi\, c}{\sqrt{2}}\, \left(\frac{1}{\varepsilon}\right)^{\frac{1}{d-1}}.\qedhere$$
\end{minipage} 
\hspace{1cm}\begin{minipage}[r]{0.3\linewidth}
\begin{center}
\includegraphics[width=0.6\linewidth]{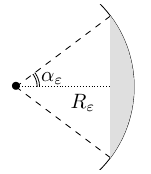}
\end{center}
\end{minipage}
\end{proof}

\paragraph{Inside an $\varepsilon$-cap}

\begin{figure}
\centering
\includegraphics[width=0.495\linewidth]{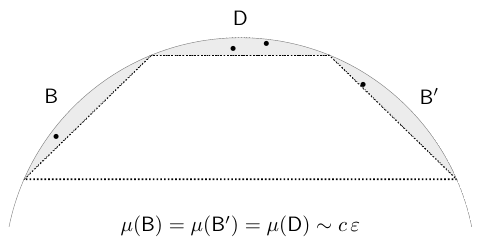}
\includegraphics[width=0.495\linewidth]{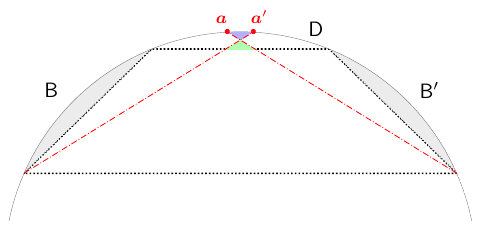}\includegraphics[width=0.495\linewidth]{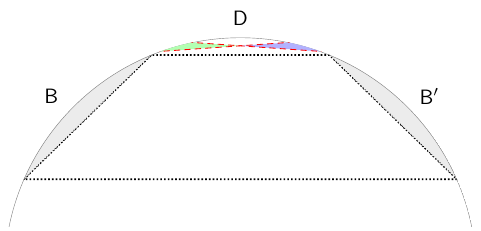}
\caption[The event $A_\polytopeC$ and its sub-events]{(Top) Definition of the subsets $\polytopeB$, $\polytopeB'$ and $\polytopeD$, and event $A_\polytopeC$.
(Bottom) Two sub-events with different $f_0(\polytopeQ_n)$ that occur with strictly positive probability (see the proof of \Cref{prop:ProbaA}).}
\label{fig:InsideCapVariance}
\end{figure}

To lower bound the variance of $f_0(\polytopeQ_n)$, we will consider certain local events (one for each $\varepsilon$-cap) that occur with strictly positive probability $p$, and are independent.
We then show that, up to a constant, the variance is lower bounded by $p\cdot m_\varepsilon$.
We now make this strategy precise.

\begin{definition}
Let $\polytopeC$ be an $\varepsilon$-cap. Let \defn{$\polytopeB$}, \defn{$\polytopeD$} and \defn{$\polytopeB'$} be the three maximal independent sub-caps contained in $\polytopeC$ labeled from left to right (see \Cref{fig:InsideCapVariance} Top). 
Let \defn{$A_\polytopeC$} be the event that among $X_1, \dots, X_n$ there is exactly one point in $\polytopeB$, one point in $\polytopeB'$, two points lie in $\polytopeD$, and all the other points lie outside the cap $\polytopeC$.
\end{definition}
Note that $\polytopeB$, $\polytopeD$ and $\polytopeB'$ have the same measure and that the sets   $\B\ssm\polytopeB$, $\B\ssm\polytopeB'$ and $\B\ssm\polytopeD$ are convex. 
We will show that $\mu(\polytopeB) = \mu(\polytopeB') = \mu(\polytopeD) \simeq c\, \varepsilon$ for some constant $c > 0$. 
We compute $\P(A_\polytopeC)$ and $\Var \bigl(f_0(\polytopeQ_n) \,\big|\, A_\polytopeC\bigr)$. The latter is the variance of $f_0(\polytopeQ_n)$, conditioned on the event $A_\polytopeC$, where all $X_i$ except the four ones that  achieve the event $A_\polytopeC$ are fixed.

\begin{proposition}\label{prop:ProbaA}
Let $d\geq 4$, $c_0 > 0$ and $\varepsilon = c_0\frac{\log n}{n}$.  
For an $\varepsilon$-cap $\polytopeC$, we have
\[\P(A_\polytopeC) \geq c\, (\log n)^4 \, n^{-c_0}.\]
\end{proposition}

\begin{proof}
First, we compute $\mu(\polytopeB) = \mu(\polytopeB') = \mu(\polytopeD)$.
The half-angle $\theta_\varepsilon$ spanned by $\polytopeB$ is $\frac{1}{3}$ of the half-angle $\alpha_\varepsilon$ spanned by $\polytopeC$ (see the proof of \Cref{prop:m_espsilon}).
It follows that $1 - \cos \theta_\varepsilon \sim 1 - \cos\frac{\alpha_\varepsilon}{3} \sim \frac{1}{3^2}\frac{\alpha_\varepsilon^2}{2} \sim \frac{1}{9}(1-\cos\alpha_\varepsilon)$, where $\cos\alpha_\varepsilon$ and $\cos\theta_\varepsilon$ are the radius of the cap $\polytopeC$ and $\polytopeB$, respectively.
Since, by \Cref{lem:FloatingRadius}, the measure of a cap of radius $R$ is equal to $c\, (1 - R^\frac{2}{d-1})$ for some $c > 0$, we obtain $\mu(\polytopeB) \sim c\, 9^{-\frac{2}{d-1}}(1 - \cos\alpha_\varepsilon)^{\frac{2}{d-1}} \sim 9^{-\frac{2}{d-1}}\mu(\polytopeC) = 9^{-\frac{2}{d-1}}\, \varepsilon$.

As all $X_i$ are independent, we have:
\begin{align*}
\P(A_\polytopeC) &= \sum_{\{i, j, k, \ell\} \subseteq [n]} \P(X_i\in \polytopeB)\P(X_j\in \polytopeB')\P(X_k\in \polytopeD)\P(X_\ell\in \polytopeD)\prod_{g\notin\{i, j, k, \ell\}}\bigl(1 - \P(X_g \in\polytopeC)\bigr) \\
&= \binom{n}{4} \bigl(\mu(\polytopeB)\bigr)^4 (1-\varepsilon)^{n-4}
~~\sim \frac{c_0^4}{4! \cdot 9^{\frac{8}{d-1}}}\, (\log n)^4 \, n^{-c_0}.
\end{align*}
For the last estimate, we have used that $\varepsilon \to 0$ implies that $\log(1-\varepsilon) = -\varepsilon + o(\varepsilon)$, which yields $(1-\varepsilon)^n = \exp\left(n\cdot \left(-c_0 \frac{\log n}{n} + o\left(\frac{\log n}{n}\right) \right)\right) \sim n^{-c_0}$.
\end{proof}

We can now bound the conditional variance
$\Var \bigl(f_0(\polytopeQ_n) \,\big|\, A_\polytopeC\bigr)$  away from $0$ independently of $\varepsilon$:
 
\begin{proposition}\label{prop:PositiveVarianceForEventA}
Let $d\geq 4$.  There exists a constant $c > 0$ such that $\Var \bigl(f_0(\polytopeQ_n) \,\big|\, A_\polytopeC\bigr) \geq c$ for all $\varepsilon>0$ sufficiently small and for every $\varepsilon$-cap $\polytopeC$.
\end{proposition}


\begin{proof}
Suppose that $X_1, \dots, X_n$ achieve the event $A_\polytopeC$.
Let $\b b, \b b'$ and $\b d_1, \b d_2$ be the points among $X_1, \dots, X_n$ which lie in $\polytopeB$, $\polytopeB'$ and $\polytopeD$, respectively.
As $\B\ssm\polytopeB$ is convex, $\b b$ is a vertex of $\polytopeQ_n = \conv(X_i ; i\in [n])$.
For the same reason, $\b b'$ and at least one of $\b d_1, \b d_2$ are vertices of $\polytopeQ_n$.

Let $v$ be the number of vertices of $\polytopeQ_n$ that are different from $\b b$, $\b b'$, $\b d_1$ and $\b d_2$.
By the previous paragraph, we have $f_0(\polytopeQ_n) \in \{v + 3,v+4\}$, depending on whether $\conv(\b b, \b b', \b d_1, \b d_2)$ is a triangle or a quadrangle.
The idea is to show that the probabilities $\P\bigl(\conv(\b b, \b b', \b d_1, \b d_2) \text{ is a triangle}\bigr)$ and $\P\bigl(\conv(\b b, \b b', \b d_1, \b d_2) \text{ is a quadrangle}\bigr)$ are strictly positive.
It will follow that $f_0(\polytopeQ_n)$ is not only determined by $v$, which implies that $\Var \bigl(f_0(\polytopeQ_n) \,\big|\, A_\polytopeC\bigr)$ is strictly positive.

Consider two lines passing through the center of $\polytopeD$ and separating $\polytopeB$ from $\polytopeB'$ (see \textcolor{red}{red dashed} lines in \Cref{fig:InsideCapVariance}, Bottom left): 
These lines divide $\polytopeD$ into four subsets and two of these are separated from $\polytopeB$ and $\polytopeB'$ by these lines. 
We will call these subsets the ``top region'' (\textcolor{blue}{blue}) and the ``bottom region'' (\textcolor{green}{green}). 
If $\b d_1$ lies in the top region and $\b d_2$ in the bottom region, then $\b d_2 \in \conv(\b b, \b b', \b d_1)$ (for any $\b b\in \polytopeB$, $\b b'\in \polytopeB'$). Due to convexity, this can be verified by considering the extreme points of each of the four regions (yielding $3 \cdot 3\cdot 2\cdot 2 $ cases) and checking that any extreme point of the green region is contained in the triangle formed by any collection of extreme points of the other three regions. 
This implies that $\conv(\b b, \b b', \b d_1, \b d_2)$ is always a triangle in this case.
\begin{minipage}[c]{0.55\linewidth}
Let $R\eqdef R_\varepsilon$ and $R_\polytopeD$ be the radius of the cap $\polytopeC$ and $\polytopeD$, respectively. Moreover, set $ 1 - H=\frac{1}{2}(1 - R_\polytopeD)$. A computation shows (see figure on the right) that the radius of the cap lying above the line between $\b a$ and $\b a'$ inside the \textcolor{blue}{top region} (see \Cref{fig:InsideCapVariance}, bottom left) is equal to $y = \frac{2H - RH^2 -R}{1 + H^2 - 2RH}$.
\end{minipage} 
\hspace{1cm}\begin{minipage}[r]{0.4\linewidth}
\begin{flushleft}
\includegraphics[width=0.9\linewidth]{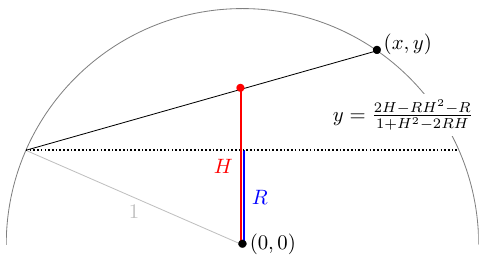}
\end{flushleft}
\end{minipage}

As seen in the proof of \Cref{prop:ProbaA}, we have $1-R_\polytopeD \sim \frac{1}{9}(1 - R)$.
Further computations show $1-y \sim \frac{1 + R}{35 R + 613}(1 - R)$ and hence $\mu(\text{top region}) \sim \frac{1}{324^{\frac{d-1}{2}}}\,\varepsilon$ when $\varepsilon \to 0$. A similar reasoning applies to the \textcolor{green}{bottom region}, showing that $\mu(\text{bottom region}) \sim c'\, \varepsilon$ for a certain constant $c' > 0$.
Thus, the probability that $\conv(\b b, \b b', \b d_1, \b d_2)$ is a triangle is at least $\frac{\mu(\text{top region})\mu(\text{bottom region})}{\mu(\polytopeD)^2}$ which is lower bounded by some strictly positive constant, independently of $\varepsilon$.

Similarly, if $\b d_1$ is in the right region (\textcolor{blue}{blue}) and $\b d_2$ in the left region (\textcolor{green}{green}) of \Cref{fig:InsideCapVariance} (Bottom Right), then $\conv(\b b, \b b', \b d_1, \b d_2)$ is a quadrangle (for any $\b b\in \polytopeB$, $\b b'\in \polytopeB'$).
Similarly as above, this occurs with a positive probability that can be bounded away from $0$ (independently of $\varepsilon$).

Consequently, $\Var \bigl(f_0(\polytopeQ_n) \,\big|\, A_\polytopeC\bigr)$ is lower bounded by a positive constant.
\end{proof}

\paragraph{Lower bound on the variance}

We lower bound the variance with \Cref{prop:m_espsilon,prop:ProbaA,prop:PositiveVarianceForEventA}.
Remark that, unlike \Cref{prop:FloatingBody}, the lower bound in this corollary holds for \emph{all $c_0 > 0$}. 

\begin{corollary}\label{cor:VarianceLower}
Let $d\geq 4$. 
For any $c_0 > 0$, there exists a constant $c > 0$ such that:
$$\Var f_0(\polytopeQ_n) \geq c \, n^{\frac{1}{d-1} - c_0}\quad \text{if } n\to+\infty. $$
\end{corollary}

\begin{proof}
Let $\epsilon=c_0\frac{\log(n)}{n}$ and let $\c C = (\polytopeC_1, \dots, \polytopeC_{m_\varepsilon})$ be a collection of $m_\varepsilon$ independent $\varepsilon$-caps. 

For $\beta$-distributed points $X_1, \dots, X_n$, we let $\b X \eqdef \bigl(X_i ~;~ X_i \notin \polytopeC \text{ for all } \polytopeC\in \c C~\text{with } \b1(A_\polytopeC) =1\bigr)$, where $\b 1(A)$ is the indicator function of the event~$A$.
Intuitively, $\b X$ is the random variable consisting of the points of $(X_1, \dots, X_n)$ that are not involved in any occurring event $A_\polytopeC$ (for $\polytopeC\in \c C$).

We recall the law of total variance: if $Y$ and $Z$ are random variables (with $\Var Y < + \infty$), then $\Var Y = \E\bigl(\Var (Y ~\big|~ Z)\bigr) + \Var\bigl(\E (Y ~\big|~ Z)\bigr)$, where on the right-hand side the expectation and variance are conditioned on $Z$.
In particular: $\Var Y \geq \E\bigl(\Var (Y ~\big|~ Z)\bigr)$.
Applying this inequality to $Y = f_0(\polytopeQ_n)$ and $Z = \b X$ yields 
$\Var f_0(\polytopeQ_n) \geq \E\bigl(\Var (f_0(\polytopeQ_n) ~\big|~ \b X)\bigr)$.

To compute $\Var (f_0(\polytopeQ_n) ~\big|~ \b X)$, first note that for independent caps $\polytopeC, \polytopeC'\in \c C$ the events $A_{\polytopeC}$ and~$A_{\polytopeC'}$ are independent. 
Moreover, if $A_\polytopeC$ and $A_{\polytopeC'}$ are events with $\b1(A_\polytopeC)=\b1(A_{\polytopeC'})=1$, then moving the points that are witnesses for $A_\polytopeC$ does not affect which of the points that witness $A_{\polytopeC'}$ are vertices of $\polytopeQ_n$ and vice versa (as, for each cap $\polytopeC$, the sub-sets $\B\ssm\polytopeB$, $\B\ssm\polytopeB'$ and $\B\ssm\polytopeD$ are convex). 
Using this independence structure and \Cref{prop:PositiveVarianceForEventA}, we conclude that there exists $c > 0$ with:
$$\Var (f_0(\polytopeQ_n) ~\big|~ \b X) ~\geq~ \sum_{\polytopeC\in \c C} \Var(f_0(\polytopeQ_n) ~\big|~ A_\polytopeC)\, \b1(A_\polytopeC) ~\geq~ c \, \sum_{\polytopeC\in\c C} \b 1(A_\polytopeC).$$

Finally, by \Cref{prop:m_espsilon,prop:ProbaA}, using $\varepsilon = c_0\frac{\log n}{n}$, we see that there exists $c' > 0$ with:
$$\begin{array}{rclcr}
\Var f_0(\polytopeQ_n) &\geq& \E\bigl(\Var (f_0(\polytopeQ_n) ~\big|~ \b X)\bigr) && \\
&\geq& c\, \sum_{\polytopeC\in\c C} \P(A_\polytopeC) ~\geq~ c\, m_\varepsilon\, \P(A_{\polytopeC_1}) && \\
&\geq& c'\, \left(\frac{n}{\log n}\right)^{\frac{1}{d-1}}\,\, (\log n)^4\, n^{-c_0}. &&
\end{array}$$
As $(\log n)^{4 - \frac{1}{d-1}} > 1$ for large $n$, we can remove this term to get the claimed formula.
\end{proof}


\paragraph{$\varepsilon$-floating body and $\varepsilon$-visible region}
In order to give an accurate upper bound for the variance of~$f_0(\polytopeQ_n)$, we will need to understand what happens when we ``add a point'' to $\polytopeQ_n$, and use the so-called \emph{Efron-Stein jackknife inequality} (see the next paragraph).
Firstly, we want to measure how close the vertices of $\polytopeQ_n$ are to the boundary of the disk: To this end, we use the $\varepsilon$-floating body.

\begin{figure}
\centering
\includegraphics[width=0.28\linewidth]{Figures/Floating_body/cap_definition.pdf}
\includegraphics[width=0.32\linewidth]{Figures/Floating_body/cap.pdf}
\includegraphics[width=0.3\linewidth]{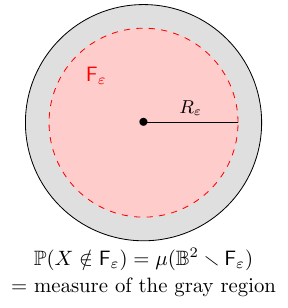}
\caption[$\varepsilon$-caps and the $\varepsilon$-floating body]{(Left) Cap induced by $\b p \in \B$.
(Middle) The $\varepsilon$-floating body is the complement of all~$\varepsilon$-caps.
(Right) With high probability, all the vertices and edges of $\polytopeQ_n$ are in the gray region.}
\label{fig:CapStuffs}
\end{figure}

\begin{definition}
Let $\varepsilon > 0$. 
The \defn{$\varepsilon$-floating body}, see \Cref{fig:CapStuffs}, is the complement of all $\varepsilon$-caps, \ie \defn{$\FBody$} $\eqdef \B \ssm \bigcup_{\polytopeC ~\varepsilon\text{-cap}} \polytopeC = \{\b p\in \B ~;~ \mu(\polytopeC_{\b p}) > \varepsilon\}$.
\end{definition}

The reader should think of it as a disk that is contained in $\polytopeQ_n$ with very high probability (for $n$ large enough): this helps us rule out cases where $\polytopeQ_n$ is ``far'' from the circle $\mathbb{S}^1$.
To be precise, suppose that $\FBody\subseteq \polytopeQ_n$, then all vertices of $\polytopeQ_n$ are in $\varepsilon$-caps.
Reciprocally, the number of vertices of $\polytopeQ_n$ within any $\varepsilon$-cap is at least 1, since $\FBody$ is not contained in the interior of~$\B\ssm\polytopeC$. 
The next lemmas will show that roughly $\frac{\mu(\polytopeC)}{\mu(\B\ssm\FBody)}\,\E \bigl(f_0(\pol_n)\bigr) \sim c\, (\log n)^{\frac{1}{d-1}}$.

\begin{remark}
The name \emph{floating} body comes from the following idea: Consider an object of arbitrary shape (here a disk) that is made out of a material with high buoyancy (\eg foam). Next, immerse this object in water and  roll it until every part that can become wet is wet.
The part of the object immersed at a given moment is the cap, and the part that remains forever dry is the floating body.
\end{remark}

\begin{lemma}
Let $d\geq 4$. The $\varepsilon$-floating body is a disk of radius $R_\varepsilon$. If  $\varepsilon\to 0$, then $1 - R_\varepsilon \sim c\, \varepsilon^{\frac{2}{d-1}}$.
\end{lemma}

\begin{proof}
By rotational symmetry of the $\beta$-distribution, the floating body is a (possibly empty) disk (for all $\varepsilon$).
If $\b p$ is on the boundary of $\FBody$, then $\mu(\polytopeC_{\b p}) = \varepsilon$ and \Cref{lem:FloatingRadius} implies the claim.
\end{proof}

\begin{lemma}\label{lem:FloatingOutside}
Let $d\geq 4$. 
 The measure of the region outside the $\varepsilon$-floating body satisfies the following:
$$\mu(\B \ssm \FBody) \sim c\, \varepsilon^{1 - \frac{1}{d-1}}\quad \text{if } \varepsilon \to 0.$$
\end{lemma}

\begin{proof}
By \Cref{lem:FloatingRadius}, the region $\B\ssm\FBody$ is an annulus of inner radius $R_\varepsilon$ and outer radius $1$.
Hence, if $\varepsilon\to 0$, its measure according to the $\beta$-distribution is equal to:
$$\mu(\B\ssm\FBody) = C\, \int_{R_\varepsilon}^1 2\pi\, (1-r^2)^\beta \, r\mathrm{d}r ~~=~ \frac{\pi\, C}{\beta+1}\, (1-R_\varepsilon^2)^{\beta+1} ~~\sim~ \frac{2^{\beta+1}\pi\, C}{\beta + 1}\, (1-R_\varepsilon)^{\beta+1}.$$
Using $\beta = \frac{d}{2} - 2$ and \Cref{lem:FloatingRadius}, we get the claimed formula.
\end{proof}

Last but not least, we need to introduce the notion of visibility from a point.

\begin{figure}
\centering
\includegraphics[width=0.27\linewidth]{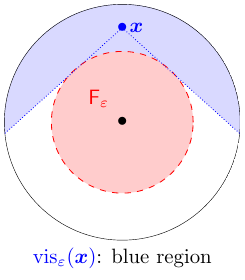}
\includegraphics[width=0.27\linewidth]{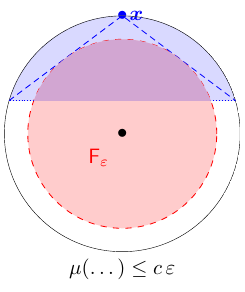}
\includegraphics[width=0.44\linewidth]{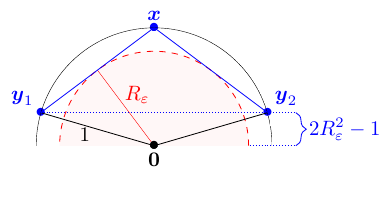}
\caption[Visibility region and the upper bound of its measure]{(Left) The $\varepsilon$-visible region from $\b x$.
(Middle) The biggest visibility region is achieved when $\b x\in\mathbb{S}^1$. It is covered by the blue-shaded cap.
(Right) The radius of this cap is $2R_\varepsilon^2-1$.}
\label{fig:VisibilityStuffs}
\end{figure}

\begin{definition}
Let $\varepsilon > 0$.
Let $\b x\in \B$. The \defn{$\varepsilon$-visible region from $\b x$} is the following subset of the disk: \defn{$\vis \b x$} $\eqdef \{\b y\in \B ~;~ [\b x, \b y]\cap\FBody = \emptyset\}$, see \Cref{fig:VisibilityStuffs} (Left).
\end{definition}

\begin{lemma}\label{lem:VisibilityProba}
For $\b x\in \B$, we have
$\mu\bigl(\vis\b x\bigr) \sim c\, \varepsilon$, if $\varepsilon \to 0$.
\end{lemma}

\begin{proof}
The measure of $\vis \b x$ is maximized when $\b x\in \mathbb{S}^1$, see \Cref{fig:VisibilityStuffs} (Middle).
In this case, $\vis \b x$ is included in the cap $\polytopeC$ obtained by taking the two tangents to the disk $\FBody$ passing through~$\b x$, and joining their points of intersection with $\mathbb{S}^1$, namely $\b y_1, \b y_2$ in \Cref{fig:VisibilityStuffs} (Right).
Looking at the kite defined by $\b 0, \b x, \b y_1$ and $\b y_2$, a simple computation shows that the radius defining this cap is $2R_\varepsilon^2-1$ if $\varepsilon$ is small enough, see \Cref{fig:VisibilityStuffs} (Right).
According to the proof of \Cref{lem:FloatingRadius}, we get that $\mu(\polytopeC)$ is of order $\bigl(1-(2R_\varepsilon^2-1)\bigr)^{\beta+\frac{3}{2}}$.
As $R_\varepsilon \to 1$, we get $1-R_\varepsilon^2 = (1+R_\varepsilon)(1-R_\varepsilon) \sim 2(1-R_\varepsilon)$.
Thus: $\mu(\polytopeC) \sim c\, \bigl(2(1-R_\varepsilon^2)\bigr)^{\beta+\frac{3}{2}} \sim 4^{\beta+\frac{3}{2}}\, c\, (1-R_\varepsilon)^{\beta+\frac{3}{2}} \sim c'\, \varepsilon$, where the last equivalence is \Cref{lem:FloatingRadius}.
\end{proof}

We now show that, for large $n$ and $\varepsilon > 0$ \emph{not too small}, the $\varepsilon$-floating body is contained in $\polytopeQ_n$ with very high probability.

\begin{proposition}\label{prop:FloatingBody}
Let $d\geq 4$, $s>0$, $c_0 = \frac{1}{d-1} + s$ and let $\varepsilon = c_0\frac{\log n}{n}$. We have
$$\P(\FBody \not\subseteq \polytopeQ_n) \leq c\, n^{-s} ~~~~\text{ when }~~~~ n\to +\infty.$$
\end{proposition}

\begin{proof}
Let $X_1, \dots, X_n$ be i.i.d. according to the density $f_{2, \frac{d}{2} - 2}$, and let $\polytopeQ_n = \conv(X_1, \dots, X_n)$.

\noindent
\begin{minipage}[c]{0.6\linewidth}
If $\FBody\not\subseteq\polytopeQ_n$, then there are two possibilities:

$(a)$ $\polytopeQ_n\subseteq\FBody$;

$(b)$ there exists $X_i\notin \FBody$ and an edge $\pol[e]$ incident to $X_i$ that  intersects $\FBody$.

Case $(a)$ corresponds to $X_j\in \FBody$ for all~$j\in [n]$.
For case $(b)$, let $\pol[C]_i$ be the $\varepsilon$-cap whose boundary line is parallel to $\pol[e]$ (see figure on the right). Since $\pol[e]$ is an edge, we have $X_j\notin \pol[C]_i$ for all $j\in [n]\setminus \{i\}$.

As $X_1, \dots, X_n$ are independent, we get
\end{minipage} 
\hspace{-0cm}\begin{minipage}[c]{0.38\linewidth}
\begin{center}
\includegraphics[width=0.5\linewidth]{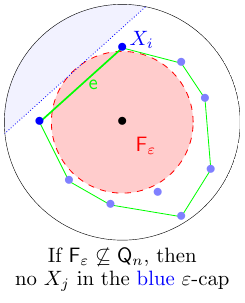}
\end{center}
\end{minipage}

$$\begin{array}{rclcl}
\P(\FBody\not\subseteq\polytopeQ_n) & = & \P(\polytopeQ_n\subseteq\FBody) & + & \P\bigl(\exists i\in [n],\, X_i\notin \FBody ~\text{ and }~\forall j\neq i,\, X_j\notin \pol[C]_i \bigr) \\
 & \leq & \prod_{i=1}^n \P(X_i\in\FBody) & + & \sum_{i = 1}^n \P(X_i\notin\FBody)\,\prod_{j\in [n]\setminus\{i\}}\P(X_j \notin \pol[C]_i).
\end{array}$$

For the asymptotics, when $n\to +\infty$, we use the fact that by \Cref{lem:FloatingOutside} there exists a constant $c_1$ such that $\P(X_i \notin \FBody) \leq c_1\, \varepsilon^{1 - \frac{1}{d-1}}$. Since $\P(X_j \notin\pol[C]_i) = 1-\mu(\pol[C]_i) = 1-\varepsilon$.
We get
$$\begin{array}{rclcl}
\P(\FBody\not\subseteq\polytopeQ_n) & \leq & \left(1 - c_1\, \varepsilon^{1 - \frac{1}{d-1}}\right)^n & + & n\, \left(c_1\, \varepsilon^{1 - \frac{1}{d-1}}\right)\, \left(1 - \varepsilon\right)^{n-1}.
\end{array}$$

Using $\varepsilon = c_0\frac{\log n}{n}$ and $\log(1-x) \sim -x$ when $x \to 0$, for the first term on the right side, we get
$\left(1 - c_1\, \varepsilon^{1 - \frac{1}{d-1}}\right)^n ~\sim~ \exp\left(-c_1 c_0^{1 - \frac{1}{d-1}} \, (\log n)^{1 - \frac{1}{d-1}} \, n^{\frac{1}{d-1}}\right)$.
And, for the second term on the right-hand side, we get:
$n\, \left(c_1\, \varepsilon^{1 - \frac{1}{d-1}}\right)\, \left(1 - \varepsilon\right)^{n-1} ~\sim~ c_1\, (\log n)^{1 - \frac{1}{d-1}}\, \exp\left(\left(\frac{1}{d-1} - c_0\right)\log n\right)$.
The second term is easily seen to be asymptotically larger than the first term when $n\to +\infty$.
Hence, simplifying the last exponential, we get
$$\P(\FBody\not\subseteq\polytopeQ_n) \leq c\, (\log n)^{1 - \frac{1}{d-1}}\, n^{\frac{1}{d-1} - c_0}.$$
By choosing $c_0 = \frac{1}{d-1} + s$, we can ensure that  $\P(\FBody\not\subseteq\polytopeQ_n) \leq c\, n^{-s}$ for any fixed $s > 0$.
\end{proof}

\paragraph{First order difference and Efron-Stein jackknife inequality}
We now explain the main tools used to find an upper bound for the variance of $f_0(\polytopeQ_n)$.
To make things more precise, we first need to introduce some further notation.
In order to make these notions reusable in other contexts, we chose to give a quite general presentation. 

\begin{figure}
\centering
\includegraphics[width=0.8\linewidth]{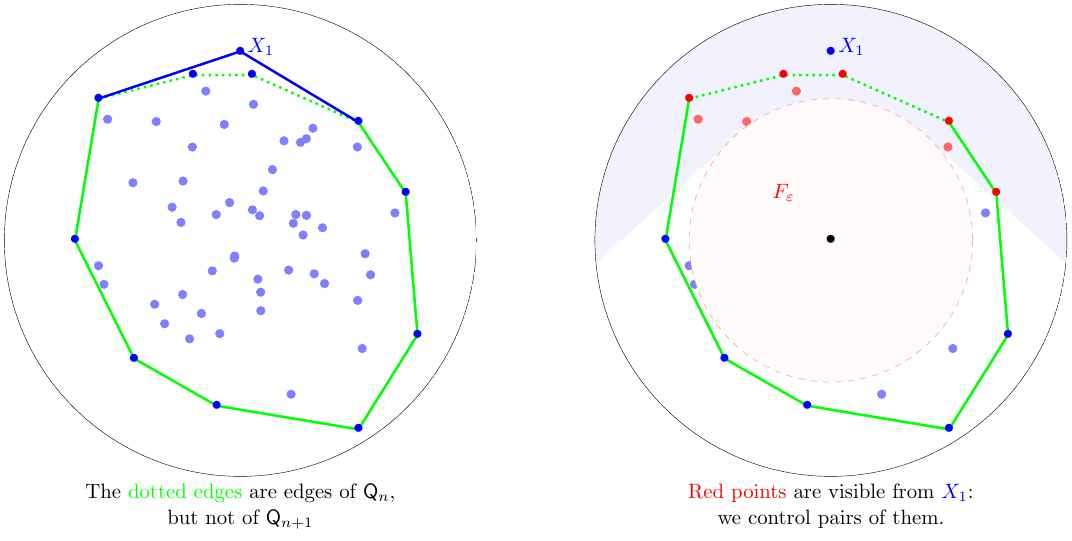}
\caption[First order difference and its upper bound]{(Left) Adding a new point $X_1$ to $\polytopeQ_n$ (\textcolor{green}{green}) creates 2 new edges (\textcolor{blue}{blue}), and destroys 3 old edges (\textcolor{green}{green} dotted), so $D_1 f_1(\polytopeQ_n) = 2 - 3 = -1$ for this sample.
(Right) The points in \textcolor{red}{red} are $\varepsilon$-visible from $X_1$ but not in the $\varepsilon$-floating body $\FBody$.}
\label{fig:D1}
\end{figure}

Let $\X$ be a \emph{Polish space}, for example, $\X$ can be a discrete space or a Euclidean space as $\R^2$.
Let~$f:\bigcup_{n=1}^\infty \X^n\to \R$ be a measurable function on the set of all (finite) ordered point configurations in $\X$.
In our setting, $\X = \B$ and $f = f_0(\polytopeQ_n)$.

The function $f$ is \defn{symmetric} if for any permutation $\sigma\in \c S_n$ and any $(x_1,\ldots,x_n)\in \X^n$, we have 
$f(x_1,\ldots,x_n) = f(x_{\sigma(1)}, \ldots, x_{\sigma(n)})$.
For $i\in [n]$ and $\b x = (x_1, \dots, x_n)\in \X^n$, we let \defn{$\b x^{i}$}~$\eqdef (x_1, \dots, x_{i-1}, x_{i+1}, \dots, x_n)$ be the vector $\b x$ with its $i^{\text{th}}$ coordinate removed.

\begin{definition}\label{def:FirstOrderDifference}
The \defn{first order difference operator} of $f:\bigcup_{n=1}^\infty \X^n\to \R$ with respect to $i$ is 
$$D_i f(\b x) = f(\b x) - f(\b x^i) \quad \text{for }\b x\in \X^n.$$
If the choice of $i$ is irrelevant (\eg if $f$ is symmetric), we write \defn{$D f(\b x)$} instead of $D_i f(\b x)$.
\end{definition}

The idea of the first order difference $Df_0(\polytopeQ_n)$ is to measure the effect of the removal of a point from $n$ randomly chosen points on the number of vertices (or edges) of the convex hull of the random points: how many vertices have been lost or gained? (See \Cref{fig:D1} (left): we count edges in figures, as it is easier to draw.)
The \emph{Efron-Stein jackknife inequality} provides an upper bound for the variance of certain random variables via the just introduced first order difference operators:

\begin{theorem}[{\cite{EfronStein}, \cite[Sec. 3.2]{Reitzner:EfronStein}}]\label{thm:EfronSteinJackknife}
Let $X_1, \dots, X_n, X_{n+1}$ be i.i.d. (according to a measure in a fixed convex body). Let $\polytope{K}_n = \conv(X_1, \dots, X_n)$ and $\polytope{K}_{n+1} = \conv(X_1, \dots, X_n, X_{n+1})$. Then for any real-valued symmetric function $f$, we have
$\Var f(\polytope{K}_n)\leq (n+1)\, \E\Bigl(\bigl(f(\polytope{K}_{n+1})-f(\polytope{K}_n)\bigr)^2\Bigr)$.
\end{theorem}

Applying this inequality to our setting directly yields:
\begin{equation}\label{eq:EfronStein}
\Var f_0(\polytopeQ_n) \leq (n+1)\, \E\Bigl( \bigl(D f_0(\polytopeQ_{n+1})\bigr)^2\Bigr).
\end{equation}

\paragraph{Upper bound on the variance}
Providing an upper bound for the variance (using \Cref{eq:EfronStein}) therefore  amounts to proving an upper bound for the second moment of $D f_0(\polytopeQ_n)$. 
As the proof of the central limit theorem will also require upper bounds for the fourth moment, we will provide upper bounds for all moments of $D f_0(\polytopeQ_n)$.

\begin{theorem}\label{thm:MomentsFirstOrderDifference}
Let $d\geq 4$ and let $p$ be a positive integer.
There exists $c > 0$ such that:
$$\E\bigl(|D f_0(\polytopeQ_n)|^p\bigr) \leq c\, (\log n)^{p + 1 -\frac{1}{d-1}}\, \left(\frac{1}{n}\right)^{1 - \frac{1}{d-1}} \quad \text{if } n\to+\infty.$$
\end{theorem}

\begin{proof}
Since  $\FBody\subseteq \polytopeQ_n$ or $\FBody\not\subseteq \polytopeQ_n$, we can split the expectation as follows:
\begin{equation*}\label{eq:splitMoment}
\E\bigl(|D f_0(\polytopeQ_n)|^p\bigr) = \E\bigl(|D f_0(\polytopeQ_n)|^p \,\big|\, \FBody\not\subseteq \polytopeQ_n\bigr)\P(\FBody\not\subseteq \polytopeQ_n) ~\,+\,~ \E\bigl(|D f_0(\polytopeQ_n)|^p\,\big|\, \FBody\subseteq \polytopeQ_n\bigr) \P(\FBody\subseteq \polytopeQ_n).
\end{equation*}

By \Cref{prop:FloatingBody}, for $s > 0$ and $c_0 \geq \frac{1}{d-1}+s$, there exists a constant $c > 0$ such that $\P(\FBody\not\subseteq \polytopeQ_n) \leq c\, n^{-s}$ for $\varepsilon = c_0 \frac{\log n}{n}$.
In addition, trivially: $|D f_0(\polytopeQ_n)| \leq n$ (it is impossible to gain or lose more than $n$ vertices).
Thus, the first term in the above equation can be bounded as follows: 
\[
\E\bigl(|D f_0(\polytopeQ_n)|^p \,\big|\, \FBody\not\subset \polytopeQ_n\bigr)\P(\FBody\not\subset \polytopeQ_n) \leq c\, n^{p - s}
\]
with $s>0$, $c_0$, $c$ and $\varepsilon$ as above.

Hence, using $\P(\FBody\subseteq \polytopeQ_n) \leq 1$, it remains to prove that $\E\bigl(|D f_0(\polytopeQ_n)|^p\,\big|\, \FBody\subseteq \polytopeQ_n\bigr)$ is upper bounded by the claimed formula.

Suppose that $\FBody\subseteq \polytopeQ_n$.
By symmetry, we only compute the first order difference operator with respect to $X_1$ in~$\polytopeQ_n$, (see \Cref{fig:D1} Left).
We distinguish two cases:
\begin{compactenum}
\item[(a)] $X_1 \in \conv(X_2, \dots, X_n)$: No vertex is removed or added and, hence, $Df_0(\polytopeQ_n) = 0$.
\item[(b)] $X_1 \notin \conv(X_2, \dots, X_n)$: Two new edges appear, and $m$ edges are deleted. An upper bound for $m$ is given by the number of vertices of $\conv(X_2, \dots, X_n)$ in the $\varepsilon$-visible region of $X_1$, see \Cref{fig:D1} Right. This yields 
$|Df_0(\polytopeQ_n)|^p = |m - 2|^p \leq \left(\sum_{i = 2}^n \b 1(X_i \in \vis X_1)\right)^p$ $= \sum_{2\leq i_1, \dots, i_p\leq n} \prod_{j} \b 1(X_{i_j} \in \vis X_1)$.
\end{compactenum}
Note that in the above product certain variables $X_{i_j}$ can be repeated, but since all variables $X_i$ are independent, we get $\P\left(\prod_{j} \b 1(X_{i_j} \in \vis X_1) = 1\right) = \mu(\vis X_1)^{\#\{i_1, \dots, i_p\}}$.
Besides, if  $X_1 \notin\conv(X_2, \dots, X_n)$, then in particular $X_1\notin\FBody$ (since we assumed $\FBody\subseteq\polytopeQ_n$).
Hence, $\P\bigl(X_1 \notin\conv(X_2, \dots, X_n)\bigr) \leq \P(X_1\notin\FBody)$.
Consequently, when $n\to+\infty$
\begin{align*}
\E\bigl(|D f_0(\polytopeQ_n)|^p\,\big|\, \FBody\subseteq \polytopeQ_n\bigr) &\leq 0 ~+~ \P(X_1\notin\FBody)\,\sum_{2\leq i_1, \dots, i_p\leq n} \mu(\vis X_1)^{\#\{i_1, \dots, i_p\}} &\\
&\leq \P(X_1\notin\FBody)\,\sum_{q = 1}^{p} \binom{n-1}{q} q^{p-q} \,\mu(\vis X_1)^{q} &\\
&\leq~ c\, \P(X_1\notin\FBody)\, \sum_{q = 1}^p n^q\, \mu(\vis X_1)^q & \text{for some } c > 0.
\end{align*}

Using \Cref{lem:FloatingOutside}, we have $\P(X_1\notin\FBody) \leq c_1\, \varepsilon^{1-\frac{1}{d-1}}$ for some constant $c_1 > 0$.
By \Cref{lem:VisibilityProba}, we get $\mu(\vis X_1) \leq c_2\, \varepsilon$, for some constant $c_2 > 0$.
So, taking $\varepsilon = c_0\frac{\log n}{n}$, the above sum is a polynomial in $(\log n)$ of degree $p$.
Thus, since $p$ is fixed, when $n\to +\infty$, there exists $c > 0$ with:
\[\pushQED{\qed} 
\E\bigl(|D f_0(\polytopeQ_n)|^p\,\big|\, \FBody\subseteq \polytopeQ_n\bigr) \leq c\, (\log n)^{p + 1 -\frac{1}{d-1}}\, \left(\frac{1}{n}\right)^{1-\frac{1}{d-1}}.
\qedhere
\]
\end{proof}

Combining \Cref{thm:EfronSteinJackknife,thm:MomentsFirstOrderDifference} with $p=2$ yields the desired upper bound for the variance.
\begin{corollary}\label{cor:VarianceUpper}
Let $d\geq 4$. When $n\to +\infty$, there exists a constant $c > 0$ such that:
$$\Var f_0(\polytopeQ_n) \leq c\, (\log n)^{3 - \frac{1}{d-1}}\, n^{\frac{1}{d-1}}.$$
\end{corollary}

\subsubsection{Central limit theorem}\label{sssec:CLT}

\paragraph{Concentration}
Knowledge about the variances and expectations of a sequence of positive random variables $(Z_n)_{n\in \N}$ yields concentration theorems (\eg via Chebyshev's inequality). 
If the variance is smaller than the square of the expectation, \ie $\Var Z_n = o\bigl((\E(Z_n))^2\bigr)$ if $n\to +\infty$, then $\Var \frac{Z_n}{\E\, Z_n} \to 0$, and $(Z_n)_{n\in \N}$ is ``highly concentrated around its sequence of expectations $(\E\, Z_n)_{n\in\N}$''. 
Applying these ideas to $Z_n = f_0(\polytopeQ_n)$ yields the following concentration inequality.

\begin{corollary}\label{cor:Concentration}
Let $d\geq 4$. Then  $\Var f_0(\polytopeQ_n) = o\Bigl(\bigl(\E f_0(\polytopeQ_n)\bigr)^2\Bigr)$ if  $n\to+\infty$.
In particular, for any fixed $a > 0$, we have
$$\P\left( \left|\frac{f_0(\polytopeQ_n)}{\E f_0(\polytopeQ_n)} - 1\right| ~\geq~ a \right) ~\xrightarrow[n\to +\infty]{}~ 0.$$
\end{corollary}

\begin{proof}
Since $\E f_0(\polytopeQ_n) \sim c\, n^{\frac{1}{d-1}}$, and $\Var f_0(\polytopeQ_n) \leq c'\, (\log n)^{3 - \frac{1}{d-1}}\, n^{\frac{1}{d-1}}$,  \Cref{cor:expectation,cor:VarianceUpper} imply that $\Var f_0(\polytopeQ_n) = o\Bigl(\bigl(\E f_0(\polytopeQ_n)\bigr)^2\Bigr)$.
Applying Chebyshev's inequality for fixed $a > 0$ yields:
$$\P\left( \left|\frac{f_0(\polytopeQ_n)}{\E f_0(\polytopeQ_n)} - 1\right| \geq a \right) ~=~ \P\left( \bigl|f_0(\polytopeQ_n) - \E f_0(\polytopeQ_n)\bigr| \geq a\cdot \E f_0(\polytopeQ_n) \right) ~\leq~ \frac{\Var f_0(\polytopeQ_n)}{a^2 \,\bigl(\E f_0(\polytopeQ_n)\bigr)^2} \to 0.\qedhere$$
\end{proof}


\paragraph{A tool to control the distance to the Gaussian distribution}
In order to prove a central limit theorem, we introduce a powerful tool from \cite{Shao-BerryEsseen}, which simplifies a previously known criterion from \cite{ReyPecati}.
The main idea is to control the Kolmogorov distance between the standard normal distribution and a certain statistic on random polytopes (in our case the number of vertices) through certain quantities, similar to those defined to establish the variance in \Cref{sssec:Variance}. 
Recall that the \defn{Kolmogorov distance} between (real) random variables $X$ and $Y$ is the supremum of the difference of their cumulative distribution functions: \defn{$d_{\text{Kol}} (X, Y)$} $= \sup_{z\in \R} |\P(X \leq z) - \P(Y \leq z)|$.

We now describe the method explicitly.
As, once more, we want to enable the readers to apply the method to their own problems, we introduce the tool from \cite{Shao-BerryEsseen} in the general setting of a \emph{Polish space} $\X$ (\eg a discrete space or a Euclidean space) and a symmetric measurable function on the set of all point configurations on $\X$ that is, $f : \bigcup_{n = 1}^\infty \X^n \to \R$.
For $\b x = (x_1, \dots, x_n)$ and $i< j$, we let\linebreak \defn{$\b x^{ij}$} $= (x_1, \dots, x_{i-1}, x_{i+1}, \dots, x_{j-1}, x_{j+1}, \dots, x_n)$ be the tuple $\b x$ without its $i^{\text{th}}$ and $j^{\text{th}}$~coordinates.

\begin{definition}
The \defn{second order difference operator} of $f : \bigcup_{n=1}^\infty \X^n \to \R$ with respect to $i, j$ is:
$$D_{ij} f(\b x) = f(\b x) - f(\b x^i) - f(\b x^j) + f(\b x^{ij})\quad \text{for } \b x \in \X^n.$$
\end{definition}
Intuitively, $D_if(\b x)$ measures the effect of the removal of the $i$\textsuperscript{th} component of $\b x$ on $f$, while $D_{ij}f(\b x)$ measures not only the effect of the removal of the $i$\textsuperscript{th} and $j$\textsuperscript{th} component of $\b x$ but also their interaction.

\begin{definition}\label{def:recombi}
Let $\b X = (X_1, \dots, X_n)$, $\b X' = (X'_1, \dots, X'_n)$ and $\widetilde{\b X}=(\widetilde{X}_1,\ldots,\widetilde{X}_n)$ be vectors of i.i.d. random variables, taking values in $\X^n$.
A \defn{recombination} of $\{\b X, \b X',\widetilde{\b X}\}$ is a vector of random variables $\b Z = (Z_1, \dots, Z_n)$, where $Z_i \in \{X_i, X'_i, \widetilde{X}_i\}$ for $i\in[n]$.
\end{definition}

For a symmetric measurable function $f:\bigcup_{n=1}^\infty\X^n \to \R$ and vectors of i.i.d. random variables $\b X, \b X'$ and $\widetilde{\b X}$ taking values in $\X^n$, let \defn{$\gamma_1$}, \defn{$\gamma_2$}, \defn{$\gamma_3$} be defined as follows:
\begin{align*}
\gamma_1(f) &\eqdef \E\Bigl(|D f(\b X)|^4\Bigr), \\
\gamma_2(f) &\eqdef \sup_{(\b Y,\b Z)}\E\Bigl(\b 1\bigl(D_{12}f(\b Y)\neq 0\bigr)\, D_1 f(\b Z)^4\Bigr), \\
\gamma_3(f) &\eqdef \sup_{(\b Y,\b Y',\b Z)}\E\Bigl(\b 1\bigl(D_{12}f(\b Y)\neq 0\bigr)\, \b 1\bigl(D_{13}f(\b Y')\neq 0\bigr)\, D_2f(\b Z)^4\Bigr),
\end{align*}
where 
the suprema run over all $(\b Y,\b Z)$ resp. $(\b Y,\b Y',\b Z)$ that are recombinations of $\{\b X,\b X',\widetilde{\b X}\}$.

With these definitions at hand, we can finally explicitly state the main tool.

\begin{theorem}[{\cite{ReyPecati}, \cite[Cor. 2.7]{Shao-BerryEsseen}}]\label{thm:ToolCLT}
Let $X_1,\ldots,X_n$ be i.i.d. random variables, taking values in a Polish space $\X$.
For a symmetric measurable function $f:\bigcup_{n=1}^\infty \X^n\to \R$, let $W = f(X_1,\ldots,X_n)$  with $\E(W) < \infty$ and $\E(W^2) < \infty$.

Let $U \sim \c N(0,1)$ be a standard Gaussian random variable.
Then there exists $c > 0$ such that:
\begin{equation}\label{eq:Tool}
d_{\mathrm{Kol}}\left(\frac{W - \E(W)}{\sqrt{\Var W}} \,,\, U\right) ~\leq~ c\,\, \frac{1}{\Var W } \left(\sqrt{n}\gamma_1(f) \,+\, n \sqrt{\gamma_2(f)} \,+\, n\sqrt{n}\sqrt{\gamma_3(f)}\right).
\end{equation}
\end{theorem}

We apply \Cref{thm:ToolCLT} to $f = f_0(\conv(\cdot))$, showing that the right-hand side of \eqref{eq:Tool} tends to $0$. 


\paragraph{Controlling $\Var$, $\gamma_1$, $\gamma_2$ and $\gamma_3$ for $f_0(\polytopeQ_n)$}

Firstly, $\Var f_0(\polytopeQ_n)$ and $\gamma_1(f_0(\polytopeQ_n))$ have been addressed in \Cref{thm:MomentsFirstOrderDifference,cor:VarianceLower} (with $p = 4$): for all $c_0 > 0$, there exists $c, c' > 0$ such that:
$$\Var\bigl(f_0(\polytopeQ_n)\bigr) \geq c\, n^{\frac{1}{d-1} - c_0} ~~\text{ and }~~ \sqrt{n}\gamma_1\bigl(f_0(\polytopeQ_n)\bigr) \leq c'\, (\log n)^{5 - \frac{1}{d-1}} \left(\frac{1}{n}\right)^{\frac{1}{2} - \frac{1}{d-1}}.$$

To control $\gamma_2$ and $\gamma_3$, we need to understand the interaction of two points in our configuration.

\begin{lemma}\label{lem:DoubleVisibility}
Let $d\geq 4$. For $\varepsilon > 0$, and $X, Y\in \B$, there exists $c > 0$ such that
\[\P(\vis X \cap \vis Y \ne \emptyset) \leq c\, \varepsilon.\]
\end{lemma}

\begin{proof}
Let $P_r$ and $P_\ell$ be the rightmost and leftmost points of the arc $(\vis X) \cap \mathbb{S}^1$. 
If $\vis X\cap \vis Y\ne \emptyset$, then, in particular, $Y\in \bigl(\vis X\cup\vis P_r\cup \vis P_\ell\bigr)$.
With \Cref{lem:VisibilityProba} we conclude that there exists $c > 0$ with $\P(\vis X\cap \vis Y\ne \emptyset) \leq 3\, c\, \varepsilon$.
\end{proof}

\begin{remark}
Since  $\vis X\cap \vis Y\neq \emptyset$ for any $Y\in \vis X$ and  since $\mu(\vis X) \geq c'\, \varepsilon$ by  \Cref{lem:VisibilityProba}, we also have  $\P(\vis X \cap \vis Y \ne \emptyset) \geq c'\, \varepsilon$ for some $c' > 0$.
\end{remark}

\begin{proposition}\label{prop:gamma2}
Let $d\geq 4$. 
There exists $c > 0$ such that, when $n \to +\infty$, we have:
$$\gamma_2\bigl(f_0(\polytopeQ_n)\bigr) \leq c\, (\log n)^{6 - \frac{1}{d-1}}\, \left(\frac{1}{n}\right)^{2 - \frac{1}{d-1}}.$$
\end{proposition}

\begin{proof}
Throughout this proof, let $\varepsilon = c_0 \frac{\log n}{n}$ for some $c_0 > 0$ and let $\b X = (X_1, \dots, X_n)$,\linebreak $\b X' = (X'_1, \dots, X'_n)$ and $\widetilde{\b X}=(\widetilde{X}_1,\ldots,\widetilde{X}_n)$ be random vectors as in \Cref{def:recombi}, picked according to the $\beta$-distribution on $\mathbb{B}^2$. 
To simplify notation, we set $f(\b X) \eqdef f_0(\conv(\b X))$. 

Let $\b Y$ and $\b Z$ be recombinations of $\b X,\b X'$ and $\widetilde{\b X}$. 
In the following, we condition on the floating body being contained in the polygon at stake. 
More precisely, let $A$ be the event that $\FBody\subseteq \bigcap_{\b W\in\{\b Y, \b Z\}}\conv(W_3, \dots, W_n)$.
Using \Cref{prop:FloatingBody} and the union bound, we conclude that for any $s > 0$, there exists $c' > 0$ such that $1 - \P(A) \leq c'\, (n-2)^{-s}\leq c'n^{-s}$.
On the complement of the event $A$, we use the trivial bound $D f(\b Z)^4 \leq n^4$.
On the event $A$, we argue as follows: 
If $\vis Y_1 \cap \vis Y_2 = \emptyset$, then $\conv(\b Y) \ssm \conv(\b Y^1)$ and $\conv(\b Y) \ssm \conv(\b Y^2)$ do not intersect, and hence $D_{12} f(\b Y) = 0$.
 \Cref{lem:DoubleVisibility} thus implies that $\P(D_{12} f(\b Y) \ne 0) \leq c\, \varepsilon$, for some $c > 0$.
Putting these arguments together, we obtain
$$\begin{array}{rcll}
\gamma_2(f) &\leq& \P(A)\, \E\Bigl(\b 1\bigl(D_{12}f(\b Y)\neq 0\bigr)\, D_2f(\b Z)^4 \,\,\Big|\,\, A\Bigr) &+~ \bigl(1-\P(A)\bigr) \, n^4 \\
&\leq& 1\cdot \P\bigl(\vis Y_1 \cap \vis Y_2 \ne \emptyset\bigr) \,\E\bigl(|D f(\b Z)|^4\bigr) &+~ c'\, n^{4-s} \\
&\leq& c\, \varepsilon\, \E\bigl(|D f(\b Z)|^4\bigr) &+~ c'\, n^{4-s}.
\end{array}$$

Finally, using \Cref{thm:MomentsFirstOrderDifference} (with $p = 4$) and $\varepsilon = c_0 \frac{\log n}{n}$, we obtain the claimed upper bound.
\end{proof}

\begin{proposition}\label{prop:gamma3}
Let $d\geq 4$. 
There exists $c > 0$ such that, when $n \to +\infty$, we have:
$$\gamma_3\bigl(f_0(\polytopeQ_n)\bigr) \leq c\, (\log n)^{7 - \frac{1}{d-1}}\, \left(\frac{1}{n}\right)^{3 - \frac{1}{d-1}}.$$
\end{proposition}

\begin{proof}
The proof is very similar to the one of \Cref{prop:gamma2} (up to multiplying once by $\varepsilon$), and we will use the same notation as therein. 
In particular, let $\b Y$, $\b Y'$, $\b Z$ be recombinations of $\b X,\b X'$ and $\widetilde{\b X}$, and let $\varepsilon = c_0 \frac{\log n}{n}$ for some $c_0 > 0$.
Here, we consider the event $A$ of the floating body $\FBody$ being contained in $\bigcap_{\b W\in\{\b Y, \b Y', \b Z\}}\conv(W_4, \dots, W_n)$.
Again by \Cref{prop:FloatingBody}, for any $s > 0$, there exists $c' > 0$ such that: $1 - \P(A) \leq c'\, n^{-s}$. 


As before, if $\vis Y_1 \cap \vis Y_2 = \emptyset$, then $D_{12} f(\b Y) = 0$,
and similarly for $\b Y'$.
So, if $D_{12} f(\b Y) \ne 0$ and $D_{12} f(\b Y') \ne 0$, then $\vis Y_1 \cap \vis Y_2 \ne \emptyset$ and $\vis Y'_1 \cap \vis Y'_2 \ne \emptyset$.
As these two events are independent, by \Cref{lem:DoubleVisibility}, they occur jointly with probability less than $c\, \varepsilon^2$ for $c > 0$.
Thus:
$$\begin{array}{rcll}
\gamma_3(f) &\leq& \P(A)\, \E\Bigl(\b 1\bigl(D_{12}f(\b Y)\neq 0\bigr)\, \b 1\bigl(D_{13}f(\b Y')\neq 0\bigr)\, D_2f(\b Z)^4 \,\,\Big|\,\, A\Bigr) &+~ \bigl(1-\P(A)\bigr) \, n^4 \\
&\leq& 1\cdot \P\bigl(\vis Y_1 \cap \vis Y_2 \ne \emptyset\bigr) \, \P\bigl(\vis Y'_1 \cap \vis Y'_2 \ne \emptyset\bigr) \,\E\bigl(|D f(\b Z)|^4\bigr) &+~ c'\, n^{4-s} \\
&\leq& c\, \varepsilon^2\, \E\bigl(|D f(\b Z)|^4\bigr) &+~ c'\, n^{4-s}.
\end{array}$$

Finally, using \Cref{thm:MomentsFirstOrderDifference} (with $p = 4$) and $\varepsilon = c_0 \frac{\log n}{n}$, the claimed upper bound follows.
\end{proof}

\begin{corollary}\label{cor:CLT}
Let $d\geq 4$, $X_1, \dots, X_n$ be $\beta$-distributed points in the disk $\B$ with $\beta = \frac{d}{2} - 2$, and let $\polytopeQ_n = \conv(X_1, \dots, X_n)$.
Let $U \sim \c N(0,1)$ be a standard Gaussian random variable.
Then:
$$d_{\mathrm{Kol}}\left(\frac{f_0(\polytopeQ_n) - \E f_0(\polytopeQ_n)}{\sqrt{\Var f_0(\polytopeQ_n)}} \,,\, U\right) \xrightarrow[n \to +\infty]{} 0.$$
\end{corollary}

\begin{proof}
By \Cref{thm:MomentsFirstOrderDifference,cor:VarianceLower,prop:gamma2,prop:gamma3}, we have that for $f_0(\polytopeQ_n)$, the quantities $\frac{1}{\Var}$, $\sqrt{n}\gamma_1$, $n\gamma_2$ and $n\sqrt{n}\gamma_3$ can all be upper bounded by terms of the form $(\log n)^a \left(\frac{1}{n}\right)^b$ with $b$ as follows (recall that, in \Cref{cor:VarianceLower}, there is no restriction on $c_0 > 0$):
\begin{center}
\begin{tabular}{c||c|c|c|c}
quantity & $\frac{1}{\Var(\cdot)}$ & $\sqrt{n}\gamma_1$ & $n\sqrt{\gamma_2}$ & $n\sqrt{n}\sqrt{\gamma_3}$ \\
\hline
$b$ & $\frac{1}{d-1} - c_0$ for any $c_0 > 0$ & $\frac{1}{2} - \frac{1}{d-1}$ & $\frac{-1}{2(d-1)}$ & $\frac{-1}{2(d-1)}$
\end{tabular}.
\end{center}

Hence the sum $\frac{1}{\Var} \left(\sqrt{n}\gamma_1 \,+\, n \sqrt{\gamma_2} \,+\, n\sqrt{n}\sqrt{\gamma_3}\right)$ is upper bounded by $c\, (\log n)^a\, \left(\frac{1}{n}\right)^{\frac{1}{2(d-1)}-c_0}$ for some $c > 0$, any $c_0 > 0$, and $a = \frac{7}{2} - \frac{1}{2(d-1)}$.
By \Cref{thm:ToolCLT}, choosing $c_0$ sufficiently small guarantees that the Kolmogorov distance at stake tends to $0$.
\end{proof}

\begin{remark}
This corollary holds for any $\beta$-distribution in the plane, with $\beta > 0$.
\end{remark}

We end this section by providing the proof of \Cref{thm:CLTmain}


\begin{proof}[Proof of \Cref{thm:CLTmain}]
By \Cref{cor:LnIsF0Qn}, the random variables $L(\omega,\polytope{P}_n)$, $f_0^{\text{up}}(\polytopeQ_n)$ and $f_0^{\text{low}}(\polytopeQ_n)$ have the same distribution.
Since $f_0(\polytopeQ_n)= f_0^{\text{up}}(\polytopeQ_n) + f_0^{\text{low}}(\polytopeQ_n)$, we get: $\E L(\omega,\pol_n)=\frac{1}{2} \E f_0(\polytopeQ_n)$.
We need to control the variance of $L(\omega,\pol_n)$. 
We are going to prove that $\Var L(\omega,\pol_n)\sim \frac{1}{2} \Var f_0(\polytopeQ_n)$ by showing that $f_0^{\text{up}}(\polytopeQ_n)$ and $f_0^{\text{low}}(\polytopeQ_n)$ are ``almost independent''.
Firstly:

$$\Var f_0(\polytopeQ_n) = \Var f_0^{\text{up}}(\polytopeQ_n) \,+\, \Var f_0^{\text{low}}(\polytopeQ_n) \,+\, \mathrm{Cov}\Bigl(f_0^{\text{up}}(\polytopeQ_n),\, f_0^{\text{low}}(\polytopeQ_n)\Bigr). $$

For $1\leq k\leq n$, let $V_j^{\text{up}}$ (resp. $V_j^{\text{low}}$) be the event that the point $X_j$ is an upper (resp. lower) vertex (\ie it has an outer normal vector with positive, resp. negative, second coordinate). 
Then:
$$\mathrm{Cov}\Bigl(f_0^{\text{up}}(\polytopeQ_n),\, f_0^{\text{low}}(\polytopeQ_n)\Bigr) = \sum_{j, k = 1}^n\mathrm{Cov}\Bigl(\b 1(V_j^{\text{up}}),\,\b 1(V_k^{\text{low}})\Bigr).$$
Once more, we control the vertices using the floating body.

Suppose $\FBody\subseteq\polytopeQ_n$.
If $X_j$ and $X_k$ (for $k\ne j$) are not contained in a common $\varepsilon$-cap, then the events ``$X_j$ is a vertex of $\polytopeQ_n$'' and ``$X_k$ is a vertex of $\polytopeQ_n$'' are independent (because\linebreak $X_k\in \conv(X_i ~;~ i\ne k)$ is equivalent to $X_k\in \conv(X_i ~;~ i \ne k \text{ and } i\ne j)$).
Thus, if both events $V^{\text{up}}_j$ and $V^{\text{low}}_k$ occur, then $X_j$ and $X_k$ lie in an $\varepsilon$-cap which contains upper and lower vertices, in particular, this $\varepsilon$-cap is contained in the visibility region of $(-1, 0)$ or of $(1, 0)$.

Finally, we use that the covariance of $\bigl(f_0^{\text{up}}(\polytopeQ_n),\, f_0^{\text{low}}(\polytopeQ_n)\bigr)$, conditioned on $\FBody\not\subseteq \polytopeQ_n$, is smaller than $\binom{n}{2}$.
Using $\varepsilon = c_0\frac{\log n}{n}$, \Cref{lem:VisibilityProba,prop:FloatingBody}, there exists $c > 0$ satisfying:
$$\left|\mathrm{Cov}\Bigl(f_0^{\text{up}}(\polytopeQ_n),\, f_0^{\text{low}}(\polytopeQ_n)\Bigr)\right| \,\leq\, \binom{n}{2} \P(\FBody\not\subseteq\polytopeQ_n) \,+\, \binom{n}{2} \Bigl(2\mu\bigl(\vis (-1, 0)\bigr)\Bigr)^2 \P(\FBody\subseteq\polytopeQ_n) \,\leq\, c\, (\log n)^2.$$

As $(\log n)^2 = o\bigl(\Var f_0(\polytopeQ_n)\bigr)$, and $\Var f_0^{\text{up}}(\polytopeQ_n) = \Var f_0^{\text{low}}(\polytopeQ_n)$, due to their distribution being identical, we get that $\Var L(\b\omega, \pol_n) = \Var f_0^{\text{up}}(\polytopeQ_n) \sim \frac{1}{2}\Var f_0(\polytopeQ_n)$.

It remains to show that $f_0^{\text{up}}(\polytopeQ_n)$ admits a central limit theorem.
Unfortunately, this cannot be deduced directly from the central limit theorem for $f_0^{\text{up}}(\polytopeQ_n)$.
However, the attentive reader will easily notice that replacing $f_0$ by $f_0^{\text{up}}$ in \Cref{sssec:CLT} (and \Cref{thm:MomentsFirstOrderDifference}) only changes the constants, not the dependencies in $n$.
Consequently, one can control the Kolmogorov distance between (the normalization of) $f_0^{\text{up}}(\polytopeQ_n)$ and a standard normal distributed variable, finally proving \Cref{thm:CLTmain}.
\end{proof}

\subsection*{Some open problems}

We finish by proposing open problems which naturally extend what we discussed in this section.

\begin{problem}
Compute the higher moments of $L_n$, the length of the coherent path (for fixed $\b c$ and $\b\omega$) on random polytopes.
Equivalently, compute the higher moments of the number of vertices of $\beta$-distributed polygons.
\end{problem}

\begin{problem}
Sample points $X_1, \dots, X_n$ at random, uniformly on the $(d-1)$-sphere, and construct $\pol_n = \conv(X_1, \dots, X_n)$.
For $\b c = \b e_1$, let $N^{\text{coh}}_\ell$ be the number of coherent paths of length $\ell$ on $\pol_n$.
Study the probability that $N^{\text{coh}}_\ell)_\ell$ is unimodal (conjecturally, it tends to $1$ when $n \to+\infty$).
\end{problem}

\begin{problem}
Sample points $X_1, \dots, X_n$ at random, uniformly on the $(d-1)$-sphere, and construct $\pol_n = \conv(X_1, \dots, X_n)$.
Study the number of coherent paths on $\pol_n$ (distribution, expectation, variance, central limit theorem).
\end{problem}

\begin{problem}
Sample points $X_1, \dots, X_n$ at random, uniformly on the $(d-1)$-sphere, and construct $\pol_n = \conv(X_1, \dots, X_n)$.
Study the number and the length of $\b e_1$-monotone paths on $\pol_n$.
\end{problem}

\begin{problem}
Extend the results of this paper to similar probability distributions, and especially to the polar case: random polytopes defined as $\pol^{\circ}_n = \{\b x\in \R^d ~;~ \inner{\b x, \b a_i} \leq 1 \text{ for all } i\in [n]\}$ where the facet normals $\b a_1, \dots, \b a_n$ are chosen uniformly at random on the sphere $\mathbb{S}^{d-1}$.
\end{problem}


\subsection{Cheat sheet of formulas}\label{ssec:CheatSheet}
Here is a quick overview of the main formulas and notations appearing throughout \Cref{sec:RandomCase}.
We use $c$ to denote positive constants (independent of $\beta$, $d$, $n$, $\varepsilon$) which do \textbf{not} need to be equal from one line to another.

\paragraph{Probabilistic model} (\Cref{ssec:ProbabilisticModel}) \vspace{0.1cm}

$X_1, \dots, X_n$, random points in $\B$, i.i.d., $\beta$-distributed, usually for $\beta_d = \frac{1}{2}d - 2$ \vspace{0.1cm}

$f_{2, \beta_d}(\b x) = C_{2, \beta_d} \Bigl( 1 - \|\b x\|\Bigr)^{\beta_d}$ for $\ x\in \B$ \vspace{0.1cm}

$\mu(A)$: measure of $A\subseteq\B$ according to the density function $f_{2, \beta_d}$ \vspace{0.1cm}

$\polytopeQ_n = \conv(X_1, \dots, X_n)$: random $\beta$-polygon \vspace{0.1cm}

$f_0(\polytopeQ_n)$: number of vertices of $\polytopeQ_n$ \vspace{0.1cm}

\paragraph{Expectation, Variance}(\Cref{sssec:Expectation,sssec:Variance})\vspace{0.1cm}

$\E\bigl(f_0(\polytopeQ_n)\bigr) \sim c\, n^{\frac{1}{d-1}}$ \vspace{0.1cm}

$\varepsilon = c_0 \frac{\log n}{n}$: it is a definition \vspace{0.1cm}

$1 - R_\varepsilon \sim c\, \varepsilon^{\frac{2}{d-1}}$: radius of $\varepsilon$-cap and floating body \vspace{0.1cm}

$m_\varepsilon \sim c\, \left(\frac{1}{\varepsilon}\right)^{\frac{1}{d-1}}$: maximal number of independent (\ie disjoint) $\varepsilon$-caps \vspace{0.1cm}

$\P(A_\polytopeC) \geq c\, (\log n)^4 \, n^{-c_0}$: probability of having 4 points ``correctly placed'' in an $\varepsilon$-cap \vspace{0.1cm}

$\P(X \notin \FBody) = \mu(\B \ssm \FBody) \sim c\, \varepsilon^{1 - \frac{1}{d-1}}$: measure of the of the part outside the floating body \vspace{0.1cm}

$\mu\bigl(\vis\b x\bigr) \sim c\, \varepsilon$: measure of the visibility region \vspace{0.1cm}

$\P(\FBody \not\subseteq \polytopeQ_n) \leq c\, n^{-s}$, for any $s > 0$ (requires $c_0 = \frac{1}{d-1} + s$) \vspace{0.1cm}

$\E\bigl(|D f_0(\polytopeQ_n)|^p\bigr) \leq c\, (\log n)^{p + 1 -\frac{1}{d-1}}\, \left(\frac{1}{n}\right)^{1 - \frac{1}{d-1}}$ where $D$ is (any) first order difference operator \vspace{0.1cm}

$c\, n^{\frac{1}{d-1} - c_0} ~\leq~ \Var f_0(\polytopeQ_n) ~\leq~ c'\, (\log n)^{3 - \frac{1}{d-1}}\, n^{\frac{1}{d-1}}$ for any $c_0 > 0$

(the lower bound does not require that $\FBody \subseteq \polytopeQ_n$, \ie $c_0 > 0$ can be arbitrarily small) \vspace{0.1cm}

\paragraph{Central limit theorem}(\Cref{sssec:CLT}) \vspace{0.1cm}

$f(\b X) = f_0(\conv(\b X))$ \vspace{0.1cm}

$\gamma_1(f) = \E\Bigl(|D f(\b X)|^4\Bigr) ~\leq~ c'\, (\log n)^{5 - \frac{1}{d-1}} \left(\frac{1}{n}\right)^{1 - \frac{1}{d-1}}$ \vspace{0.1cm}

$\gamma_2(f) = \sup_{(\b Y,\b Z)}\E\Bigl(\b 1\bigl(D_{12}f(\b Y)\neq 0\bigr)\, D_1 f(\b Z)^4\Bigr) ~\leq~ c\, (\log n)^{6 - \frac{1}{d-1}}\, \left(\frac{1}{n}\right)^{2 - \frac{1}{d-1}}$ \vspace{0.1cm}

$\gamma_3(f) = \sup_{(\b Y,\b Y',\b Z)}\E\Bigl(\b 1\bigl(D_{12}f(\b Y)\neq 0\bigr)\, \b 1\bigl(D_{13}f(\b Y')\neq 0\bigr)\, D_2f(\b Z)^4\Bigr) ~\leq~ c\, (\log n)^{7 - \frac{1}{d-1}}\, \left(\frac{1}{n}\right)^{3 - \frac{1}{d-1}}$ \vspace{0.1cm}

$d_{\text{Kol}}\left(\frac{f_0(\polytopeQ_n) - \E f_0(\polytopeQ_n)}{\sqrt{\Var f_0(\polytopeQ_n)}} \,,\, U\right) ~\leq~ c\, (\log n)^{\frac{7}{2} - \frac{1}{2(d-1)}}\, \left(\frac{1}{n}\right)^{\frac{1}{2(d-1)}} ~\xrightarrow[n\to+\infty]{} 0$, where $U \sim \c N(0, 1)$

\appendix
\section{Algorithms for monotone paths and coherent paths}\label{sec:Algos}


The present section recalls the algorithms used for computing (coherent) monotone paths.
We do not aim at giving a complete presentation of the notion, but rather at introducing it for the non-expert readers.
(Note that this section is not included in the version of this article published in the \emph{International Mathematics Research Notices}, one can only cite it from the ArXiv version.)

\paragraph{Monotone paths}

The naive idea for counting $\b c$-monotone paths on a polytope $\pol$ according to their lengths is to enumerate (asking your favorite library of your favorite programming language) all monotone paths of the graph $G_{\pol, \b c}$, and to record their lengths.
Actually, the usual algorithm to enumerate monotone paths in a graph can directly be adapted to sort them by length.
The time complexity is $O(f_1)$, where~$f_1$ is the number of edges of $\pol$.

The algorithm works as follows:
First, find a \emph{topological order} on $G_{\pol, \b c}$, \ie label the vertices of $G_{\pol, \b c}$ with $\b v_1, \dots, \b v_n$ such that $i < j$ if $\b v_i \to \b v_j$. Given this order, for any $1\leq i\leq n$ and for any outgoing arc $\b v_i\to \b v_j$, the number of monotone paths of length $\ell$ using the edge $\b v_i\to\b v_j$ equals the number of paths of length $\ell - 1$ from $\b v_1$ to $\b v_i$.
This algorithm can be formalized as follows:

\RestyleAlgo{ruled}
\SetKwComment{Comment}{/* }{ */}
\begin{algorithm}[hbt!]
\caption{Counting $\b c$-monotone paths on $\pol$ according to length}

\KwData{Polytope $\pol\subseteq\R^d$ with $n$ vertices, direction $\b c\in \R^d$}
\KwResult{List $L$ such that $L(\ell)$ is the number of $\b c$-monotone paths of length $\ell$ on $\pol$}

\vspace{0.1cm}

$G_{\pol, \b c} \gets \bigl($directed graph of $\pol$ where $\b u\to \b v$ iff $\inner{\b u, \b c} < \inner{\b v, \b c}\bigr)$\;

Relabel the nodes of $G_{\pol, \b c}$ according to a topological order from $\b v_1$ to $\b v_n$\;

\For{$i$ from $1$ to $n$}{
$L_i \gets [0, \dots, 0]$ \Comment*[r]{the list $L_i$ is indexed from $0$ to $n$}
}


$L_{1}(0) \gets 1$\;

\vspace{0.1cm}

\For{$i$ from $1$ to $n-1$}{
\For{$j$ such that there is an out-going edge $\b v_i \to \b v_j \in G_{\pol, \b c}$}{
$L_{j}(\ell) \gets L_{j}(\ell) + L_{i}(\ell-1)$ ~~~~for $\ell\in [1, n]$
}}

\vspace{0.1cm}

\KwOut{$L_n$}
\end{algorithm}

\paragraph{Coherent paths}
There are two usual ways to check if a monotone path is coherent.
Recall that a $\b c$-monotone path on $\pol$ is coherent if and only if there is $\b \omega$ to capture it, \ie
if the upper faces of the projection of $\pol$ onto the plane spanned by $\b c$ and $\b\omega$ is (the projection of) the monotone path.

\vspace{0.27cm}

\noindent
\begin{minipage}[c]{0.46\linewidth}
\underline{First method:}
Let $(\b u_1, \dots, \b u_r)$ be the vertices of the monotone path to be tested for coherence.
Being the upper faces of a\linebreak 2-dimensional polytope $\pol_{\b c, \b \omega}$ is easy to check inductively: suppose we know that $\b u_1, \dots, \b u_i$ are upper vertices of $\pol_{\b c, \b \omega}$ for some given $i < r$, then the next upper vertex is the neighbor $\b v$ of $\b u_i$ which maximizes the slope in the plane $\b c$ and $\b \omega$, \ie we need to find $\b u_{i+1}$ that maximizes $\frac{\inner{\b \omega, \b v - \b u_i}}{\inner{\b c, \b v - \b u_i}}$ under the conditions that $\b u_i\b v$ is an edge of $\pol$ and $\inner{\b u_i, \b c} < \inner{\b v, \b c}$.
This yields several inequalities (maximizing the slope amounts to be greater than all other slopes) which are linear in $\b\omega$.
We can hence gather all these linear inequalities to make a cone and the path is coherent if and only if this cone is full-dimensional.
\end{minipage} 
\hspace{0.03\linewidth}
\begin{minipage}[c]{0.51\linewidth}
\begin{center}
\includegraphics[width=0.95\linewidth]{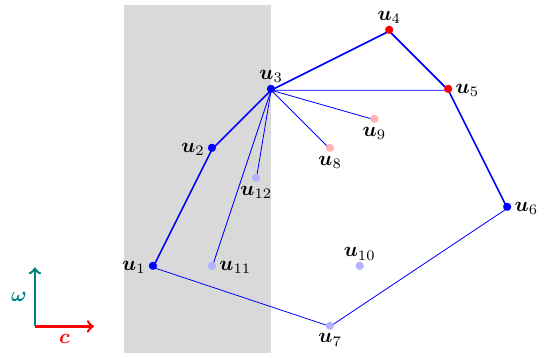}
\addcontentsline{lof}{figure}{13~~~ Finding the upper path of a polygonal projection}
\end{center}
Figure 13:~ Projection of $\pol$ onto a plane (not all edges are drawn).
If $\b u_1, \b u_2, \b u_3$ are upper vertices, then to compute the next vertex, one lists the right-neighbors of $\b u_3$ (not in gray: $\b u_4, \b u_5, \b u_8, \b u_9$), and finds the slope-maximizer from $\b u_3$: here~$\b u_4$.
\label{fig:Algo2}
\end{minipage}
\vspace{0.27cm}

The algorithm is easy to write.
It creates a cone in dimension $d$ with $O(f_1)$ inequalities, where~$f_1$ is the number of edges of $\pol$.
As far as we know, even if improvements exist, the time complexity remains bound to finding a vector in the interior of such a cone, hence it is $O(d^2\, f_1)$.

\underline{Second method:}
Let $\pol[M]_{\pol, \b c}$ be the monotone path polytope of polytope~$\pol$ and direction $\b c$, see \cite{BilleraSturmfels-FiberPolytope,Athanasiadis_2000,BlackDeLoera2021monotone} for definitions, or \cite[Section 2.1]{poullot2024verticesmonotonepathpolytopes} for constructions of monotone path polytopes.
The idea is that for $\b \omega, \b \omega'\in \R^d$, the vertices of $\pol[M]_{\pol, \b c}^{\b\omega}$ and $\pol[M]_{\pol, \b c}^{\b\omega'}$ are the same if and only if the coherent paths captured by $\b\omega$ and $\b\omega'$ are the same.
The next algorithm uses this key concept.

The time complexity is driven by the computation of the monotone path polytope.
One can compute it via $n-2$ Minkowski sums in $\R^d$ (for $d$-polytope with $n$ vertices).
Even though costly, most polytope libraries carry efficiently implemented algorithms which are enough for our use-case.


\begin{algorithm}[hbt!]
\caption{Checking if a monotone path is coherent}

\KwData{Polytope $\pol\subseteq\R^d$, direction $\b c\in \R^d$, and $\b c$-monotone path $(\b u_1, \dots, \b u_r)$ on $\pol$}
\KwResult{Boolean answer to ``$(\b u_1, \dots, \b u_r)$ is coherent?'', and certificate $\b\omega$ if true}

\vspace{0.1cm}

$\pol[C] \gets \cone\Bigl\{\b \omega\in \R^d ~;~ \forall i\in[r-1],~\forall \b v \text{ improving neighbor of } \b u_i,~ \frac{\inner{\b \omega, \b u_{i+1} - \b u_i}}{\inner{\b c, \b u_{i+1} - \b u_i}} \geq \frac{\inner{\b \omega, \b v - \b u_i}}{\inner{\b c, \b v - \b u_i}}\Bigr\}$

\vspace{0.1cm}

\KwOut{Boolean: $\dim\pol[C] = d$; certificate: $\b\omega \in \text{interior}(\pol[C])$}
\end{algorithm}


\begin{algorithm}[h!]
\caption{Finding all the coherent paths via the monotone path polytope}

\KwData{Polytope $\pol\subseteq\R^d$, direction $\b c\in \R^d$}
\KwResult{List of all $\b c$-coherent paths on $\pol$ (with certificates)}


$\pol[M]_{\pol, \b c}$, \texttt{Paths}, \texttt{Certificates} $\gets$ monotone path polytope of $\pol$ for the direction $\b c$, $[\,]$, $[\,]$

\For{$\b v$ a vertex of $\pol[M]_{\pol, \b c}$}{
$\b\omega \gets$ a vector in the interior of the normal cone of $\b v$

\texttt{Paths} is appended the list of upper vertices of the projection of $\pol$ in the plane $(\b c, \b \omega)$

\texttt{Certificates} is appended $\b\omega$
}


\KwOut{\texttt{Paths}, \texttt{Certificates}}
\end{algorithm}













\bibliographystyle{alpha}
\bibliography{Biblio}

@article{BilleraSturmfels-FiberPolytope,
	author = {Billera, Louis J. and Sturmfels, Bernd},
	journal = {Anals of Mathematics},
	number = {135},
	pages = {527--549},
	publisher = {Springer},
	title = {Fiber Polytopes},
	year = {1992}
}

@article{Athanasiadis_2000,
	author = {Christos A. Athanasiadis and Jes{\'{u}}s A. De Loera and Victor Reiner and Francisco Santos},
	doi = {10.1006/eujc.1999.0319},
	journal = {European Journal of Combinatorics},
	month = {jan},
	number = {1},
	pages = {19--47},
	publisher = {Elsevier {BV}},
	title = {Fiber Polytopes for the Projections between Cyclic Polytopes},
	volume = {21},
	year = 2000
}

@article{Black-PhD,
    author = {Alexander Black},
    title = {Monotone Paths on Polytopes: Combinatorics and Optimization},
    journal = {PhD thesis},
    year = {2024}
}

@article{BlackDeLoera2021monotone,
      title={Monotone Paths on Cross-Polytopes}, 
      author={Alexander Black and Jesús A. De Loera},
      year={2023},
      eprint={2102.01237},
      archivePrefix={arXiv},
      primaryClass={math.CO},
    journal={Discrete \& Computational Geometry},
pages={1245 -- 1265},
issue=4,
volume={70},
doi={10.1007/s00454-023-00563-4}
}

@article{Braenden,
    author = {Petter Br{\"a}nden},
    title = {Unimodality, log-concavity, real-rootedness and beyond},
    journal = {Handbook of Combinatorics},
    year = {2015},
pages={437--483}
}

@article{EfronStein,
author = {Bradley Efron and Charles M. Stein},
title = {{The Jackknife Estimate of Variance}},
volume = {9},
journal = {The Annals of Statistics},
number = {3},
publisher = {Institute of Mathematical Statistics},
pages = {586 -- 596},
year = {1981},
doi = {10.1214/aos/1176345462},
URL = {https://doi.org/10.1214/aos/1176345462}
}

@misc{ManeckeSanyalSo2019shypersimplices,
      title={S-hypersimplices, pulling triangulations, and monotone paths}, 
      author={Sebastian Manecke and Raman Sanyal and Jeonghoon So},
      year={2020},
      journal={The electronic journal of combinatorics},
      volume={27},
      number={3},
}

@article{KabluchkoThaeleZaporozhets:Beta,
title = {Beta polytopes and Poisson polyhedra: $f$-vectors and angles},
journal = {Advances in Mathematics},
volume = {374},
pages = {107333},
year = {2020},
issn = {0001-8708},
doi = {https://doi.org/10.1016/j.aim.2020.107333},
url = {https://www.sciencedirect.com/science/article/pii/S0001870820303613},
author = {Zakhar Kabluchko and Christoph Th{\"a}le and Dmitry Zaporozhets}
}

@article{KabluchkoTemesvariThale-IntrinsicVolumeBetaPolytopes,
    AUTHOR = {Kabluchko, Zakhar and Temesvari, Daniel and Th\"ale,
              Christoph},
     TITLE = {Expected intrinsic volumes and facet numbers of random
              beta-polytopes},
   JOURNAL = {Math. Nachr.},
  FJOURNAL = {Mathematische Nachrichten},
    VOLUME = {292},
      YEAR = {2019},
    NUMBER = {1},
     PAGES = {79--105},
      ISSN = {0025-584X,1522-2616},
   MRCLASS = {52A22 (52B11 53C65 60D05)},
  MRNUMBER = {3909222},
MRREVIEWER = {Jan\ Rataj},
       DOI = {10.1002/mana.201700255},
       URL = {https://doi.org/10.1002/mana.201700255},
}

@article{Reitzner:EfronStein,
author = {Matthias Reitzner},
title = {{Random polytopes and the Efron--Stein jackknife inequality}},
volume = {31},
journal = {The Annals of Probability},
number = {4},
publisher = {Institute of Mathematical Statistics},
pages = {2136 -- 2166},
year = {2003},
doi = {10.1214/aop/1068646381},
URL = {https://doi.org/10.1214/aop/1068646381}
}

@misc{poullot2024verticesmonotonepathpolytopes,
      title={Vertices of the monotone path polytopes of hypersimplicies}, 
      author={Germain Poullot},
      year={2024},
      eprint={2411.14102},
      archivePrefix={arXiv},
      primaryClass={math.CO},
      url={https://arxiv.org/abs/2411.14102}, 
}

@book {AE-Landau,
    AUTHOR = {Amann, Herbert and Escher, Joachim},
     TITLE = {Analysis. {I}},
      NOTE = {Translated from the 1998 German original by Gary Brookfield},
 PUBLISHER = {Birkh\"auser Verlag, Basel},
      YEAR = {2005},
     PAGES = {xiv+426},
      ISBN = {3-7643-7153-6},
   MRCLASS = {00-01 (26-01)},
  MRNUMBER = {2105047},
       DOI = {10.1007/b137107},
       URL = {https://doi.org/10.1007/b137107},
}

@article{Loday2004-Associahedron,
    author = {Loday, {J}ean-{L}ouis},
    title = {Realization of the {S}tasheff polytope},
    journal = {Archiv der Mathematik},
    year = 2004,
    doi = {10.1007/s00013-004-1026-y},
    volume = {83},
    number = {3},
    pages = {267-278}
}

@article{PilaudSantosZiegler-CelebratingLodayAssociahedron,
    author = {Pilaud, Vincent and Santos, Francisco and Ziegler, Günter M.},
    title = {Celebrating {L}oday’s associahedron},
    journal = {Archiv der Mathematik},
    year = 2023,
    doi = {10.1007/s00013-023-01895-6},
    volume = {121},
    number = {5},
    pages = {559-601}
}

@article{Nelson2017-MaximalChainsTamariLattice,
title = {A recursion on maximal chains in the Tamari lattices},
journal = {Discrete Mathematics},
volume = {340},
number = {4},
pages = {661-677},
year = {2017},
issn = {0012-365X},
doi = {https://doi.org/10.1016/j.disc.2016.11.030},
url = {https://www.sciencedirect.com/science/article/pii/S0012365X16303958},
author = {Luke Nelson},
}

@misc{DahlbergFishel2024-MaximalChainsGraphAssociahedra,
      title={Maximal chains in lattices from graph associahedra: Tamari to the weak order}, 
      author={Samantha Dahlberg and Susanna Fishel},
      year={2024},
      eprint={2409.13898},
      archivePrefix={arXiv},
      primaryClass={math.CO},
      url={https://arxiv.org/abs/2409.13898}, 
}

@article{Shao-BerryEsseen,
title = {Berry-Esseen bounds for functionals of independent random variables},
journal = {Stochastic Processes and their Applications},
volume = {183},
pages = {104574},
year = {2025},
issn = {0304-4149},
doi = {https://doi.org/10.1016/j.spa.2025.104574},
url = {https://www.sciencedirect.com/science/article/pii/S0304414925000158},
author = {Qi-Man Shao and Zhuo-Song Zhang}
}

@article{ReyPecati,
author = {Rapha{\"e}l Lachi{\`e}ze-Rey and Giovanni Peccati},
title = {{New Berry–Esseen bounds for functionals of binomial point processes}},
volume = {27},
journal = {The Annals of Applied Probability},
number = {4},
publisher = {Institute of Mathematical Statistics},
pages = {1992 -- 2031},
year = {2017},
doi = {10.1214/16-AAP1218},
URL = {https://doi.org/10.1214/16-AAP1218}
}

@article{Reitzner-CLT,
    author ={Reitzner, Matthias},
title={Central limit theorems for random polytopes},
journal= {Probability Theory and Related Fields},
year={2005},
pages={483--507},
volume={133},
doi={10.1007/s00440-005-0441-8},
}

@article{KellyTolle1981-ExpectedNbVerts,
   AUTHOR = {Kelly, Douglas G. and Tolle, Jon W.},
    TITLE = {Expected number of vertices of a random convex polyhedron},
  JOURNAL = {SIAM J. Algebraic Discrete Methods},
 FJOURNAL = {Society for Industrial and Applied Mathematics. Journal on
             Algebraic and Discrete Methods},
   VOLUME = {2},
     YEAR = {1981},
   NUMBER = {4},
    PAGES = {441--451},
     ISSN = {0196-5212},
  MRCLASS = {52A22 (05A15 52A25 60D05 90C15)},
 MRNUMBER = {634368},
MRREVIEWER = {V.\ K.\ Ohanyan},
      DOI = {10.1137/0602047},
      URL = {https://doi.org/10.1137/0602047},
}

@book{Borgwardt1987-SimplexMethod,
   AUTHOR = {Borgwardt, Karl-Heinz},
    TITLE = {The simplex method: {A} probabilistic analysis},
   SERIES = {Algorithms and Combinatorics: Study \& Research Texts},
   VOLUME = {1},
PUBLISHER = {Springer-Verlag, Berlin},
     YEAR = {1987},
    PAGES = {xii+268},
     ISBN = {3-540-17096-0},
  MRCLASS = {90C05},
 MRNUMBER = {868467},
MRREVIEWER = {J\"urgen\ K\"ohler},
      DOI = {10.1007/978-3-642-61578-8},
      URL = {https://doi.org/10.1007/978-3-642-61578-8},
}

@article{KarpSitnik2010-LogConcavityHypergeometricFunctions,
    AUTHOR = {Karp, Dmitrii B. and Sitnik, Sergei M.},
     TITLE = {Log-convexity and log-concavity of hypergeometric-like functions},
   JOURNAL = {J. Math. Anal. Appl.},
  FJOURNAL = {Journal of Mathematical Analysis and Applications},
    VOLUME = {364},
      YEAR = {2010},
    NUMBER = {2},
     PAGES = {384--394},
      ISSN = {0022-247X,1096-0813},
   MRCLASS = {33C05 (26A51)},
  MRNUMBER = {2576190},
MRREVIEWER = {\'Arp\'ad\ Baricz},
       DOI = {10.1016/j.jmaa.2009.10.057},
       URL = {https://doi.org/10.1016/j.jmaa.2009.10.057},
}

@article{KalmykovKarp2017-LogConcavityAndTuranInequalityHypergeometricFunctions,
    AUTHOR = {Kalmykov, Sergey I. and Karp, Dmitrii B.},
     TITLE = {Log-concavity and {T}ur{\'a}n-type inequalities for the generalized hypergeometric function},
   JOURNAL = {Anal. Math.},
  FJOURNAL = {Analysis Mathematica},
    VOLUME = {43},
      YEAR = {2017},
    NUMBER = {4},
     PAGES = {567--580},
      ISSN = {0133-3852,1588-273X},
   MRCLASS = {33C20 (26A51 26D15)},
  MRNUMBER = {3738361},
MRREVIEWER = {Min-Jie\ Luo},
       DOI = {10.1007/s10476-017-0503-z},
       URL = {https://doi.org/10.1007/s10476-017-0503-z},
}

@incollection{KleeMinty-SimplexAlgorithm,
    AUTHOR = {Klee, Victor and Minty, George J.},
     TITLE = {How good is the simplex algorithm?},
 BOOKTITLE = {Inequalities, {III} ({O}. {S}hisha, editor)},
     PAGES = {159--175},
YEAR = {1972},
 PUBLISHER = {Academic Press, New York-London},
}

@article{BlanchardDeLoeraLouveaux-LengthMonotonePath,
    AUTHOR = {Blanchard, Mo\"ise and De Loera, Jes\'us A. and Louveaux,
              Quentin},
     TITLE = {On the length of monotone paths in polyhedra},
   JOURNAL = {SIAM J. Discrete Math.},
  FJOURNAL = {SIAM Journal on Discrete Mathematics},
    VOLUME = {35},
      YEAR = {2021},
    NUMBER = {3},
     PAGES = {1746--1768},
      ISSN = {0895-4801,1095-7146},
   MRCLASS = {52B12 (05C20 52B55 90C05 90C49)},
  MRNUMBER = {4297821},
MRREVIEWER = {Ruriko\ Yoshida},
       DOI = {10.1137/20M1315646},
       URL = {https://doi.org/10.1137/20M1315646},
}

@article{AthanasiadisDeLoeraZhang-EnumerativeProblems,
    AUTHOR = {Athanasiadis, Christos A. and De Loera, Jes\'us A. and Zhang,
              Zhenyang},
     TITLE = {Enumerative problems for arborescences and monotone paths on
              polytope graphs},
   JOURNAL = {J. Graph Theory},
  FJOURNAL = {Journal of Graph Theory},
    VOLUME = {99},
      YEAR = {2022},
    NUMBER = {1},
     PAGES = {58--81},
      ISSN = {0364-9024,1097-0118},
   MRCLASS = {05C30 (05C10)},
  MRNUMBER = {4371479},
MRREVIEWER = {Robert\ Cimikowski},
       DOI = {10.1002/jgt.22725},
       URL = {https://doi.org/10.1002/jgt.22725},
}

@article{AthanasiadisSantos-MonotonePathsOnZonotopes,
    AUTHOR = {Athanasiadis, Christos A. and Santos, Francisco},
     TITLE = {Monotone paths on zonotopes and oriented matroids},
   JOURNAL = {Canad. J. Math.},
  FJOURNAL = {Canadian Journal of Mathematics. Journal Canadien de
              Math\'ematiques},
    VOLUME = {53},
      YEAR = {2001},
    NUMBER = {6},
     PAGES = {1121--1140},
      ISSN = {0008-414X,1496-4279},
   MRCLASS = {52C40 (20F55 52B05 52B12)},
  MRNUMBER = {1863845},
MRREVIEWER = {Winfried\ Hochst\"attler},
       DOI = {10.4153/CJM-2001-042-3},
       URL = {https://doi.org/10.4153/CJM-2001-042-3},
}

@article{AthanasiadisEdelmanReiner-MonotonePathsOnPolytopes,
    AUTHOR = {Athanasiadis, Christos A. and Edelman, Paul H. and Reiner,
              Victor},
     TITLE = {Monotone paths on polytopes},
   JOURNAL = {Math. Z.},
  FJOURNAL = {Mathematische Zeitschrift},
    VOLUME = {235},
      YEAR = {2000},
    NUMBER = {2},
     PAGES = {315--334},
      ISSN = {0025-5874,1432-1823},
   MRCLASS = {52B11 (05C40 52B22)},
  MRNUMBER = {1795510},
MRREVIEWER = {Tara\ S.\ Holm},
       DOI = {10.1007/s002090000152},
       URL = {https://doi.org/10.1007/s002090000152},
}

@article{AngelHolroydRomikVirag-RandomSortingNetworks,
    AUTHOR = {Angel, Omer and Holroyd, Alexander E. and Romik, Dan and
              Vir\'ag, B\'alint},
     TITLE = {Random sorting networks},
   JOURNAL = {Adv. Math.},
  FJOURNAL = {Advances in Mathematics},
    VOLUME = {215},
      YEAR = {2007},
    NUMBER = {2},
     PAGES = {839--868},
      ISSN = {0001-8708,1090-2082},
   MRCLASS = {60C05 (05C25 05E10 68P10)},
  MRNUMBER = {2355610},
MRREVIEWER = {Yoshiharu\ Kohayakawa},
       DOI = {10.1016/j.aim.2007.05.019},
       URL = {https://doi.org/10.1016/j.aim.2007.05.019},
}

@article{Dauvergne-ArchimedianLimitRandomSortingNetwork,
    AUTHOR = {Dauvergne, Duncan},
     TITLE = {The {A}rchimedean limit of random sorting networks},
   JOURNAL = {J. Amer. Math. Soc.},
  FJOURNAL = {Journal of the American Mathematical Society},
    VOLUME = {35},
      YEAR = {2022},
    NUMBER = {4},
     PAGES = {1215--1267},
      ISSN = {0894-0347,1088-6834},
   MRCLASS = {60C05 (05A05 68P10)},
  MRNUMBER = {4467309},
MRREVIEWER = {Hai\ Yan\ Chen},
       DOI = {10.1090/jams/993},
       URL = {https://doi.org/10.1090/jams/993},
}

@article{BlackSanyal-FlagPolymatroids,
    AUTHOR = {Black, Alexander and Sanyal, Raman},
     TITLE = {Underlying flag polymatroids},
   JOURNAL = {Adv. Math.},
  FJOURNAL = {Advances in Mathematics},
    VOLUME = {453},
      YEAR = {2024},
     PAGES = {Paper No. 109835, 42},
      ISSN = {0001-8708,1090-2082},
   MRCLASS = {90C57 (05A05 05B35 52B12 52B40 90C27)},
  MRNUMBER = {4776323},
MRREVIEWER = {Brahim\ Chaourar},
       DOI = {10.1016/j.aim.2024.109835},
       URL = {https://doi.org/10.1016/j.aim.2024.109835},
}

@book{Dantzig1963-LinearProgrammingAndExtensions,
    AUTHOR = {Dantzig, George B.},
     TITLE = {Linear programming and extensions},
 PUBLISHER = {Princeton University Press, Princeton, NJ},
      YEAR = {1963},
     PAGES = {xvi+625},
   MRCLASS = {90.50},
  MRNUMBER = {201189},
MRREVIEWER = {P.\ L.\ Hammer},
}

@article{MallowsVanderbei-YoungTableauxAsSums,
    AUTHOR = {Mallows, Colin and Vanderbei, Robert J.},
     TITLE = {Which {Y}oung tableaux can represent an outer sum?},
   JOURNAL = {J. Integer Seq.},
  FJOURNAL = {Journal of Integer Sequences},
    VOLUME = {18},
      YEAR = {2015},
    NUMBER = {9},
     PAGES = {Article 15.9.1, 8},
      ISSN = {1530-7638},
   MRCLASS = {05E10},
  MRNUMBER = {3401114},
MRREVIEWER = {Ricardo\ Mamede},
       DOI = {10.1088/2058-7058/18/9/19},
       URL = {https://doi.org/10.1088/2058-7058/18/9/19},
}

@misc{AraujoBlackBurcroffGaoKruegerMcDonough-RealizableStandardYoungTableaux,
      title={Realizable Standard Young Tableaux}, 
      author={Igor Araujo and Alexander Black and Amanda Burcroff and Yibo Gao and Robert A. Krueger and Alex McDonough},
      year={2023},
      eprint={2302.09194},
      archivePrefix={arXiv},
      primaryClass={math.CO},
      url={https://arxiv.org/abs/2302.09194}, 
}

@article{BlackDeLoeraKaferSanita-SimplexMethodOn01polytopes,
  title={On the Simplex Method for 0/1-Polytopes},
  author={Alexander Black and Jes{\'u}s A. De Loera and Sean Kafer and Laura Sanit{\`a}},
  journal={Mathematics of Operations Research},
  year={2021},
}

@article{Thale-CentralLimitTheorem,
    AUTHOR = {Th{\"a}le, Christoph},
     TITLE = {Central limit theorem for the volume of random polytopes with vertices on the boundary},
   JOURNAL = {Discrete Comput. Geom.},
  FJOURNAL = {Discrete \& Computational Geometry. An International Journal of Mathematics and Computer Science},
    VOLUME = {59},
      YEAR = {2018},
    NUMBER = {4},
     PAGES = {990--1000},
      ISSN = {0179-5376,1432-0444},
   MRCLASS = {52A22 (60D05 60F05)},
  MRNUMBER = {3802312},
MRREVIEWER = {Elena\ Villa},
       DOI = {10.1007/s00454-017-9862-2},
       URL = {https://doi.org/10.1007/s00454-017-9862-2},
}

@article{BesauRosenThale-Projective,
    AUTHOR = {Besau, Florian and Rosen, Daniel and Th{\"a}le, Christoph},
     TITLE = {Random inscribed polytopes in projective geometries},
   JOURNAL = {Math. Ann.},
  FJOURNAL = {Mathematische Annalen},
    VOLUME = {381},
      YEAR = {2021},
    NUMBER = {3-4},
     PAGES = {1345--1372},
      ISSN = {0025-5831,1432-1807},
   MRCLASS = {52A22 (52A55 58B20 60D05 60F05)},
  MRNUMBER = {4333417},
MRREVIEWER = {Matthias\ Reitzner},
       DOI = {10.1007/s00208-021-02257-9},
       URL = {https://doi.org/10.1007/s00208-021-02257-9},
}

@article{ThaleTurchiWespi-IntrinsicVolumes,
    AUTHOR = {Th{\"a}le, Christoph and Turchi, Nicola and Wespi, Florian},
     TITLE = {Random polytopes: central limit theorems for intrinsic
              volumes},
   JOURNAL = {Proc. Amer. Math. Soc.},
  FJOURNAL = {Proceedings of the American Mathematical Society},
    VOLUME = {146},
      YEAR = {2018},
    NUMBER = {7},
     PAGES = {3063--3071},
      ISSN = {0002-9939,1088-6826},
   MRCLASS = {52A22 (60D05 60F05)},
  MRNUMBER = {3787367},
MRREVIEWER = {Ilya\ S.\ Molchanov},
       DOI = {10.1090/proc/14000},
       URL = {https://doi.org/10.1090/proc/14000},
}

@article {BilleraKapranovSturmfels-CellularStrings,
    AUTHOR = {Billera, Louis J. and Kapranov, Mikhail M. and Sturmfels, Bernd},
     TITLE = {Cellular strings on polytopes},
   JOURNAL = {Proc. Amer. Math. Soc.},
  FJOURNAL = {Proceedings of the American Mathematical Society},
    VOLUME = {122},
      YEAR = {1994},
    NUMBER = {2},
     PAGES = {549--555},
      ISSN = {0002-9939,1088-6826},
   MRCLASS = {52B40 (52B99 55P99 57T30)},
  MRNUMBER = {1205482},
MRREVIEWER = {Sergey\ Yuzvinsky},
       DOI = {10.2307/2161048},
       URL = {https://doi.org/10.2307/2161048},
}

@book {Ziegler-LecturesOnPolytopes,
    AUTHOR = {Ziegler, G\"unter M.},
     TITLE = {Lectures on polytopes},
    SERIES = {Graduate Texts in Mathematics},
    VOLUME = {152},
 PUBLISHER = {Springer-Verlag, New York},
      YEAR = {1995},
     PAGES = {x+370},
      ISBN = {0-387-94365-X},
   MRCLASS = {52Bxx},
  MRNUMBER = {1311028},
MRREVIEWER = {Margaret\ M.\ Bayer},
       DOI = {10.1007/978-1-4613-8431-1},
       URL = {https://doi.org/10.1007/978-1-4613-8431-1},
}

@book {Gruenbaum-ConvexPolytopes,
    AUTHOR = {Gr{\"u}nbaum, Branko},
     TITLE = {Convex polytopes},
    SERIES = {Pure and Applied Mathematics},
    VOLUME = {Vol. 16},
      NOTE = {With the cooperation of Victor Klee, M. A. Perles and G. C.
              Shephard},
 PUBLISHER = {Interscience Publishers John Wiley \& Sons, Inc., New York},
      YEAR = {1967},
     PAGES = {xiv+456},
   MRCLASS = {52.10},
  MRNUMBER = {226496},
MRREVIEWER = {G.\ T.\ Sallee},
}
\label{sec:biblio}

\end{document}